\documentclass[smallcondensed,envcountsect]{svjour3}
\smartqed

\usepackage{geometry}                
\usepackage{graphicx}
\usepackage{amssymb}
\usepackage{amsmath}
\usepackage{color}

\usepackage{amsthm}
\usepackage{bm}
\usepackage{epstopdf}
\usepackage{accents}
\usepackage{stmaryrd}
\usepackage{adjustbox}

\definecolor{redcol}{rgb}{1.,0.,0.0} 
\definecolor{lnkcol}{rgb}{0.,0.,0.0} 
\usepackage[pdftex,pdfpagelabels=true,bookmarksdepth=3]{hyperref}
\hypersetup{  colorlinks=true,
              citecolor=lnkcol,%
              linkcolor=lnkcol,%
              urlcolor=lnkcol,
            pdfborder=false,
            bookmarks=true,
            bookmarksnumbered=true,
            pdftitle = {The BR1 Scheme is Stable for the Compressible Navier-Stokes Equations},
pdfsubject = {},
pdfauthor = {G.Gassner,D.Kopriva,F.Hindenlang,A.Winters},
pdfkeywords = {}
 }
\usepackage[all]{hypcap}  
\usepackage[dvips]{thumbpdf}

\usepackage{cleveref} 

  \crefname{chapter}{Chap.}{Chaps.}       
  \Crefname{chapter}{Chapter}{Chapters}
  \crefname{section}{Sec.}{Secs.}
  \Crefname{section}{Section}{Sections}
  \crefname{table}{Table}{Tables}
  \Crefname{table}{Table}{Tables}
  \crefname{figure}{Fig.}{Figs.}
  \Crefname{figure}{Figure}{Figures}
  \crefname{equation}{}{}
  \Crefname{equation}{Equation}{Equations}
  \crefname{algorithm}{Alg.}{Algs.}
  \Crefname{algorithm}{Algorithm}{Algorithms}
  \crefname{thm}{Thm.}{Thms.}
  \Crefname{thm}{Theorem}{Theorems}
  \crefname{lem}{Lemma}{Lemmas}
  \Crefname{lem}{Lemma}{Lemmas}
  \crefname{cor}{Corr.}{Corrs.}
  \Crefname{cor}{Corrollary}{Corrollaries}
  \crefname{prop}{Prop.}{Props.}
  \Crefname{prop}{Proposition}{Propositions}
  \crefname{rem}{Rem.}{Rems.}
  \Crefname{rem}{Remark}{Remarks}
  \crefname{appendix}{Appendix}{Appendices}
  \Crefname{appendix}{Appendix}{Appendices}
  
\usepackage{stackengine}
\stackMath

\usepackage{accsupp} 
\newcommand*{\llbrace}{%
  \BeginAccSupp{method=hex,unicode,ActualText=2983}%
    \textnormal{\usefont{OMS}{lmr}{m}{n}\char102}%
    \hspace*{-3pt}
    \textnormal{\usefont{OMS}{lmr}{m}{n}\char102}%
  \EndAccSupp{}%
}
\newcommand*{\rrbrace}{%
  \BeginAccSupp{method=hex,unicode,ActualText=2984}%
    \textnormal{\usefont{OMS}{lmr}{m}{n}\char103}%
        \hspace*{-3pt}
    \textnormal{\usefont{OMS}{lmr}{m}{n}\char103}%
  \EndAccSupp{}%
}
\makeatletter
\def\munderbar#1{\underline{\sbox\tw@{$#1$}\dp\tw@\z@\box\tw@}}
\makeatother

\newcommand\iprod[1]{\left\langle #1\right\rangle} 				
\newcommand\inorm[1]{\left |\left| #1\right|\right|}		
\newcommand\iprodN[1]{\left\langle #1\right\rangle_{N}}			
\newcommand\spacevec[1]{\accentset{\,\rightarrow}{#1}}					
\newcommand\contravec[1]{\tilde{ #1}}					
\newcommand\contraspacevec[1]{\spacevec{\tilde{#1}}}
\newcommand\statevec[1]{\mathbf #1}					
\newcommand\contrastatevec[1]{\tilde{\mathbf #1}}

\newcommand\bigmatrix[1]{\mathcal #1}

\newcommand{\R}{\ensuremath{\bm{\hat{\mathcal{R}}_1}}}
\newcommand{\T}{\ensuremath{\bm{\hat{\mathcal{T}}}}}

\newcommand{\uavg}{\overline{\lVert \vec{v} \rVert^2}}

\newcommand{\textblue}[1]{\textcolor{black}{#1}}
\newcommand{\textred}[1]{\textcolor{black}{#1}}
\newcommand{\shat}{\ensuremath{\hat{s}}} 

\newcommand\acclrvec[1]{\accentset{\,\leftrightarrow}{#1}}  
\newcommand\bigstatevec[1]{\acclrvec{{\mathbf #1}}}
\newcommand\biggreekstatevec[1]{\acclrvec{\boldsymbol #1}}
\newcommand\bigcontravec[1]{\acclrvec{\tilde{\mathbf #1}}}

\newcommand\overRe{\frac{1}{{\operatorname{Re} }}}
\newcommand\twooverRe{\frac{2}{{\operatorname{Re} }}}
\newcommand\dS{\,\operatorname{dS} }
\newcommand\BI{\,\operatorname{BI} }
\newcommand\BL{\,\operatorname{BL} }
\newcommand\BR{\,\operatorname{BR} }
\newcommand\PBT{\,\operatorname{PBT} }
\newcommand\burg{\mathrm{Burgers}}
\newcommand\stand{\mathrm{Standard}}
\newcommand\ec{\mathrm{ec}}
\newcommand\ent{{\,\epsilon}} 
\newcommand\interiorfaces{{\mathrm{interior}\atop\mathrm{faces}}}
\newcommand\boundaryfaces{{\mathrm{Boundary}\atop\mathrm{faces}}}

\newcommand\mmatrix[1]{\underbar{#1}}				
\newcommand{\jump}[1]{\left\llbracket#1\right\rrbracket}   
\newcommand{\average}[1]{\left\{\!\!\left\{#1\right\}\!\!\right\}}   
\newcommand{\partialderivative}[2]{\frac{\partial #1}{\partial #2}}  
\newcommand{\avg}[1]{\ensuremath{\llbrace#1\rrbrace}}

\theoremstyle{plain}
\theoremstyle{remark}
\newtheorem{rem}{Remark}

\numberwithin{equation}{section}

\begin{document}

\title{The BR1 Scheme is Stable for the Compressible Navier-Stokes Equations}
\author{Gregor~J.~Gassner \and Andrew~R.~Winters \and Florian~J.~Hindenlang \and David~A.~Kopriva}
\institute{Gregor J. Gassner (\email{ggassner@math.uni-koeln.de})
           \and Andrew R. Winters \at Mathematical Institute, University of Cologne, Cologne, Germany\\ 
           Florian J. Hindenlang \at Max Planck Institute for Plasma Physics, Boltzmannstra{\ss}e 2, D-85748 Garching, Germany \\
           David A. Kopriva \at  Department of Mathematics, The Florida State University, Tallahassee, FL 32306, USA}

\date{Received: date / Accepted: date}
\maketitle

\begin{abstract}
\textblue{We show how to modify the original Bassi and Rebay  scheme (BR1) [\textit{F. Bassi and S. Rebay, A High Order Accurate Discontinuous Finite Element Method for the Numerical Solution of the Compressible Navier-Stokes Equations, Journal of Computational Physics, 131:267--279, 1997}] to get a provably stable discontinuous Galerkin collocation spectral element method (DGSEM) with Gauss-Lobatto (GL) nodes for the compressible Navier-Stokes equations (NSE) on three dimensional curvilinear meshes. }

\textblue{Specifically, we show that the BR1 scheme can be provably stable if the metric identities are discretely satisfied, a two-point average for the metric terms is used for the contravariant fluxes in the volume, an entropy conserving split form is used for the advective volume integrals, the auxiliary gradients for the viscous terms are computed from gradients of entropy variables, and the BR1 scheme is used for the interface fluxes.}  

\textblue{Our analysis shows that even with three dimensional curvilinear grids, the BR1 fluxes do not add artificial dissipation at the interior element faces. Thus, the BR1 interface fluxes preserve the stability of the discretization of the advection terms and we get either energy stability or entropy-stability for the linear or nonlinear compressible NSE, respectively.}

\end{abstract}


\maketitle

\section{Introduction}
\textblue{In \cite{Bassi&Rebay:1997:B&F97}, Bassi and Rebay introduced a now popular discontinuous Galerkin (DG) approximation to the compressible Navier-Stokes equations (NSE).} \textblue{They extended the DG method introduced for the approximation} of first order hyperbolic conservation laws by Reed and Hill in 1973 \cite{ReedHill} and extended and popularized with the series of Cockburn and Shu et al., e.g.  \cite{cockburn1991,cockburn2001,CockburnHou&Shu,Cockburn1998199} to diffusion problems.
The method is of Galerkin finite element type, however the approximation space is piecewise polynomial with discontinuities across element interfaces in contrast to classic continuous Galerkin finite element methods. This discontinuous ansatz automatically localizes the data dependencies of the scheme as well as allows the introduction of approximate Riemann solver fluxes in surface integrals of the weak form to connect neighboring elements. Riemann solvers introduce artificial dissipation for advection problems in a natural way. The dissipation behavior of the DG approximation of advection is such that dissipation is very low for well resolved scales, but very high for coarsely resolved scales \cite{gassner2011}. This means that the scheme naturally damps away small scale oscillations for advection dominated problems and it is this difference to classical continuous Galerkin finite element methods, which are virtually dissipation free, that makes DG advantageous. 

\textblue{Bassi and Rebay rewrote the second order partial differential equations (PDEs) into an extended system of first order PDEs by introducing the gradient as a new unknown. After rewriting the system}, the standard DG ansatz analogous to the discretization of first order hyperbolic conservation laws is used. The only difference is that now two additional surface integral terms appear in the weak form and suitable choices are needed in addition to the approximate Riemann solvers of the advection terms. In the BR1 scheme, the simplest of all is used: the arithmetic mean of the viscous fluxes and the arithmetic mean of the solution (needed for the auxiliary gradient formulation). The resulting scheme is arguably the most simple variant to date for the discretization of second order terms. 

\textblue{Although popular, the Bassi and Rebay approach, and particularly the BR1 interface approximation, has never been shown to be stable for the compressible Navier-Stokes equations. When used for purely elliptic problems}, the equivalent DG scheme suffers from some bad properties \cite{Arnoldetal2002,arnold2000} such as: 1) a widened stencil, which results in a higher fill in of the operator matrix, which in turn is bad for memory consumption and efficiency of iterative linear algebra solvers; 2) the BR1 scheme is consistent, but not stable and 3) due to the symmetry of the numerical fluxes of the BR1 scheme, the convergence behavior has an odd-even behavior with non-optimal convergence. The results for pure elliptic PDEs have left a bad impression of the BR1 scheme. Many other variants for the discretization of second order terms have arisen since, e.g. \cite{Arnoldetal2002,arnold2000}, which mostly overcome the observed deficiencies of the BR1 scheme. 

When applying the BR1 scheme to the unsteady compressible Navier-Stokes Equations (NSE), the disadvantages observed in purely elliptic problems don't seem so dramatic anymore. In fact, some of the issues are naturally resolved when considering an explicit time discretization for the unsteady compressible NSE: 1) the wider stencil does not matter within a Runge-Kutta discretization, as the effective stencil size is increased anyway throughout the Runge-Kutta stages, plus the scheme is purely explicit, so no fill in considerations are necessary; 2) in combination with a de-aliased DG discretization of the advection operator, stability issues have not been observed when turbulence is severely under-resolved \cite{Gassner:2013qf}; 3) due to the upwind nature of the Riemann solver for the advection terms, there is no symmetry of the numerical fluxes at the interface and no negative influence on the observed convergence rate when considering high Reynolds number (advection dominated) flows, see e.g. \cite{FLD:FLD3943,Hindenlang201286}. In addition, the scheme is generic in the sense that it is independent of the underlying form or structure of the viscous terms, as only arithmetic means are used at the interfaces. It is also parameter free, i.e. no particular choice of any penalty constant is necessary. This makes the BR1 scheme very simple to implement for generic linear and nonlinear viscous terms on general unstructured or structured curvilinear grids. Furthermore, in the experience of the authors and analyzed in \cite{viscous_flux_choice}, the BR1 discretization allows for a relatively large explicit time step in comparison to other schemes. 

\textblue{Unfortunately, DG approximations to the Euler and Navier-Stokes equations are known to sometimes fail due to aliasing instabilities, e.g. \cite{Gassner:2013qf}. To make the approximation more robust, production codes add ad hoc stabilization procedures such as overintegration and filtering, e.g. \cite{Mengaldo201556,Gassner:2013qf}.} 

\textblue{The popularity of the scheme, plus the fact that ad hoc stabilisation is often necessary, shows that there is a need for careful analysis and construction of a provably stable version. }

\textblue{In this paper we present the full analysis of a provably stable Bassi and Rebay type approximation for the compressible NSE on curvilinear grids in three space dimensions. We present a semi-discrete analysis in a step by step manner}, i.e. the energy analysis for the linearized compressible NSE in the first part, \cref{sec:linearanalysis}, and the entropy analysis for the nonlinear equations in the second part,  \cref{sec:nonlinearanalysis}. Both parts have a similar structure: we first introduce the continuous and discrete analysis for a one dimensional scalar equation \textblue{to introduce the basic steps and concepts. This is in preparation for the extension} to the general three dimensional NSE on unstructured curvilinear grids.\\
\\ 
\textblue{From the analysis we show that a Bassi and Rebay type approximation is stable if
\begin{itemize}
\item[\textbullet] the metric identities are discretely satisfied,
\item[\textbullet]  the metric terms in the construction of the contravariant fluxes are incorporated as a two-point average,
\item[\textbullet]  the advective volume terms are discretised in a split form manner,
\item[\textbullet]  for nonlinear problems the viscous fluxes are computed in terms of the entropy variables and its gradients,
\item[\textbullet]  and the BR1 scheme is used for the interface fluxes, i.e. arithmetic means are used.
\end{itemize}}

In preparation of the two main parts of this work, we first introduce the discontinuous Galerkin collocation spectral element (DGSEM) with Gauss-Lobatto (GL) nodes in \cref{sec:standard_DGSEM}, we present the proofs for linear equations in \cref{sec:linearanalysis}, the proofs for the nonlinear equations in \cref{sec:nonlinearanalysis}, and we draw our conclusions in the last section,  \cref{sec:conclusions}.

\subsection{Nomenclature}

The analysis and proofs in this work are quite technical. For clarity, we collect notation that we use throughout this work here.\vspace{0.5cm}

\begin{tabular}{ll}
$\mathbb{P}^{N}$ & Space of polynomials of degree $\leqslant N$\\[0.05cm]
$\mathbb{I}^{N}$& Polynomial Interpolation operator\\[0.05cm]
$(x,y,z)$ &Physical space coordinates\\[0.05cm]
$\left(\xi,\eta,\zeta\right)$ &Reference space coordinates\\[0.05cm]
$\spacevec v $ & Vector in three dimensional space\\[0.05cm]
$\spacevec n= n_{1}\hat x + n_{2}\hat y + n_{3}\hat z$ & Cartesian space normal vector\\[0.05cm]
$\hat n= \hat n^{1}\hat \xi + \hat n^2\hat \eta + \hat n^3\hat \zeta$ & Reference space normal vector\\[0.05cm]
$\statevec u$ &Continuous quantity\\[0.05cm]
$\statevec {U}$ & Polynomial approximation\\[0.05cm]
$\bigstatevec f,\, \bigcontravec f$& Block vector of Cartesian flux and contravariant flux\\[0.05cm]
$\mmatrix B$& Matrix\\[0.05cm]
$\bigmatrix B$& Block matrix\\[0.05cm]
\end{tabular}

\section{The DGSEM for Compressible Viscous Flows}

Compressible viscous flows are modelled by the Navier-Stokes equations,
\begin{equation}
\label{eq:nse}
{\statevec u_t} + \sum\limits_{i = 1}^3 {\frac{{\partial {\statevec f_{i}}}}{{\partial {x_i}}} = \overRe\sum\limits_{i = 1}^3 {\frac{{\partial \statevec f_{v,i}\left( {\statevec u,{\nabla _x}\statevec u} \right)}}{{\partial {x_i}}}} }.
\end{equation}
The state vector contains the conservative variables
\begin{equation}\statevec u = \left[ {\begin{array}{*{20}{c}}
  \rho  \\ 
  {\rho \spacevec v} \\ 
  {\rho E} 
\end{array}} \right] = \left[ {\begin{array}{*{20}{c}}
  \rho  \\ 
  {\rho v_1} \\ 
  {\rho v_2} \\ 
  {\rho v_3} \\ 
  {\rho E} 
\end{array}} \right].\end{equation}
In standard form, the components of the advective flux are
\begin{equation}
\statevec f_{1}  = \left[ {\begin{array}{*{20}c}
   {\rho v_1}  \\
   {\rho v_1^2  + p}  \\
   {\rho v_1\,v_2}  \\
   {\rho v_1\,v_3}  \\
   {\rho v_1\,H}  \\

 \end{array} } \right]\quad \statevec f_{2}  = \left[ {\begin{array}{*{20}c}
   {\rho v_2}  \\
   {\rho v_2\,v_1}  \\
   {\rho v_2\,v_2  + p}  \\
   {\rho v_2\,v_3}  \\
   {\rho v_2\,H}  \\

 \end{array} } \right]\quad \statevec f_{3}  = \left[ {\begin{array}{*{20}c}
   {\rho v_3}  \\
   {\rho v_3\,v_1}  \\
   {\rho v_3\,v_2}  \\
   {\rho v_3\,v_3  + p}  \\
   {\rho v_3\,H}  \\

 \end{array} } \right],
\end{equation}
where\begin{equation}
 H = E + \frac{p}
{\rho }\quad E = e + \frac{1}
{2}\left| {\spacevec v} \right|^2 \quad e = \frac{1}
{{\gamma  - 1}}\frac{p}
{\rho }.
\end{equation}
The equations have been scaled with respect to free stream reference values so that the Reynolds number is
\begin{equation}
\mathrm{Re} = \frac{\rho_{\infty}V_{\infty}L}{\mu_{\infty}},
\end{equation}
where $L$ is the length scale and $V_{\infty}$ is the free-stream velocity. Additionally, the Mach number and Prandtl numbers are 
\begin{equation}{\mathrm{M}_\infty } = \frac{{{V_\infty }}}{{\sqrt {\gamma \mathrm{R}{T_\infty }} }},\quad \Pr  = \frac{{{\mu _\infty }{C_p}}}{{{\lambda _\infty }}}\,.\end{equation}
The viscous fluxes written in terms of the primitive variables are
\begin{equation}
  \label{eq:navierstokes_difffluxes}
  \begin{aligned}
  \statevec f_{v,1} &= \left[   0 \quad 
  {{\tau _{11}}} \quad 
  {{\tau _{12}}} \quad 
  {{\tau _{13}}} \quad 
  \left(\left( {\sum}_{j=1}^3 v_j\tau _{1j}\right) + \lambda\frac{\partial T}{\partial x}\right)\;
  \right]^T\,, \\ 
\statevec   f_{v,2} &= \left[
  0 \quad 
  {{\tau _{21}}} \quad 
  {{\tau _{22}}} \quad 
  {{\tau _{23}}} \quad 
  \left( \left( {\sum}_{j=1}^3 v_j\tau _{2j}\right) + \lambda\frac{\partial T}{\partial y}\right)\;
  \right]^T \,,\\ 
\statevec   f_{v,3} &= \left[ 
  0 \quad 
  {{\tau _{31}}} \quad 
  {{\tau _{32}}} \quad 
  {{\tau _{33}}} \quad 
  \left( \left({\sum}_{j=1}^3 v_j\tau _{3j}\right) + \lambda\frac{\partial T}{\partial z}\right)\;
  \right]^T  \,,
  \end{aligned} 
\end{equation}
where
\begin{equation}
\tau _{ij}  = \mu \left( {\frac{{\partial v_j }} {{\partial x^i }} + \frac{{\partial v_i }} {{\partial x^j }}} \right) - \frac{2}{3}\mu \left( {\nabla_{x} \cdot\spacevec v} \right)\delta _{ij} \,,\quad \lambda=\frac{\mu}{{(\gamma  - 1)\Pr \mathrm{M}_\infty ^2}}\,,
\end{equation}
and the temperature is
\begin{equation}
T = \gamma \mathrm{M}_{\infty}^{2}\frac{p}{\rho}\,.
\end{equation}

For a compact notation that will simplify the analysis, we define \emph{block vectors} (with the double arrow)
\begin{equation}
\bigstatevec{f} =
 \left[ {\begin{array}{*{20}{c}}
  {{\statevec f_1}} \\ 
  {{\statevec f_2}} \\ 
  {{\statevec f_3}} 
\end{array}} 
\right]\,,
\end{equation}
and the spatial gradient of a state as
\begin{equation}
\spacevec{\nabla}_{x} \statevec u = \left[ {\begin{array}{*{20}{c}}
  {{\statevec u_x}} \\ 
  {{\statevec u_y}} \\ 
  {{\statevec u_z}} 
\end{array}} 
\right]\,.
\end{equation}
The dot product of two block vectors is defined by
\begin{equation}
\bigstatevec f \cdot \bigstatevec g = \sum\limits_{i = 1}^3 {{{\statevec f}_i}^T{{\statevec g}_i}}\,. 
\end{equation}
Finally, the dot product of a block vector with a vector is a state vector,
\begin{equation} 
\spacevec g\cdot\bigstatevec f  = \sum\limits_{i = 1}^3 {{{ g}_i}{{\statevec f}_i}}\,. 
\end{equation}

With this notation the divergence of a flux is defined as
\begin{equation}
\spacevec\nabla_{x}  \cdot \bigstatevec f = \sum\limits_{i = 1}^3 {\frac{{\partial {\statevec f_i}}}{{\partial {x_i}}}}\, ,
\end{equation}
which allows us to write the Navier-Stokes equations compactly as
\begin{equation}
  {{\statevec u}_t} + {\spacevec\nabla _x} \cdot \bigstatevec f = \overRe{\spacevec\nabla _x} \cdot {\bigstatevec f_v}\left( {\statevec u,\spacevec\nabla_x\statevec u} \right)  \; .
\end{equation}

As part of the approximation procedure, it is customary to represent the solution gradients as a new variable to get a first order system of equations
\begin{equation}
 \begin{split}
  \label{eq:navierstokes_mixed}
  {{\statevec u}_t} + {\spacevec\nabla _x} \cdot \bigstatevec f &= \overRe{\spacevec\nabla _x} \cdot {{\bigstatevec{\statevec f}}_v}\left( {\statevec u,\bigstatevec q} \right)   \\
  \bigstatevec q &= {\spacevec\nabla _x}\statevec u\,.  
\end{split} 
\end{equation}

To set up the standard spectral element approximation, one subdivides the physical domain, $\Omega$, into $K$ non-overlapping and conforming hexahedral elements, $e^{k}$, $k=1,2,\ldots,K$. These elements can have curved faces if necessary to accurately approximate the geometry.

So that the equations can be approximated by a Legendre spectral element method, they are re-written in computational space on the reference element $E=[-1,1]^{3}$. Each element is mapped from the reference element with a mapping $\spacevec x = \spacevec X(\spacevec \xi\,)$, where $\spacevec X = X\hat x + Y\hat y + Z\hat x$ and the hats represent unit vectors. Similarly, the reference element space is represented by $\spacevec \xi = \xi\hat\xi + \eta\hat\eta + \zeta\hat\zeta$.

From the transformation, we define the three covariant basis vectors
\begin{equation}
\spacevec{a}_{i}=\partialderivative{\spacevec{X}}{\xi^{i}}\,,\quad i=1,2,3,
\end{equation}
and (volume weighted) contravariant vectors, formally written as
\begin{equation}
\mathcal{J}\spacevec{a}^{\,i}=\spacevec{a}_{j}\times\spacevec{a}_{k}\,, \quad (i,j,k)\;\text{cyclic}\,,
\end{equation}
where
\begin{equation}
\mathcal{J}=\spacevec{a}_{i}\cdot\left(\spacevec{a}_{j}\times \spacevec{a}_{k}\right),\quad (i,j,k)\;\text{cyclic}
\end{equation}
is the Jacobian of the transformation. Elements with curved faces will have basis vectors and Jacobian that vary within an element. However, even with curved elements, the contravariant coordinate vectors satisfy the metric identities \cite{Kopriva:2006er}
\begin{equation}
\sum\limits_{i = 1}^3 {\frac{{\partial \left( J{a^i_n} \right)}}{{\partial {\xi ^i}}} = 0}\,,\quad n=1,2,3\,. 
\label{eq:MetricIdentities}
\end{equation}

In terms of the reference space variables, the gradient of a function, $g$, is
\begin{equation}
{\spacevec\nabla _x}g = \frac{1}{\mathcal{J}}\sum\limits_{i = 1}^3 {\mathcal J{{\spacevec a}^{\,i}}\frac{{\partial g}}{{\partial {\xi ^i}}}} \label{eq:GradientTransformation} 
\end{equation}
and the divergence of a vector, $\spacevec g$, is
\begin{equation}
{\spacevec\nabla _x} \cdot \spacevec g = \frac{1}{\mathcal{J}}\sum\limits_{i = 1}^3 {\frac{\partial }{{\partial {\xi ^i}}}\left( {\mathcal J{{\spacevec a}^{\,i}} \cdot \spacevec g} \right)}. 
\label{eq:TransformedDiv}
\end{equation}

The transformation of the gradient and divergence can be written in terms of block vectors and block matrices of the form
\begin{equation}
\bigmatrix{B}  = \left[ {\begin{array}{*{20}{c}}
  {{\mmatrix B_{11}}}&{{\mmatrix B_{12}}}&{{\mmatrix B_{13}}} \\ 
  {{\mmatrix B_{21}}}&{{\mmatrix B_{22}}}&{{\mmatrix B_{23}}} \\ 
  {{\mmatrix B_{31}}}&{{\mmatrix B_{32}}}&{{\mmatrix B_{33}}} 
\end{array}} \right],
\end{equation}
where each $\mmatrix B_{ij}$ is a $5\times 5$ matrix.
In terms of block vectors and matrices, the transformation of the gradient \cref{eq:GradientTransformation} is
\begin{equation}
\spacevec \nabla_{x}\statevec u=\left[ 
{\begin{array}{*{20}{c}}
  {{{\statevec u}_x}} \\ 
  {{{\statevec u}_y}} \\ 
  {{{\statevec u}_z}} 
\end{array}} \right] = \frac{1}{\mathcal J}\left[ {\begin{array}{*{20}{c}}
  {\mathcal Ja_1^1\,{\mmatrix I_5}}&{\mathcal Ja_1^2\,{\mmatrix I_5}}&{\mathcal Ja_1^3\,{\mmatrix I_5}} \\ 
  {\mathcal Ja_2^1\,{\mmatrix I_5}}&{\mathcal Ja_2^2\,{\mmatrix I_5}}&{\mathcal Ja_2^3\,{\mmatrix I_5}} \\ 
  {\mathcal Ja_3^1\,{\mmatrix I_5}}&{\mathcal Ja_3^2\,{\mmatrix I_5}}&{\mathcal Ja_3^3\,{\mmatrix I_5}} 
\end{array}} \right]\left[ {\begin{array}{*{20}{c}}
  {{{\statevec u}_\xi }} \\ 
  {{{\statevec u}_\eta }} \\ 
  {{{\statevec u}_\zeta }} 
\end{array}} \right] = \frac{1}{\mathcal J} \bigmatrix M{\spacevec\nabla _\xi }\statevec u,
\end{equation}
where $\mmatrix I_{5}$ is the $5\times 5$ identity matrix.

For the divergence, \cref{eq:TransformedDiv}, each component is
\begin{equation}
\frac{\partial }{{\partial {\xi ^i}}}\left( {J{{\spacevec a}^{\,i}} \cdot \spacevec g} \right) = \frac{\partial }{{\partial {\xi ^i}}}\left( {\mathcal Ja_1^i{g_{1}} + \mathcal Ja_2^i{g_{2}} + \mathcal Ja_3^i{g_{3}}} \right),
\end{equation}
so for a block vector $\bigstatevec g$,
\begin{equation}
\spacevec\nabla_{x}  \cdot \bigstatevec g = \frac{1}{\mathcal{J}}\spacevec\nabla_{\xi}\cdot\left({\bigmatrix M^T}\bigstatevec g\right).
\end{equation}

Finally, if we define the contravariant block vector 
\begin{equation}
\bigcontravec f = \bigmatrix M^{T}\bigstatevec f\,,
\end{equation}
the transformed system of the Navier-Stokes equations \cref{eq:navierstokes_mixed} is compactly written as
\begin{equation}
\begin{split}
 \mathcal J {{\statevec u}_t} + {\spacevec\nabla _\xi } \cdot {\bigcontravec f} &= \overRe{\spacevec\nabla _\xi } \cdot {{\bigcontravec f}_v}\left( {\statevec u,{\bigstatevec q}} \right) \\
 \mathcal J \bigstatevec q &= \bigmatrix M  {\spacevec\nabla _\xi }\statevec u 
\end{split} 
\label{eq:xFormedNSEquations}
\end{equation}
on the reference element.

The spectral element approximation is derived from weak forms of the equations
\cref{eq:xFormedNSEquations}. Let us define the inner product on the reference element for state vectors
\begin{equation}{{{\left\langle {\statevec v,\statevec u} \right\rangle}_E} = \int\limits_{E}{\statevec u^{T}\statevec vd\xi d\eta d\zeta } }\,.\end{equation} 
Similarly, for block vectors,
\begin{equation}
\left\langle {\bigstatevec f,\bigstatevec g} \right\rangle_{E} = \int\limits_{E} {\sum\limits_{i = 1}^3 {\statevec f_i^T{\statevec g_i}} d\xi d\eta d\zeta }\,. 
\end{equation}
Since there should be no confusion in context, we will usually leave off the subscript $E$.
The weak forms that serve as the starting point of the approximation are created by multiplying each equation by an appropriate test function and integrating over the element. After integration by parts, the weak form of \cref{eq:xFormedNSEquations} reads as
\begin{equation}
 \begin{gathered}
  \left\langle {\mathcal{J}\statevec u,\boldsymbol\phi } \right\rangle + \int_{\partial E} {\boldsymbol\phi^{T}\left\{ {\bigcontravec f - \overRe{{\bigcontravec f}_v}} \right\} \cdot \hat n\dS}  - \left\langle {\bigcontravec f,{\spacevec\nabla _\xi }\boldsymbol\phi } \right\rangle = -\overRe\left\langle {{{\bigcontravec f}_v},\spacevec\nabla _\xi  \boldsymbol\phi} \right\rangle \hfill \\
  \left\langle {\mathcal{J}\bigstatevec q,\biggreekstatevec \psi } \right\rangle = \int_{\partial E} {{{\statevec u}^T}\left\{ {{\bigmatrix M^T}\biggreekstatevec \psi } \right\} \cdot \hat n\,\dS}  - \left\langle {\statevec u,\spacevec\nabla _\xi   \cdot \left( {{\bigmatrix M^T}\biggreekstatevec \psi } \right)} \right\rangle. \hfill \\ 
 \end{gathered} 
\end{equation}

\subsection{The Spectral Element Approximation}
\label{sec:standard_DGSEM}

To get spectral accuracy, we approximate the state vector by polynomials of degree $N$, which we represent as $\statevec U \in\mathbb{P}^{N}(E)$. The polynomials can be written in terms of the Legendre basis functions, or equivalently in terms of the Lagrange basis with nodes at the Legendre Gauss or Gauss-Lobatto points with nodal values $\statevec U_{nml}$, $n,m,l=0,1,\ldots,N$. We write the interpolation of a function $g$ through those nodes as $G=\mathbb{I}^{N}(g)$. Fluxes are also approximated with polynomials of degree $N$, represented nodally, and computed from the nodal values of the state and gradients. Derivatives are approximated by exact differentiation of the polynomial interpolants. Differentiation and interpolation do not commute, however, so $\left(\mathbb{I}^{N}(g)\right)'\ne \mathbb{I}^{N}\left(g'\right)$ \cite{CHQZ:2006,Kopriva:2009nx}.

The geometry and metric terms are also approximated with polynomials of degree $N$. Most importantly, the metric terms are computed so that the discrete metric identities \cite{Kopriva:2006er}
\begin{equation}
\sum\limits_{i = 1}^3 {\frac{{\partial {\mathbb{I}^N}\left( Ja^i_n \right)}}{{\partial {\xi ^i}}} = 0} \,,\quad n=1,2,3
\label{eq:DiscreteMetricIdentities}
\end{equation}
are satisfied. This is ensured if the metric terms are computed as
\begin{equation}
Ja_n^i =  - {{\hat x}_i} \cdot {\nabla _\xi } \times \left( {{\mathbb{I}^N}\left( {{X_l}{\nabla _\xi }{X_m}} \right)} \right)\,,\quad i = 1,2,3,\;n = 1,2,3,\;(n,m,l)\;\text{cyclic}\,.
\end{equation}
This definition ensures free stream preservation discretely \cite{Kopriva:2006er} and has already been shown to be important for stability \cite{Kopriva2016274}.

Integrals and inner products are approximated by a Legendre-Gauss or Legendre-Gauss-Lobatto quadrature. We write the quadrature in one space dimension as
\begin{equation}
 \int_{ - 1}^1 {g\left( \xi  \right)d\xi }\approx\int\limits_N {g(\xi)\,d\xi }  \equiv \sum\limits_{n = 0}^N {g\left( {{\xi _n}} \right){\omega_n}}  \,, 
\end{equation}
where $\omega_{n}$, $n=0,\ldots,N$ are the quadrature weights.
Tensor product extension is used for multiple dimensions.
We write the discrete inner product between two functions $f$ and $g$ in three space dimensions as
\begin{equation}
{\left\langle {f,g} \right\rangle_N} = \sum\limits_{n,m,l = 0}^N {{f_{nml}}{g_{nml}}{\omega_n}{\omega_m}{\omega_l}} \equiv \sum\limits_{n,m,l = 0}^N {{f_{nml}}{g_{nml}}{\omega_{nml}}}, 
\end{equation}
where  $f_{nml}=f\left(\xi_{n},\eta_{m},\zeta_{l}\right)$, etc.
An important fact to remember here is that from the definition, 
\begin{equation}
{\left\langle {{\mathbb{I}^N}(f),V} \right\rangle_N} = {\left\langle {f,V} \right\rangle_N},
\end{equation}
for any $V\in\mathbb{P}^{N}$.

The Gauss type quadratures are chosen because of their high precision. For polynomials $U$ and $V$, 
\begin{equation}
{\left\langle {U,V} \right\rangle_N} = \left\langle {U,V} \right\rangle\quad \forall \;UV \in {\mathbb{P}^{2N + \delta }},
\end{equation}
where $\delta = 1$ for the Gauss and $\delta = -1$ for the Gauss-Lobatto quadrature \cite{CHQZ:2006}. The exactness of \emph{both} Gauss and Gauss-Lobatto quadrature leads to a summation-by-parts formula,
\begin{equation}
{\int\limits_N {UV'dx}  = \left. {UV} \right|_{ - 1}^1 - \int\limits_N {U'Vdx} },
\label{eq:SumByParts_DAK}
\end{equation}
which extends to all space dimensions \cite{gassner2010}. 

In this work, we restrict ourselves to Gauss-Lobatto quadrature, where the boundary nodes are included. Including the boundary nodes simplifies the construction of stable surface terms, since no interpolation of volume data to the element surface is necessary. From summation-by-parts \cref{eq:SumByParts_DAK}, we have the discrete extended Gauss Law \cite{Kopriva:2017yg}
\begin{equation}
 {\iprodN{\spacevec\nabla_\xi  \cdot \bigcontravec F,V} = \int_{\partial E ,N} {\left(\bigcontravec F \cdot \hat n\right)\,V \dS}  - \iprodN{\bigcontravec F,\spacevec\nabla_\xi V}},
 \label{eq:DiscreteGreens_DAK}
\end{equation}
  where $\hat n$ is the reference space unit outward normal at the faces of $E$ and
\begin{equation}
 \begin{split}
 \label{eq:discrete_surfint}
 \int_{\partial E,N} {\left(\bigcontravec F \cdot \hat n\right)\,\dS} 
 &= \sum\limits_{j,k = 0}^N {\left. {{\omega_{jk}}{\contravec F^{1}}\left( {\xi    ,{\eta_j},{\zeta_k}} \right)} \right|_{\xi   =  - 1}^1}  
  + \sum\limits_{i,k = 0}^N {\left. {{\omega_{ik}}{\contravec F^{2}}\left( {{\xi_i},\eta    ,{\zeta_k}} \right)} \right|_{\eta  =  - 1}^1}  
  + \sum\limits_{i,j = 0}^N {\left. {{\omega_{ij}}{\contravec F^{3}}\left( {{\xi_i},{\eta_j},\zeta    } \right)} \right|_{\zeta =  - 1}^1} 
\\&\equiv\int\limits_N {\left. {{\contravec F^{1}}d\eta d\zeta } \right|} _{\xi    =  - 1}^1 
       + \int\limits_N {\left. {{\contravec F^{2}}d\xi d\zeta  } \right|} _{\eta   =  - 1}^1 
       + \int\limits_N {\left. {{\contravec F^{3}}d\xi d\eta   } \right|} _{\zeta  =  - 1}^1.
 \end{split}
\end{equation}
With this notation, the approximations give the discrete weak forms of the DGSEM
\begin{equation}
 \begin{gathered}
  \left\langle {J\statevec U,\boldsymbol\phi } \right\rangle_{N} 
  + \int_{\partial E,N} {\boldsymbol\phi^{T}\left\{ {\bigcontravec F - \overRe{{\bigcontravec F}_v}} \right\} \cdot \hat n\dS}  
  - \left\langle{\bigcontravec F,{\spacevec\nabla _\xi }\boldsymbol\phi } \right\rangle_{N}
  = -\overRe\left\langle {{{\bigcontravec F}_v},\spacevec\nabla_\xi \boldsymbol\phi} \right\rangle_{N} \hfill \\
  \left\langle J\bigstatevec Q, \biggreekstatevec \psi \right\rangle_{N} = 
  \int_{\partial E,N} {{{\statevec U}^T}\left(\left( {{\bigmatrix M^T}\biggreekstatevec \psi } \right) \cdot \hat n\right) \dS}  
  - \left\langle {\statevec U,\spacevec\nabla_\xi  \cdot \mathbb{I}^{N}\left( {{\bigmatrix M^T}\biggreekstatevec \psi } \right)} \right\rangle_{N}, 
 \end{gathered}
\end{equation}
where the test functions are restricted to polynomials in $\mathbb{P}^{N}$.

Elements are coupled through the boundary terms by way of \emph{numerical fluxes}, which we represent as $\contrastatevec F^{*}$, $\contrastatevec F^{*}_{v}$ and $\statevec U^{*}$ 
\begin{equation}
 \begin{gathered}
\left\langle {J\statevec U,\boldsymbol\phi } \right\rangle_{N} 
+ \int_{\partial E,N} {\boldsymbol\phi^{T}\left\{ {\contrastatevec F^{*} - \overRe{{\contrastatevec F}^{*}_v}} \right\} \dS}  
- \left\langle {\bigcontravec F,{\spacevec\nabla _\xi }\boldsymbol\phi } \right\rangle_{N} 
= -\overRe\left\langle {{{\bigcontravec F}_v},\spacevec\nabla_\xi \boldsymbol\phi} \right\rangle_{N} \hfill \\
\left\langle {J\bigstatevec Q,\biggreekstatevec \psi } \right\rangle_{N} 
= \int_{\partial E,N} {{{\statevec U}^{*,T}}\left(\left( {{\bigmatrix M^T}\biggreekstatevec \psi } \right) \cdot \hat n\right)\dS}  
- \left\langle {\statevec U,\spacevec\nabla_\xi  \cdot \mathbb{I}^{N}\left( {{\bigmatrix M^T}\biggreekstatevec \psi } \right)} \right\rangle_{N}\,. \hfill \\ 
 \end{gathered}
\end{equation}
Applying the discrete extended Gauss law \cref{eq:DiscreteGreens_DAK} to the equation for $\bigstatevec Q$ gives the final weak form of the DGSEM for the compressible Navier-Stokes equations
\begin{equation}
 \begin{gathered}
\left\langle {J\statevec U,\boldsymbol\phi } \right\rangle_{N} 
+ \int_{\partial E,N} {\boldsymbol\phi^{T}\left\{ {\contrastatevec F^{*} - \overRe{{\contrastatevec F}^{*}_v}} \right\} \dS}  
- \left\langle {\bigcontravec F,{\spacevec\nabla _\xi }\boldsymbol\phi } \right\rangle_{N} 
= -\overRe\left\langle {{{\bigcontravec F}_v},\spacevec\nabla_\xi \boldsymbol\phi} \right\rangle_{N} \hfill \\
\left\langle {J\bigstatevec Q,\biggreekstatevec \psi } \right\rangle_{N}
= \int_{\partial E,N} {\left\{{{\statevec U}^{*}}-\statevec U\right\}^{T}\left(\left( {{\bigmatrix M^T}\biggreekstatevec \psi } \right) \cdot \hat n\right)\dS}  
- \left\langle {\spacevec\nabla_{\xi}\statevec U, {{\bigmatrix M^T}\biggreekstatevec \psi }} \right\rangle_{N}\,. \hfill \\ 
 \end{gathered} 
 \label{eq:StandardDGSEM}
\end{equation}

The numerical advective flux $\contrastatevec F^{*}$ is usually computed with an approximate Riemann solver such as the Lax-Friedrichs or Roe solvers. The coupling functions for the viscous terms include the Bassi-Rebay (BR1), Local DG (LDG), Interior Penalty (IP), and others. See \cite{arnold2000} for a review of the variety of solvers used in practice for the viscous terms.

The approximation with an upwind Riemann solver for the advective flux and the BR1 scheme for the viscous terms is usually stable in practice, at least for well-resolved flows. Examples include two and three dimensional computations, e.g. \cite{FLD:FLD3943,Hindenlang201286,Daniel-A.-Nelson:2014ve}. Often, however, some kind of filtering is applied to ensure stability \cite{Frank2016132} or the volume integrals are ``overintegrated'', i.e., evaluated with quadratures $M>N$ \cite{Kirby2003,Mengaldo201556}.

In fact, despite its use in applications, \emph{the approximation \cref{eq:StandardDGSEM} is not necessarily even linearly stable}. Immediate possible contributors to instability are the use of inexact quadrature in the volume terms, or the choice of the numerical fluxes. We show in the following that with the appropriate approximation of the advective and diffusion terms, the approximation is stable using the BR1 scheme for the surface viscous terms. Then, starting from the stable approximation, we will see precisely what contributes to instability in  \cref{eq:StandardDGSEM} for underresolved flows.

\section{Linear Stability Analysis}
\label{sec:linearanalysis}

\subsection{BR1 is Stable: Linear Scalar Advection-Diffusion in 1D}

To motivate and to provide an outline of the concrete steps used in the analysis of the Navier-Stokes equations, we start with an initial boundary value problem for the simpler linear variable coefficient advection-diffusion equation
\begin{equation}
  BVP\;\left\{ 
  \begin{gathered}
  {u_t} + {(au)_x} = {\left( {b{u_x}} \right)_x}\,,\quad x \in [0,L] \hfill \\
  u(0,t) = 0 \hfill \\
  {u_x}(L,t) = 0 \hfill \\
  u(x,0) = {u_0}(x) \hfill \\ 
  \end{gathered}  \right.
  \label{eq:advDiffBVP}
\end{equation}
with $a = \text{const.} > 0$ and $b=b(x)>0$. The choice $a>0$ is solely to simplify the writing of the boundary conditions. It is not necessary in general. To guarantee stability of the advective terms with the DG approximation \cite{Gassner:2013uq}, we re-write them in split form,
\begin{equation}
{u_t} + \frac{1}{2}\left\{ {{{(au)}_x} + a{u_x} + {a_x}u} \right\} = {\left( {b{u_x}} \right)_x}\,.
\end{equation}
\textblue{Note that in the case of a constant advection speed, $a$, the split form reduces to the standard conservative form. However, we purposely leave it in split form in anticipation of the more complex proofs where the split form is necessary.}  Since $a$ is constant,
\begin{equation}
{u_t} + \frac{1}{2}\left\{ {{{(au)}_x} + a{u_x}} \right\} = {\left( {b{u_x}} \right)_x}\,.
\end{equation}
We also split the gradient of the solution as a variable itself so
\begin{equation}
 \begin{split}
  {u_t} + \frac{1}{2}\left\{ {{{(au)}_x} + a{u_x} } \right\} &= \left(bq\right)_x   \\
  q &= {u_x}\,. 
 \end{split}
\end{equation}
We then construct two weak forms
\begin{equation}
 \begin{split}
  \iprod {{u_t},\phi } + \frac{1}{2}\left\{ {\iprod{{{(au)}_x},\phi } + \iprod{a{u_x},\phi } } \right\} &= \iprod{\left(bq\right)_x,\phi }   \\
  \iprod{q,\psi } &= \iprod{{u_x},\psi }  
 \end{split},
\end{equation}
where
\begin{equation}
\iprod{u,\phi } = \int_0^L {u\phi dx}\,. 
\end{equation}

Integrating the second and last inner products in the equation for $u$ by parts yields
\begin{equation}
 \begin{split}
  \left\langle {{u_t},\phi } \right\rangle + \frac{1}{2}\left. {au\phi } \right|_0^L - \left. {bq\phi } \right|_0^L + \frac{1}{2}\left\{ {\left\langle {a{u_x},\phi } \right\rangle - \left\langle {au,{\phi _x}} \right\rangle } \right\} &=  - \left\langle {bq,{\phi _x}} \right\rangle  \\
  \left\langle {q,\psi } \right\rangle &= \left\langle {{u_x},\psi } \right\rangle \,.   
 \end{split} 
\end{equation}

\subsubsection{Continuous Energy Analysis in 1D}
The boundary value problem \cref{eq:advDiffBVP} is well-posed. When we choose $\phi=u$ and $\psi = bq$ and impose the boundary conditions,
\begin{equation}
 \begin{split}
  \frac{1}{2}\frac{d}{{dt}}{\left\| u \right\|^2} &=  - \frac{1}{2}{u^2}(L) - \left\langle {bq,{u_x}} \right\rangle   \\
 \left\langle {q,bq} \right\rangle &= \left\langle {{u_x},bq} \right\rangle \geqslant 0, 
 \end{split}
\end{equation}
where $\inorm{u}^{2}= \iprod{u,u}$.
Substituting for the $\left\langle {bq,{u_x}} \right\rangle$ term in the equation for $u$,
\begin{equation}
\frac{1}{2}\frac{d}{{dt}}{\left\| u \right\|^2} \leqslant 0.
\end{equation}
Then integrating in time over the interval $[0,T]$ and applying the initial condition gives the energy bound
\begin{equation}
\left\| {u(T)} \right\| \leqslant \left\| {{u_0}} \right\|.
\end{equation}
It is a bound of this type that we require for the stability of the discrete approximation.

\subsubsection{An Energy Stable DGSEM in 1D}
\label{sec:EnergyStableDGSEM}
To construct the DGSEM in one space dimension, we subdivide the interval into elements $e^{k}=\left[x_{k-1},x_{k}\right]\;k=1,2,\ldots,K$, where the $x_{k},\;k=0,1,\ldots,K$ are the element boundaries with $x_{0}=0$ and $x_{K}=L$. Each element is mapped onto the reference element $E=[-1,1]$ by the linear mapping
\begin{equation}
x = {x_{k - 1}} + \Delta {x_k}\frac{{\xi  + 1}}{2},
\end{equation}
where $\Delta x_{k} = x_{k}-x_{k-1}$ is the length of the element. With the change of variable, $u$ and $q$ satisfy
\begin{equation}
 \begin{gathered}
  \frac{{\Delta {x_k}}}{2}{\left\langle{{u_t},\phi } \right\rangle} + \frac{1}{2}\left\{ {{{\left\langle {{{\left( {au} \right)}_\xi },\phi } \right\rangle}} + {{\left\langle {a{u_\xi },\phi } \right\rangle}}} \right\} = {\left\langle {(bq)_{\xi},\phi } \right\rangle} \hfill \\
  \frac{{\Delta {x_k}}}{2}{\left\langle {q,\psi } \right\rangle} = {\left\langle { {u_\xi },\psi } \right\rangle}, \hfill \\ 
 \end{gathered}
 \label{eq:OrigAdvDifElemEqn} 
\end{equation}
where now
\begin{equation}
\iprod{u,\phi } = \int_{-1}^1 {u\phi d\xi}\,. 
\end{equation}

In preparation for the approximation, we integrate the first and third terms in the braces of the equation for $u$ in \cref{eq:OrigAdvDifElemEqn} by parts once, as well as the diffusion term on the right. Similarly, we integrate the equation for $q$ by parts
\begin{equation}
 \begin{gathered}
  \frac{{\Delta {x_k}}}{2}\left\langle {{u_t},\phi } \right\rangle + \left. {au\phi } \right|_{ - 1}^1 - \left. {bq\phi } \right|_{ - 1}^1 - \frac{1}{2}\left\{ { \left\langle {au,{\phi _\xi }} \right\rangle +\left\langle {u,{{\left( {a\phi } \right)}_\xi }} \right\rangle} \right\} = -{\left\langle{bq,{\phi _\xi }} \right\rangle} \hfill \\
  \frac{{\Delta {x_k}}}{2}{\left\langle {q,\psi } \right\rangle} = \left. {u \psi } \right|_{ - 1}^1 - \left\langle{u,{\psi_\xi }} \right\rangle, \hfill \\ 
 \end{gathered}
 \label{eq:OrigAdvDifElemEqn2} 
\end{equation}
We then make the following approximations
\begin{equation}
 \begin{split}
  u      &\approx U \in {\mathbb{P}^N}( - 1,1) \,,  \\
  q      &\approx Q \in {\mathbb{P}^N}( - 1,1)\,,  \\
 b \psi  &\approx {\mathbb{I}^N}\left( {b \psi } \right)\,.  
 \end{split} 
\end{equation}
Furthermore, we approximate inner products with discrete inner products that approximate the continuous ones using \emph{Gauss-Lobatto Quadrature}, and restrict $\phi$ and $\psi$ to the polynomial space. We use Gauss-Lobatto quadrature to ensure that the interpolant of the flux equals the flux of the interpolant at the element endpoints. This property will be needed later to ensure that the advective boundary terms are dissipative. Finally, elements are coupled by introducing continuous \emph{numerical fluxes} $F^{*}$ for the advective flux, $U^{*}$ for the boundary solution in the gradient equation and $Q^{*}$ for the boundary derivative. The result is the split form DGSEM approximation for the solution on an element
\begin{equation} 
\begin{split}
  \frac{{\Delta {x_k}}}{2}{\left\langle {{U_t},\phi } \right\rangle_N} 
  &+ \left. {\left\{F^{*}-bQ^{*}\phi \right\}} \right|_{ - 1}^1 
  - \frac{1}{2}\left\{ {  {{\iprod {aU,{\phi _\xi }}}_N} 
  + {{\left\langle {U,{{\left( {{\mathbb{I}^N}\left( {a\phi } \right)} \right)}_\xi }} \right\rangle}_N}} \right\} \\
  &= -{\left\langle {bQ,{\phi _\xi }} \right\rangle_N}\\
  \frac{{\Delta {x_k}}}{2}{\left\langle {Q,\psi } \right\rangle_N} 
  &= \left. {U^{*} \psi } \right|_{ - 1}^1 - {\left\langle {U,\psi_{\xi}} \right\rangle_N}. \hfill \\ 
\end{split} 
\end{equation}

Summation-by-parts applied once to each equation gives the final form of the approximation
\begin{equation}
\begin{split}
  \frac{{\Delta {x_k}}}{2}{\left\langle {{U_t},\phi } \right\rangle_N} 
  &+ \left. {\left( {{F^*} - \frac{1}{2}aU} \right)\phi } \right|_{ - 1}^1 
  - \left. {{bQ^*}\phi } \right|_{ - 1}^1 \\
  &+ \frac{1}{2}\left\{ { {{\left\langle {{{\left( {{\mathbb{I}^N}\left( {aU} \right)} \right)}_\xi },\phi } \right\rangle}_N} - {{\left\langle{U,{{\left( {{\mathbb{I}^N}\left( {a\phi } \right)} \right)}_\xi }} \right\rangle}_N}} \right\} 
  = -{\left\langle {bQ,{\phi _\xi }} \right\rangle_N} \hfill \\
  \frac{{\Delta {x_k}}}{2}{\left\langle{Q,\psi } \right\rangle_N} 
  &= \left. {\left( {{U^{*}} - U} \right) \psi } \right|_{ - 1}^1 +{\left\langle{ {U_\xi },\psi } \right\rangle_N}. \hfill \\ 
\end{split}
 \end{equation}

To show stability, we replace $\phi$ by $U$ and $\psi$ by $\mathbb{I}^{N}\left(bQ\right)$
\begin{equation}
 \begin{gathered}
  \frac{{\Delta {x_k}}}{2}{\left\langle {{U_t},U} \right\rangle_N} + \left. {\left( {{F^*} - \frac{1}{2}aU} \right)U} \right|_{ - 1}^1 - \left. {{bQ^*}U} \right|_{ - 1}^1 = - {\left\langle {bQ,{U_\xi }} \right\rangle_N} \hfill \\
  \frac{{\Delta {x_k}}}{2}{\left\langle{Q,bQ} \right\rangle_N} = \left. {\left( {{U^{*}} - U} \right)bQ} \right|_{ - 1}^1 + {\left\langle {{U_\xi },bQ} \right\rangle_N}  
 \end{gathered} 
\end{equation}
so that the element-wise energy satisfies
\begin{equation}
\begin{split}
\frac{1}{2}\frac{d}{{dt}}\frac{{\Delta {x_k}}}{2}{\left\langle {U,U} \right\rangle_N} &+ \left. {\left( {{F^*} - \frac{1}{2}aU} \right)U} \right|_{ - 1}^1 - \left. {{bQ^*}U} \right|_{ - 1}^1 - \left. {\left( {{U^{*}} - U} \right)bQ} \right|_{ - 1}^1\\&=  - \frac{{\Delta {x_k}}}{2}{\left\langle {Q,bQ} \right\rangle_N} \le 0,
\end{split}
\end{equation}
or, making the element ID explicit,
\begin{equation}
\frac{1}{2}\frac{d}{{dt}}\frac{{\Delta {x_k}}}{2}\inorm{U^{k}}_{N}^{2} + \left. {\left( {{F^*} - \frac{1}{2}aU^{k}} \right)U^{k}} \right|_{ - 1}^1 - \left. {{bQ^*}U^{k}} \right|_{ - 1}^1 - \left. {\left( {{U^{*}} - U^{k}} \right)b Q^{k}} \right|_{ - 1}^1 \leqslant 0,\end{equation}

The total energy in the domain is the sum over all of the elements
\begin{equation}
\left\| U \right\|_{N}^2 = \sum\limits_{k = 1}^{K} {\frac{{\Delta {x_k}}}{2}\left\| {{U^k}} \right\|_N^2}. 
\end{equation} 
Since the numerical quantities $F^{*}$, $Q^{*}$ and $U^{*}$ are continuous by construction, the total energy satisfies
\begin{equation}
\frac{1}{2}\frac{d}{{dt}}\inorm{U}_{N}^{2}  \leqslant \BL - \BR + \BI ,
\end{equation}
with the left and right boundary terms  
\begin{equation}
\BR = {\left. {\left\{ {\left( {{F^*} - \frac{1}{2}a{U^K}} \right){U^K} - \left[ {{bQ^*}{U^K} + {U^{*}}b{Q^K} - {U^K}b{Q^K}} \right]} \right\}} \right|_{\xi  = 1}}\,,
\end{equation}
\begin{equation}
\BL = {\left. {\left\{ {\left( {{F^*} - \frac{1}{2}a{U^1}} \right){U^1} - \left[ {{bQ^*}{U^1} + {U^{*}}b{Q^1} - {U^1}b{Q^1}} \right]} \right\}} \right|_{\xi  =  - 1}}\,,
\end{equation}
and the inner boundary term
\begin{equation}\BI = \sum\limits_\interiorfaces {  \left( {{F^*}\jump{ U}  - \frac{1}{2}a\jump{ {U^2}} } \right) - b\left\{{Q^*}\jump{ U}  + {{U^{*}}\jump{ Q}  - \jump{ UQ} }\right\}} \,,
\label{eq:LinScalBI}
\end{equation}
where 
\begin{equation}
 \jump{V}=V^{R}-V^{L}=V^{k+1}|_{\xi=-1} -  V^{k}|_{\xi=1}\label{eq:jump1D}
\end{equation} 
is the usual jump in the argument across the interface between element $k,k+1$.

We now examine each of the boundary contributions in turn. If we apply boundary states
\begin{equation}
{F^*} = 0,\quad {U^{*}} = 0,\quad {Q^*} = Q
\end{equation}
on the left, then
\begin{equation}
\BL = {\left\{ { - \frac{1}{2}a{{\left( {{U^1}} \right)}^2} - \left[ {b{Q^1}{U^1} - b{Q^1}{U^1}} \right]} \right\}_{\xi  =  - 1}} =  - \frac{1}{2}a{\left. {{{\left( {{U^1}} \right)}^2}} \right|_{\xi  =  - 1}} \leqslant 0.
\end{equation}
Applying the states
\begin{equation}
{F^*} = aU,\quad {U^{*}} = U,\quad {Q^*} = 0,
\end{equation}
on the right boundary gives
\begin{equation}
\BR = {\left\{ {\frac{1}{2}a{{\left( {{U^K}} \right)}^2} - \left[ {{U^K}b{Q^K} - {U^K}b{Q^K}} \right]} \right\}_{\xi  = 1}} = \frac{1}{2}a{\left( {{U^K}} \right)^2} \geqslant 0
\end{equation}
The boundary conditions applied in this way are therefore dissipative.

We are now left with only the contribution of the jumps at the interfaces. The contribution of the advective part is non-positive. With the upwind value $F^{*}=aU^{L}$,
\begin{equation}
a\left({U^{L}\jump{ U}  - \frac{1}{2}\jump{ {U^2}} }\right)=-\frac{1}{2}a\jump{U}^{2}\leqslant 0\,,
\end{equation}
which leaves only the interface contribution of the diffusion terms to bound.

The BR1 scheme \cite{Bassi&Rebay:1997:B&F97} computes the interface values as
\begin{equation}
 \begin{gathered}
  {U^{*}}\left(U^{L},U^{R}\right) = \frac{{{U^L} + {U^R}}}{2} \equiv \avg{ U}  \hfill \\
  {Q^*}\left(Q^{L},Q^{R}\right) = \frac{{{Q^L} + {Q^R}}}{2} \equiv \avg{ Q}\,.  \hfill \\ 
 \end{gathered}
\end{equation}
Using the identity
\begin{equation}
\jump{wv}=\avg{w}\jump{v} + \avg{v}\jump{w},
\end{equation}
the interface contribution from the approximation of the diffusion terms in \cref{eq:LinScalBI} is
\begin{equation}
\avg{Q}\jump{U} + \avg{U}\jump{Q} - \jump{UQ}=0\,.
\end{equation}

With all boundary and interface terms accounted for, 
\begin{equation}
\frac{1}{2}\frac{d}{{dt}}{\left\| U \right\|_{N}^2} \leqslant 0\,,
\end{equation}
and integrated over the time interval $[0,T]$ leads to the desired bound
\begin{equation}
{\left\| {U(T)} \right\|_N} \leqslant {\left\| {{U_0}} \right\|_N}\,.
\end{equation}

In summary, we see that using a split form for the advection terms and the BR1 scheme for the diffusion, the DGSEM is stable for the advection-diffusion equation.

\subsection{BR1 is Stable: 3D Linearized Compressible Navier-Stokes Equations}
\label{sec:BR1_isstable_3Dlinearized_NSE}
We follow the same steps as in the previous section to now show that the DGSEM with the BR1 scheme is linearly stable for the compressible Navier-Stokes equations, provided that the advective terms and the physical boundary conditions are approximated stably. A roadmap for the development of well-posed problems and stable approximations was recently presented by  Nordstr\"om \cite{Nordstrom:2016jk}. Since we are interested in the influence on stability of the BR1 scheme in this paper, we will consider only the first two steps in that roadmap: symmetrization of the equations and the energy method.

We will leave the approximation of the physical boundary conditions to a future paper and assume here that they are properly posed and implemented in a stable manner. Starting from the stable variant of the DGSEM, we will then see why the standard approximation may be unstable even for smooth flows when underresolved.
 
\subsubsection{Continuous Energy Analysis in 3D}

As with the advection-diffusion approximation, we will prove stability using an energy method. An important difference between the two is in the step between \cref{eq:OrigAdvDifElemEqn} and \cref{eq:OrigAdvDifElemEqn2} which used the fact that the adjoint of a scalar is itself.

The Navier-Stokes equations linearized about a \emph{constant state} can be written in the form
\begin{equation}
{\statevec u_t} + \sum\limits_{j = 1}^3 {\frac{{\partial {\mmatrix A_j}\statevec u}}{{\partial {x_j}}}}  = \overRe\sum\limits_{i = 1}^3 {\frac{\partial }{{\partial {x_i}}}\left( {\sum\limits_{j = 1}^3 {{\mmatrix B_{ij}}\frac{{\partial \statevec u}}{{\partial {x_j}}}} } \right)},
\label{eq:CNSE}
 \end{equation}
 where $\statevec u = [\delta\rho \;\delta v_{1}\; \delta v_{2}\; \delta v_{3}\; \delta p]^{T}$ represents the perturbation from the reference values. 
The coefficient matrices $\mmatrix A_{j}$ and $\mmatrix B_{ij}$ are \emph{constant} in the linear approximation of the equations.  We use the primitive variable formulation because that is what one usually implements. The system is known to be symmetrizable by a single constant symmetrization matrix, $\mmatrix S$, and there are multiple symmetrizers \cite{Isi:A1981Lw20700001} to choose from that will enable us to apply the energy method. We write the symmetrized matrices as 
$\mmatrix A^{s}_{j}=\mmatrix S^{-1}\mmatrix A_{j}\mmatrix S = \left(\mmatrix A^{s}_{j}\right)^{T}$ and 
$\mmatrix B^{s}_{ij}=\mmatrix S^{-1}\mmatrix B_{ij}\mmatrix S = \left(\mmatrix B^{s}_{ij}\right)^{T}$.
 Explicit representations of the symmetrizer and coefficient matrices are written down in \cite{Nordstrom:2005qy}.
 

To again simplify the notation for the use in the analysis, we define 
a block vector of matrices, e.g.
\begin{equation}\bigstatevec A = \left[ {\begin{array}{*{20}{c}}
  {{\mmatrix A_1}} \\ 
  {{\mmatrix A_2}} \\ 
  {{\mmatrix A_3}} 
\end{array}} \right]
\end{equation}
and  the diagonal block matrix and full block matrix
\begin{equation}
\bigmatrix S = \left[ {\begin{array}{*{20}{c}}
  {{\mmatrix S}}\;&0&0 \\ 
  0&{{\mmatrix S}}\;&0 \\ 
  0&0&{{\mmatrix S}}\; 
\end{array}} \right]\,,\quad 
\bigmatrix B = \left[ {\begin{array}{*{20}{c}}
  {{\mmatrix B_{11}}}&{{\mmatrix B_{12}}}&{{\mmatrix B_{13}}} \\ 
  {{\mmatrix B_{21}}}&{{\mmatrix B_{22}}}&{{\mmatrix B_{23}}} \\ 
  {{\mmatrix B_{31}}}&{{\mmatrix B_{32}}}&{{\mmatrix B_{33}}} \\ 
\end{array}} \right].
\end{equation}
Then the product rule applied to the divergence of the flux in \cref{eq:CNSE} can be written as
\begin{equation}
{\spacevec\nabla _x} \cdot \bigstatevec f = \left({{\spacevec\nabla }_x} \cdot \bigstatevec A\right)\statevec u +  \left(\bigstatevec A\right)^T{\spacevec\nabla _x}\statevec u,
\end{equation}
where 
\begin{equation}
\bigstatevec f = \left[ {\begin{array}{*{20}{c}}
  {{\mmatrix A_1}\statevec u} \\ 
  {{\mmatrix A_2}\statevec u} \\ 
  {{\mmatrix A_3}\statevec u} 
\end{array}} \right]\,,\quad \left(\bigstatevec A\right)^T=\left[\mmatrix A_1 \;\mmatrix A_2 \;\mmatrix A_3\right ].
\end{equation}
The nonconservative advective form of the linearized Navier-Stokes equations can therefore be written as
\begin{equation}
\statevec u_{t} + \left({{\spacevec\nabla }_x} \cdot \bigstatevec A\right)\statevec u + \left(\bigstatevec A\right)^T{\spacevec\nabla _x}\statevec u = \overRe\spacevec\nabla_{x}\cdot \left(\bigmatrix B\spacevec\nabla_{x}\statevec u\right).
\end{equation}
Averaging the conservative and nonconservative forms gives the split form of the PDE
\begin{equation}
\statevec u_{t} + \frac{1}{2}\left\{{\spacevec\nabla _x} \cdot \bigstatevec f+\left({{\spacevec\nabla }_x} \cdot \bigstatevec A\right)\statevec u + \left(\bigstatevec A\right)^T{\spacevec\nabla _x}\statevec u\right\} = \overRe\spacevec\nabla_{x}\cdot \left(\bigmatrix B\spacevec\nabla_x\statevec u\right).
\end{equation}
\textred{Because the coefficient matrices are constant the split form is not necessary. It is presented only in preparation for the discrete analysis. In the discrete approximation the coefficient matrices are still constant, but they get multiplied by the metric terms of the curvilinear elements. This indirectly introduces variable coefficients, even for the linear NSE.} 

\textred{Proceeding with the continuous analysis, the next step is to drop the divergence of the coefficient matrices,  ${{\spacevec\nabla }_x} \cdot \bigstatevec A$, since it is zero. We will see later in the discrete analysis that this step needs additional attention, since it depends on properties of the discrete metric terms.}

As usual, we construct a weak form by multiplying the split form by a test function and integrating over the domain. In inner product notation,
\begin{equation}
\iprod{\statevec u_{t},\boldsymbol\phi} + \frac{1}{2}\left\{\iprod{\spacevec\nabla_{x}\cdot \bigstatevec f,\boldsymbol\phi}  + \iprod{\left(\bigstatevec A\right)^T{\spacevec\nabla _x}\statevec u,\boldsymbol\phi}\right\} = \overRe\iprod{\spacevec\nabla_{x}\cdot \left(\bigmatrix B\spacevec\nabla_x\statevec u\right),\boldsymbol\phi}.
\end{equation}
As before, we introduce the intermediate block vector $\bigstatevec q = \spacevec\nabla_{x} \statevec u$ to get the first order system
\begin{equation}
 \begin{split}
  \iprod{\statevec u_{t},\boldsymbol\phi} + \frac{1}{2}\left\{\iprod{\spacevec\nabla_{x}\cdot \bigstatevec f,\boldsymbol\phi}  + \iprod{\left(\bigstatevec A\right)^T{\spacevec\nabla _x}\statevec u,\boldsymbol\phi}\right\} &= \overRe\iprod{\spacevec\nabla_{x}\cdot \left(\bigmatrix B\bigstatevec q\right),\boldsymbol\phi}
   \\
  \iprod{\bigstatevec q,\biggreekstatevec\psi}  &= \iprod{\spacevec\nabla_{x} \statevec u,\biggreekstatevec\psi}   \,.
 \end{split} 
 \label{eq:SplitLinearSystem0}
\end{equation}
Then we apply the extended Gauss law \cref{eq:DiscreteGreens_DAK}
to the flux divergence terms to separate surface and volume contributions
\begin{equation}
 \begin{split}
\left\langle {{{\statevec u}_t},\boldsymbol\phi } \right\rangle &+ \int_{\partial \Omega } {\left( \frac{1}{2}\left(\bigstatevec f\cdot\spacevec n \right)
- \overRe\left(\left(\bigmatrix B{{\spacevec\nabla }_x}\statevec u\right)\cdot\spacevec n\right) \right)^{T}\boldsymbol\phi \dS}  \\&+ \frac{1}{2}\left\{ {\left\langle {{{\spacevec\nabla }_x}\statevec u,\bigstatevec{f}^{\,(T)}\left( \boldsymbol\phi  \right)} \right\rangle - \left\langle {\bigstatevec f,{{\spacevec\nabla }_x}\boldsymbol\phi } \right\rangle } \right\} \\
&=  - \overRe\left\langle {\bigmatrix B\bigstatevec q,{{\spacevec\nabla }_x}\boldsymbol\phi } \right\rangle ,
 \end{split}
 \label{eq:SplitLinearSystem}
\end{equation}
where
\begin{equation}
 \bigstatevec f^{\,(T)}\left(\boldsymbol\phi \right) = \left[ {\begin{array}{*{20}{c}}
  {{\mmatrix A^{T}_1}\boldsymbol\phi } \\ 
  {{\mmatrix A^{T}_2}\boldsymbol\phi } \\ 
  {{\mmatrix A^{T}_3}\boldsymbol\phi } 
\end{array}} \right]\,,
\end{equation}
and $\spacevec n$ is the physical space outward normal to the surface.


With suitable boundary and initial conditions, the equations are well-posed. First, we set $\boldsymbol\phi = \left(\mmatrix S^{-1}\right)^{T}\mmatrix S^{-1}\statevec u$ in \cref{eq:SplitLinearSystem}, which includes symmetrization as part of the test function. Then
\begin{equation}
 \begin{split}
\left\langle{\mmatrix S^{-1}{{\statevec u}_t},\mmatrix S^{-1}\statevec u } \right\rangle 
&+
 \int_{\partial \Omega } {\left(\frac{1}{2}{\mmatrix S^{-1}\left(\bigstatevec{f}\cdot\spacevec n\right) - \overRe\mmatrix S^{-1}\left(\left(\bigmatrix B{{\spacevec\nabla }_x}\statevec u\right)\cdot\spacevec n\right)} \right)^{T}\mmatrix S^{-1}\statevec u \dS}  
\\&
+ \frac{1}{2}\left\{ {\left\langle {{{\spacevec\nabla }_x}\statevec u,\bigstatevec f^{\,(T)}\left( \left(\mmatrix S^{-1}\right)^{T}\mmatrix S^{-1}\statevec u  \right)} \right\rangle 
- \left\langle {\bigmatrix S^{-1}\bigstatevec f,{{\spacevec\nabla }_x}\left(\mmatrix S^{-1}\statevec u\right) } \right\rangle 
} \right\} 
\\&
=  - \overRe\left\langle {\bigmatrix S^{-1}\bigmatrix B\bigstatevec q,{{\spacevec\nabla }_x}\left(\mmatrix S^{-1}\statevec u\right )} \right\rangle .
 \end{split}
 \label{eq:SplitLinearSystemStab2}
\end{equation}
Let us define the symmetric state vector as $\statevec u^{s}=\mmatrix S^{-1}\statevec u$ and examine the terms in \cref{eq:SplitLinearSystemStab2} separately.
First,
\begin{equation}
\left\langle {\mmatrix S^{-1}{{\statevec u}_t},\mmatrix S^{-1}\statevec u } \right\rangle=\frac{1}{2}\frac{d}{dt}\inorm{\statevec u^{s}}^{2}.
\end{equation}
Next, we consider the volume term
\begin{equation}
\left\langle {\bigmatrix S^{-1}\bigmatrix B\bigstatevec q,{{\spacevec\nabla }_x}\left(\mmatrix S^{-1}\statevec u \right)} \right\rangle 
= 
\left\langle{\bigmatrix B^{s}\bigstatevec q^{s},{{\spacevec\nabla }_x}\statevec u^{s} } \right\rangle.
\end{equation}
Making the changes on the boundary terms,
\begin{equation} 
\int\limits_{\partial \Omega } {\left( \frac{1}{2}{\mmatrix S^{-1}\left(\bigstatevec f\cdot\spacevec n\right) - \overRe\mmatrix S^{-1}\left(\left(\bigmatrix B{{\spacevec\nabla }_x}\statevec u\right)\cdot\spacevec n\right )} \right)^{T}\mmatrix S^{-1}\statevec u \dS}  =
 \int\limits_{\partial \Omega } 
 {\left( \frac{1}{2}\left(\bigstatevec f^{s}\cdot\spacevec n\right) 
  - \overRe\left(\left(\bigmatrix B^{s}{{\spacevec\nabla }_x}\statevec u^{s}\right)\cdot\spacevec n\right) \right)^{T}\statevec u^{s} \dS}  ,
\end{equation}
where
\begin{equation}
\bigstatevec f^{s} = \left[ {\begin{array}{*{20}{c}}
  {\mmatrix A_1^s\statevec u^{s}} \\ 
  {\mmatrix A_2^s\statevec u^{s}} \\ 
  {\mmatrix A_3^s\statevec u^{s}} 
\end{array}} \right].
\end{equation}
The most interesting terms are the volume flux terms. The solution flux term is
\begin{equation}
 \left\langle {\bigmatrix S^{-1}\bigstatevec f,{{\spacevec\nabla }_x}\left(\mmatrix S^{-1}\statevec u\right) } \right\rangle = \left\langle {\bigstatevec f^{s},{{\spacevec\nabla }_x}\statevec u^{s} } \right\rangle,
\end{equation}
whereas the test function flux term is
\begin{equation}
\begin{split}
\left\langle {{{\spacevec\nabla }_x}\statevec u,\bigstatevec f^{\,(T)}\left( \left(\mmatrix S^{-1}\right)^{T}\mmatrix S^{-1}\statevec u  \right)} \right\rangle &=
\left\langle {\bigmatrix S{{\spacevec\nabla }_x}\mmatrix S^{-1}\statevec u,{\left(\bigmatrix S^{-1}\bigstatevec f\left( \statevec u^{s}  \right)\right)^{T}}} \right\rangle
\\&=\left\langle {{\spacevec\nabla }_x}\mmatrix S^{-1}\statevec u,{\left(\bigmatrix S^{-1}\bigstatevec f\left( \statevec u^{s}  \right)\bigmatrix S\right)^{T}}\right\rangle
\\&= \left\langle {{\spacevec\nabla }_x}\statevec u^{s},\bigstatevec f^{s}\left( \statevec u^{s}  \right)\right\rangle.
\end{split}
\label{eq:ExactFluxCancellation}
 \end{equation}

Next, we set   $\biggreekstatevec\psi = \left(\bigmatrix S^{-1}\right)^{T}\bigmatrix S^{-1}\bigmatrix B\bigstatevec q$ in the second equation of \cref{eq:SplitLinearSystem0} 
\begin{equation}
\begin{split}
  \iprod{\bigstatevec q,\left(\bigmatrix S^{-1}\right)^{T}\bigmatrix S^{-1}\bigmatrix B\bigstatevec q}
  &=  \iprod{\spacevec\nabla_{x} \statevec u,\left(\bigmatrix S^{-1}\right)^{T}\bigmatrix S^{-1}\bigmatrix B\bigstatevec q}
  \\ \text{i.e.}\quad
   \iprod{\bigstatevec q^{s},\bigmatrix B^{s}\bigstatevec q^{s}}&= \iprod{\spacevec\nabla_{x} \statevec u^{s},\bigmatrix B^{s}\bigstatevec q^{s}}.
  \end{split}
\end{equation}

Gathering all the terms, the flux volume terms cancel leaving
\begin{equation}
\begin{split}
\frac{1}{2}\frac{d}{dt}\inorm{\statevec u^{s}}^{2}&+ \int_{\partial \Omega } {\left( \frac{1}{2}\left(\bigstatevec f^{s}\cdot\spacevec n\right)
   - \overRe\left(\left(\bigmatrix B^{s}{{\spacevec\nabla }_x}\statevec u^{s}\right)\cdot\spacevec n\right) \right)^{T}\statevec u^{s} \dS}   \\&= - \overRe\left\langle\bigstatevec q^{s},\bigmatrix B^{s}\bigstatevec q^{s} \right\rangle\leqslant 0.
\end{split}
\end{equation}
We see, then,  that the growth in the energy, defined as the $\mathbb{L}^{2}$ norm, is determined by the boundary integral,
\begin{equation}
\frac{1}{2}\frac{d}{dt}\inorm{\statevec u^{s}}^{2}\leqslant -\int_{\partial \Omega } {\left( \frac{1}{2}\left(\bigstatevec f^{s}\cdot\spacevec n\right)
   - \overRe\left(\left(\bigmatrix B^{s}{{\spacevec\nabla }_x}\statevec u^{s}\right)\cdot\spacevec n\right) \right)^{T}\statevec u^{s} \dS} .
\end{equation}
Integrating in time over the interval $[0,T]$,
\begin{equation}
\inorm{\statevec u^{s}(T)}^{2}\leqslant\inorm{\statevec u(0)}-\int_{0}^{T}{\int_{\partial \Omega } {\left( \left(\bigstatevec f^{s}\cdot\spacevec n\right)
   - \twooverRe\left(\left(\bigmatrix B^{s}{{\spacevec\nabla }_x}\statevec u^{s}\right)\cdot\spacevec n\right) \right)^{T}\statevec u^{s} \dS}}.
\end{equation}

To properly pose the problem, initial and boundary data must be specified. The value at $t=0$ is replaced by initial data $\statevec u_{0}$. As for the boundary terms, Ref. \cite{Nordstrom:2005qy} shows that they can be split, in characteristic fashion, into incoming and outgoing information with boundary data specified along the incoming characteristics 
\begin{equation}
\int_{\partial \Omega } {\left( \left(\bigstatevec f^{s}\cdot\spacevec n\right)
   - \twooverRe\left(\left(\bigmatrix B^{s}{{\spacevec\nabla }_x}\statevec u^{s}\right)\cdot\spacevec n\right) \right)^{T}\statevec u^{s}\dS}  =   \int_{\partial \Omega } {{\statevec w^{ + T}}{\Lambda ^ + }{\statevec w^ + }\dS}  - \int_{\partial \Omega } {{\statevec g^T}\left| {{\Lambda ^ - }} \right|\statevec g\dS} ,
\end{equation}
where $\Lambda^{+}>0$ and $\Lambda^{-}<0$. We will assume here that boundary data $\statevec g = 0$ and hence
\begin{equation}
\int_{\partial \Omega } {\left( \left(\bigstatevec f^{s}\cdot\spacevec n\right)
   - \twooverRe\left(\left(\bigmatrix B^{s}{{\spacevec\nabla }_x}\statevec u^{s}\right)\cdot\spacevec n\right) \right)^{T}\statevec u^{s}\dS}  \geqslant 0 ,
\end{equation}
so that
\begin{equation}
\inorm{\statevec u^{s}(T)}\leqslant\inorm{\statevec u^{s}_{0}}.
\end{equation}
Finally, since $\statevec u^{s} = \mmatrix S^{-1}\statevec u$, $\statevec u = \mmatrix S\statevec u^{s}$, 
\begin{equation}\frac{1}{{{{\left\| \mmatrix S \right\|}_2}}}\left\| \statevec u \right\| \leqslant \left\| {{\statevec u^s}} \right\| \leqslant {\left\| {{\mmatrix S^{ - 1}}} \right\|_2}\left\| \statevec u \right\|,\end{equation}
and therefore
\begin{equation}
\inorm{\statevec u(T)}\leqslant C\inorm{\statevec u_{0}}.
\label{eq:ContinuousBound}
\end{equation}
It is a bound like \cref{eq:ContinuousBound} that we seek for a stable discontinuous Galerkin approximation.

\subsubsection{An Energy Stable DGSEM in 3D}
\label{sec:ESDGEM_NSE3D}

In terms of the reference space variables, the linearized Navier-Stokes equations become
\begin{equation}
\mathcal{J}{{\statevec u}_t} + {{\spacevec\nabla }_\xi } \cdot \left( {{\bigmatrix M^T}\bigstatevec f} \right) = \overRe{\spacevec\nabla _\xi } \cdot \left( {\frac{1}{\mathcal{J}}{\bigmatrix M^T}\bigmatrix B\bigmatrix M{\spacevec\nabla _\xi }\statevec u} \right).
\end{equation}
Let us define, then, the gradient vector with the intermediate variable $\bigstatevec q$, and the contravariant viscous flux by $\bigcontravec f_{v}=\bigmatrix M^{T}\bigmatrix B\bigstatevec q$ to write the equations in reference space as
\begin{equation}
 \begin{split}
  \mathcal{J}{{\statevec u}_t} + {{\spacevec\nabla }_\xi } \cdot{\bigcontravec f} &= \overRe{\spacevec\nabla _\xi } \cdot  \bigcontravec f_{v}   \\
  \mathcal{J}\bigstatevec q &= \bigmatrix M{\spacevec\nabla _\xi }\statevec u\,.  
 \end{split} 
\end{equation}
With the product rule, we also construct the split advective form
\begin{equation}
 \begin{split}
  \mathcal{J}{{\statevec u}_t} + \frac{1}{2}\left\{ {{{\spacevec\nabla }_\xi } \cdot {\bigcontravec f} + \left( {{{\spacevec\nabla }_\xi } \cdot \bigcontravec A} \right)\statevec u + \left(\bigcontravec A\right)^T{{\spacevec\nabla }_\xi }\statevec u} \right\} &= \overRe{\spacevec\nabla _\xi } \cdot \bigcontravec f_{v}  \\
  \mathcal{J}\bigstatevec q &= \bigmatrix M{\nabla _\xi }\statevec u \,,   
\end{split}
\end{equation}
where $\bigcontravec  A = \bigmatrix M^{T}\bigstatevec A$.

 \textred{We note that for general curvilinear grids, the metrics terms in the matrix $\bigmatrix M$ depend on space. Hence, the transformed coefficient matrices $\bigcontravec  A$ are not constant anymore. Thus, the discrete divergence ${{{\spacevec\nabla }_\xi } \cdot \bigcontravec A}$ is not automatically zero as in the continuous case. However, with the metric identities \cref{eq:MetricIdentities} the discrete divergence of the transformed coefficient matrices is exactly zero. If \cref{eq:MetricIdentities} doesn't hold aliasing due to the erroneous discrete divergence of the transformed coefficient matrices could disrupt the stability of the method \cite{Kopriva2017}. It follows that the metric identities are a crucial ingredient for the stability of the discretization as will be seen in the analysis that follows.}

We then form the discontinuous Galerkin approximation as usual by replacing integrals with Gauss-Lobatto quadratures, boundary quantities by numerical ones, solutions by polynomial interpolants and restrict the test functions to the polynomial space. The result is the formal statement of the DG approximation (c.f. \cref{eq:StandardDGSEM})
\begin{equation}
\begin{split}
\left\langle {{J{\statevec U}_t},\boldsymbol\phi } \right\rangle_{N} &+ \int_{\partial E,N } {\left( {{\contrastatevec F}^{*} - \overRe\contrastatevec F^{*}_{v}} \right)^{T}\boldsymbol\phi \dS}  
\\&- \frac{1}{2}\left\{ {\left\langle {\statevec U,{{\spacevec\nabla }_{\xi}}\cdot\bigcontravec F^{(T)}\left( \boldsymbol\phi  \right)} \right\rangle_{N} 
+ \left\langle {\bigcontravec F,{{\spacevec\nabla }_{\xi}}\boldsymbol\phi } \right\rangle_{N} 
} \right\} \\&=  - \overRe\left\langle {\bigcontravec F_{v},{{\spacevec\nabla }_{\xi}}\boldsymbol\phi } \right\rangle_{N} ,
\end{split}
\label{eq:WeakDiscreteForm}
\end{equation}
where $\bigcontravec F^{(T)} = \mathbb{I}^{N}\left(\bigmatrix M^{T}\mathbb{I}^{N}\left(\bigstatevec f^{(T)}\left(\boldsymbol\phi\right)\right)\right)$, and
\begin{equation} 
\begin{split}
\iprod{ J\bigstatevec Q,{{\biggreekstatevec\psi}}}_{N} &= \int_{\partial E,N} {{{\statevec U}^{*,T}}\left(\left({\bigmatrix M^T}\biggreekstatevec\psi \right) \cdot \hat n\right)\,\dS}  - \left\langle {\statevec U,\spacevec\nabla  \cdot \mathbb{I}^{N}\left( {{\bigmatrix M^T}\biggreekstatevec\psi } \right)} \right\rangle_{N}
\\&
= \int_{\partial E,N} {\left\{{{\statevec U}^{*}-\statevec U}\right\}^{T}\left(\left({\bigmatrix M^T}\biggreekstatevec\psi \right) \cdot \hat n\right)\,\dS}  + \left\langle{{\spacevec\nabla }_{\xi}}{\statevec U, {{\bigmatrix M^T}\biggreekstatevec\psi }} \right\rangle_{N}.
\end{split}
\label{eq:DiffTermsDGSEM}
\end{equation}

The equation for the approximate solution $\statevec U$, \cref{eq:WeakDiscreteForm}, is commonly known as the \emph{weak form} of the approximation. If we apply the extended discrete Gauss law \cref{eq:DiscreteGreens_DAK} one more time to the inner product $\left\langle{\bigcontravec F,{{\spacevec\nabla }_{\xi}}\boldsymbol\phi } \right\rangle_{N}$, we get the algebraically equivalent \cite{gassner2010} \emph{strong form}
\begin{equation}
\begin{split}
\left\langle{{J{\statevec U}_t},\boldsymbol\phi } \right\rangle_{N} &+ \int_{\partial E,N } {\left( {\left\{{\contrastatevec F}^{*}-\frac{1}{2}\left(\bigcontravec F\cdot\hat n\right)\right\}  - \overRe\contrastatevec F^{*}_{v}} \right)^{T}\boldsymbol\phi \dS}  
\\&+ \frac{1}{2}\left\{ { \left\langle {{\spacevec\nabla }_{\xi}}\cdot{\bigcontravec F\left(\statevec U\right),\boldsymbol\phi } \right\rangle_{N} -
\left\langle {\statevec U,{{\spacevec\nabla }_{\xi}}\cdot\bigcontravec F^{(T)}\left( \boldsymbol\phi  \right)} \right\rangle_{N}} \right\} \\&=  - \overRe\left\langle {\bigcontravec F_{v},{{\spacevec\nabla }_{\xi}}\boldsymbol\phi } \right\rangle_{N} .
\end{split}
\label{eq:DSFormOfDGSEM}
\end{equation}

To assess stability of the approximation, we follow the same steps as to show well-posedness.  First we set $\biggreekstatevec{\psi} =  \left(\bigmatrix S^{-1}\right)^{T}\bigmatrix S^{-1}\bigmatrix B\bigstatevec Q$ in \cref{eq:DiffTermsDGSEM}. 
Using the fact that $\left(\bigmatrix S^{-1}\right)^{T}$ commutes with $\bigmatrix M^{T}$,
\begin{equation} 
\iprod{ J\bigmatrix S^{-1}\bigstatevec Q,\bigmatrix S^{-1}\bigmatrix B\bigstatevec Q}_{N} 
= 
\iprod{ J\bigstatevec Q^{s},\bigmatrix B^{s}\bigstatevec Q^{s}}_{N}
= \int_{\partial E,N} {\left\{{{\statevec U}^{s,*}-\statevec U^{s}}\right\}^{T}\left(\bigcontravec F^{s}_{v}  \cdot \hat n\right)\,\dS}  +\left\langle {{\spacevec\nabla }_{\xi}}{\statevec U^{s}, \bigcontravec F^{s}_{v}} \right\rangle_{N},
\end{equation}
where $\bigcontravec F^{s}_{v}=\mathbb I^{N}\left(\bigmatrix M^{T}\bigmatrix B^{s}\bigstatevec Q^{s}\right)$. Next, 
we set $\boldsymbol\phi = \left(\mmatrix S^{-1}\right)^{T}\mmatrix S^{-1}\statevec U = \left(\mmatrix S^{-1}\right)^{T}\statevec U^{s}$ in \cref{eq:DSFormOfDGSEM}.  The advective volume terms cancel, for (c.f. \cref{eq:ExactFluxCancellation})
 \begin{equation}
 \begin{split}
 \left\langle {{\spacevec\nabla }_{\xi}}\cdot{\bigcontravec F\left(\statevec U\right),\left(\mmatrix S^{-1}\right)^{T}\mmatrix S^{-1}\statevec U } \right\rangle_{N} &-
\left\langle {\statevec U,{{\spacevec\nabla }_{\xi}}\cdot\bigcontravec F^{(T)}\left( \left(\mmatrix S^{-1}\right)^{T}\mmatrix S^{-1}\statevec U  \right)} \right\rangle_{N}
\\&= 
 \left\langle {{\spacevec\nabla }_{\xi}}\cdot{\bigcontravec F^{s}\left(\statevec U^{s}\right),\statevec U^{s} } \right\rangle_{N} -
\left\langle{\statevec U^{s},{{\spacevec\nabla }_{\xi}}\cdot\bigcontravec F^{s}\left( \statevec U^{s}  \right)} \right\rangle_{N}\\&=0\,.
\end{split}
\end{equation}
Therefore,
\begin{equation}
\begin{split}
\frac{1}{2}\frac{d}{dt} \inorm{\statevec U^{s}}^{2}_{J,N}\equiv\left\langle {{J{\statevec U^{s}}_t},\statevec U^{s} } \right\rangle_{N}  =  &-  \int_{\partial E,N } {\left( {\left\{{\contrastatevec F}^{s,*}-\frac{1}{2}\left(\bigcontravec F^{s}\cdot\hat n\right)\right\}  - \overRe\contrastatevec F^{s,*}_{v}} \right)^{T}\statevec U^{s} \dS}
\\&+\overRe\int_{\partial E,N} {\left\{{{\statevec U}^{s,*}-\statevec U^{s}}\right\}^{T}\left(\bigcontravec F^{s}_{v} \cdot \hat n\right)\dS}\\&- \overRe\iprod{ J\bigstatevec Q^{s},\bigmatrix B^{s}\bigstatevec Q^{s}}_{N}\,.
\end{split}
\end{equation}
Separating the advective and viscous boundary terms, the elemental contribution to the total energy is
\begin{equation}
\begin{split}
\frac{1}{2}\frac{d}{dt} \inorm{\statevec U^{s}}^{2}_{J,N}=&- 
 \int_{\partial E,N } {{\left\{{\contrastatevec F}^{s,*}-\frac{1}{2}\left(\bigcontravec F^{s}\cdot\hat n\right)\right\}}^{T}\statevec U^{s} \dS}\\&
+  \overRe\int_{\partial E,N } {\left\{\contrastatevec F_{v}^{s,*,T}\statevec U^{s} +\statevec U^{s,*,T}\left(\bigcontravec F^{s}_{v}\cdot \hat n\right) -\statevec U^{s,T}\left(\bigcontravec F^{s}_{v}  \cdot \hat n\right)\right\} \dS}
\\&- \overRe\iprod{ J\bigstatevec Q^{s},\bigmatrix B^{s}\bigstatevec Q^{s}}_{N}.
\end{split}
\end{equation}

The total energy is found by summing over all of the elements. When all the element contributions are summed, the interior faces get contributions
from two elements. For a conforming mesh as assumed here, the contributions match pointwise. Also when the elements are conforming, the outward facing 
normals at each face point in precisely opposite directions. If we designate the element on one side of a shared face (arbitrarily chosen) as the ``master'' and the other as 
the ``slave'' then we can represent the normal contribution of the contravariant flux on the master element side as $\bigcontravec F\cdot\hat n = \contrastatevec F^\text{master}_{n}$. Since the
outward normal at a face is in either the $\hat \xi^{i}$ or $-\hat\xi^{i}$ direction, $\contrastatevec F^\text{master}_{n}$ is proportional to the contravariant flux for that coordinate, $\pm\contrastatevec F^{i}$. Along that \emph{same} direction, but with opposite sign is the contribution from the slave element side, $\contrastatevec F^\text{slave}_{n}$. The sum of the contributions then is represented in terms of the jump on the master element side
\begin{equation}
 \label{eq:jump_master_slave}
\jump{ \contrastatevec F_{n}}=\contrastatevec F^\text{slave}_{n} - \contrastatevec {F}^\text{master}_{n}\,.
\end{equation}
Note that this notation mimics the 1D case \cref{eq:jump1D}, where the ``right'' side of the interface was the slave element side at $\xi=-1$ of element $k+1$ and the ``left'' side corresponded to the master element side at $\xi=1$ of element $k$.

Let $\statevec U^{s,k}$ be the (symmetric) solution vector on element $k$. Then summing over all elements,
\begin{equation}
\begin{split}
\frac{1}{2}\frac{d}{dt} \sum\limits_{k=1}^{K}\inorm{\statevec U^{s,k}}^{2}_{J,N}=& 
 \sum\limits_\interiorfaces \int_{N } {\left\{{\contrastatevec F}_{n}^{s,*,T}\jump{\statevec U^{s}}-\frac{1}{2}\jump{\left(\contrastatevec F_{n}^{s}\right)^{T}\statevec U^{s}} \right\}\dS}\\&
- \overRe \sum\limits_\interiorfaces \int_{N } {\left\{\contrastatevec F_{v,n}^{s,*,T}\jump{\statevec U^{s}}+\statevec U^{s,*,T}\jump{\bigcontravec F^{s}_{v,n} }-\jump{\statevec U^{s,T}\bigcontravec F^{s}_{v,n}}\right\} \dS}
\\&-\overRe \sum\limits_{k=1}^{K}\iprod{ J\bigstatevec Q^{s,k},\bigmatrix B^{s}\bigstatevec Q^{s,k}}_{N}
\\& + \PBT\,,
\end{split}
\label{eq:LinearEnergySumBound}
\end{equation}
where PBT represents the physical boundary terms, which we assume are dissipative, i.e. $\PBT\leqslant 0$. Note that the interior faces always have the master element side orientation.

Numerical fluxes are used to resolve two discontinuous states $\statevec U^{s,L}$ and $\statevec U^{s,R}$ and the viscous fluxes.
The advective flux can be split according the wave directions relative to the normal of the master side
\begin{equation}
\left(\bigcontravec f^{s}\cdot \hat n\right) = \left(\left(\bigmatrix M^{T}\bigstatevec {A}^{s}\right)\cdot\hat n\right)\statevec u^{s}
=\left(\bigcontravec {A}^{s}\cdot \hat n\right)\statevec u ^{s}\;\equiv\; \tilde {\mmatrix A}_n^{s}\statevec u^{s} =  \left(\tilde {\mmatrix A}_n^{s,+}+\tilde {\mmatrix A}_n^{s,-}\right)\statevec u^{s}\,,
\end{equation}
where
\begin{equation}
\tilde{\mmatrix A}_n^{s,\pm} = \frac{1}{2}\left(\tilde{\mmatrix A}_n^{s}\pm \left| \tilde{\mmatrix A}_n^{s}\right| \right)\,.
\end{equation}
From that splitting, we can write the numerical advective flux as
\begin{equation}
\contrastatevec F_{n}^{s,*}\left(\statevec U^{s,L}, \statevec U^{s,R}\right)=\frac{\tilde{\mmatrix A}_n^{s}\statevec U^{s,L} + \tilde{\mmatrix A}_n^{s}\statevec U^{s,R}}{2}+\frac{\sigma}{2} \left| \tilde{\mmatrix A}_n^{s}\right|\left(\statevec U^{s,L} - \statevec U^{s,R}\right) = \tilde{\mmatrix A}_n^{s}\avg{\statevec U^{s}} - \frac{\sigma}{2}\left| \tilde{\mmatrix A}_n^{s}\right|\jump{\statevec U^{s}}\,.
\label{eq:LinearNumericalFlux}
\end{equation}
The fully upwind flux corresponds to $\sigma = 1$, whereas $\sigma = 0$ gives the centered flux. 

With either the upwind or central numerical flux \cref{eq:LinearNumericalFlux}, the contribution of the advective fluxes at the faces is dissipative. For any two state vectors,
\begin{equation}
\begin{split}
\jump{\statevec a^{T}\statevec b} &= \jump{\sum\limits_{m=1}^{5}a_{m}b_{m}} = \sum\limits_{m=1}^{5}\jump{a_{m}b_{m}} \\&= \sum\limits_{m=1}^{5}\left( \avg{a_{m}}\jump{b_{m}}+\jump{a_{m}}\avg{b_{m}}\right) = \avg{\statevec a}^{T}\jump{\statevec b} + \jump{\statevec a}^{T}\avg{\statevec b}\,.
\end{split}
\label{eq:KnightsFormerlyKnownAsLemma1}
\end{equation}
Therefore
\begin{equation}
\begin{split}
\jump{\left(\contrastatevec F^{s}_{n}\right)^{T}\statevec U^{s}} &= \avg{\contrastatevec F^{s}_{n}}^{T}\jump{\statevec U^{s}} + \jump{\contrastatevec F^{s}_{n}}^{T}\avg{\statevec U^{s}}
\\&=\avg{\statevec U^{s}}^{T}\tilde{\mmatrix A}_n^{s}\jump{\statevec U^{s}} + \jump{\statevec U^{s}}^{T}\tilde{\mmatrix A}_n^{s}\avg{\statevec U^{s}} 
\\&= 2\avg{\statevec U^{s}}^{T}\tilde{\mmatrix A}_n^{s}\jump{\statevec U^{s}}\,,
\end{split}
\end{equation}
so
\begin{equation}
\contrastatevec F_{n}^{s,*,T}\jump{\statevec U^{s}}-\frac{1}{2}\jump{\left(\contrastatevec F^{s}_{n}\right)^{T}\statevec U^{s}} = -\frac{\sigma}{2}\jump{\statevec U^{s}}^{T}\left| \tilde{\mmatrix A}_n^{s}\right|\jump{\statevec U^{s}}\leqslant 0\,,
\end{equation}
and the contribution of the advective interface terms to the energy in \cref{eq:LinearEnergySumBound} is nonpositive. 

We are now left with bounding the viscous surface contributions in \cref{eq:LinearEnergySumBound}. Using the BR1 scheme for the surface values for the viscous terms,
the surface contribution coming from the viscous terms at each node is
\begin{equation}
\avg{\contrastatevec F^{s}_{n}}^{T}\jump{\statevec U^{s}} + \avg{\statevec U^{s}}^{T}\jump{\contrastatevec F^{s}_{n}}-\jump{\statevec U^{s,T}\contrastatevec F^{s}_{v,n}}\,.
\end{equation}
Replacing the jump in the product using \cref{eq:KnightsFormerlyKnownAsLemma1}, the surface contribution of the viscous terms due to the BR1 scheme vanishes because
\begin{equation}
\avg{\contrastatevec F^{s}_{n}}^{T}\jump{\statevec U^{s}} + \avg{\statevec U^{s}}^{T}\jump{\contrastatevec F^{s}_{n}}- \avg{\statevec U^{s}}^{T}\jump{\contrastatevec F^{s}_{n}} - \jump{\statevec U^{s}}^{T}\avg{\contrastatevec F^{s}_{n}}=0\,.
\end{equation}

We are now in the position to write the total energy bound. We define the total energy over the domain by
\begin{equation}
\inorm{\statevec U^{s}}^{2}_{N} = \sum\limits_{k=1}^{K}\inorm{\statevec U^{s,k}}^{2}_{J,N} .
\end{equation}
Since  $\bigmatrix B^{s}\ge 0$ \cite{Nordstrom:2005qy} and $J>0$, we can also define the physical dissipation over the domain by the broken norm
\begin{equation}
\inorm{\bigstatevec Q^{s}}_{\bigmatrix B^{s},N}^{2} \equiv \sum\limits_{k=1}^{K}\iprod{J\bigstatevec Q^{s,k},\bigmatrix B^{s}\bigstatevec Q^{s,k}}_{N}\geqslant 0\,.
\end{equation}
Then
\begin{equation}
\frac{1}{2}\frac{d}{dt} \inorm{\statevec U^{s}}^{2}_{N} = -\overRe\inorm{\bigstatevec Q^{s}}_{\bigmatrix B^{s},N}^{2} - 
\frac{\sigma}{2} \sum\limits_\interiorfaces \int_{N } { \jump{\statevec U^{s}}^{T}\left| \tilde{\mmatrix A}_n^{s}\right|\jump{\statevec U^{s}}\dS}+ \PBT\,.
\label{eq:DGSEMBound}
\end{equation}

Under the assumption that the physical boundary terms are approximated stably, all terms on the right side of \cref{eq:DGSEMBound} are nonpositive. Using the equivalence of the norms, we have the stability bound
\begin{equation}
\inorm{\statevec U}_{N}\leqslant C \inorm{\statevec U_{0}}_{N}.
\end{equation}
We see in \cref{eq:DGSEMBound} that dissipation includes physical dissipation plus, for $\sigma > 0$, artificial dissipation added at element surfaces that depends on the size of the jumps in the solution. Again, as in the scalar problem, we see that the BR1 scheme itself has no effect on the energy. Stability is therefore guaranteed if the conditions at the physical boundaries are properly posed and implemented stably.

\subsection{Why is Instability Seen in Practice?}

For smooth flows, stability of the approximation actually has nothing to do with BR1, and everything to do with the approximation of the advective terms where instability might be a result of
aliasing instability.
The BR1 approximation itself does not contribute to instability, but it also does not add any stabilizing dissipation.

If the standard form of the DGSEM is used even with the symmetric equations, then the bound on the energy is (c.f. \cite{Kopriva:2017yg})
\begin{equation}
\begin{split}
\frac{1}{2}\frac{d}{dt} \inorm{\statevec U^{s}}^{2}_{N} = &-\overRe\inorm{\bigstatevec Q^{s}}_{\bigmatrix B^{s},N}^{2} - 
\frac{\sigma}{2} \sum\limits_\interiorfaces \int_{N } { \jump{\statevec U^{s}}^{T}\left| \tilde{\mmatrix A}_n^{s}\right|\jump{\statevec U^{s}}\dS} \\&+
 \frac{1}{2}\sum\limits_{k=1}^{K}\left|{\iprod {\left\{ {{\mathbb{I}^N}\left( \left(\bigcontravec A^{s}\right) \right)^T   \spacevec\nabla \statevec U^{s,k} -\spacevec\nabla  \cdot \bigcontravec F^{s}(\statevec U^{s,k})} \right\},\statevec U^{s,k}}_N }\right|+ \PBT\,.
\end{split}
\label{eq:DGSEMBoundWithAliasing}
\end{equation}
The additional volume term is due to aliasing. It is the amount by which the product rule fails to hold. The sign of the aliasing error is indeterminate, and \cref{eq:DGSEMBoundWithAliasing} shows that the physical diffusion or the dissipation associated with the Riemann solver have the ability to counterbalance the product rule error and stabilize the scheme. 
For well resolved solutions, the aliasing error will be spectrally small, making it likely that the physical and interface dissipations are sufficiently large for stabilization. For under resolved solutions, the aliasing errors may be too large for the approximation to be stable. 
For large Reynolds numbers, the physical dissipation may be too small. The artificial dissipation due to the Riemann solver might be sufficiently large, depending on which solver is chosen. (For example, a Lax-Friedrichs Riemann solver will be more dissipative than the exact upwind solver.) 
Finally, changing from the BR1 to another scheme, coupled with a more dissipative Riemann solver \emph{might} be enough to stabilize the aliasing term. But, ultimately, the key to a stable DGSEM is the stable approximation of the advective terms.

\section{Nonlinear Stability Analysis}
\label{sec:nonlinearanalysis}
We now show how to construct \emph{nonlinearly} stable discontinuous Galerkin spectral element approximations, stable in the sense that the mathematical entropy is bounded by the initial value. We first analyze the Dirichlet problem for the scalar, one dimensional Burgers equation to motivate the steps. We then do the same for the nonlinear Navier-Stokes equations, postponing a study of the physical boundary conditions as we did for the linear counterpart to a later paper. We therefore focus only on the influence of the interior element boundary approximations. We will see two differences from the linear analysis of the Navier-Stokes equations. First, the Euler advective terms do not have a split form, but are otherwise approximated by a special two-point flux. A relationship between the two point flux and split forms is known for the Burgers equation, so we motivate the relationship in the next section. The second difference is that the viscous terms are to be expressed in terms of gradients of the entropy variables rather than the solution state vector.

\subsection{BR1 is Stable: The Nonlinear Viscous Burgers Equation in 1D}
\label{sec:BR1_isstable_1D_Burger}
To motivate the analysis of the nonlinear 3D Navier-Stokes equations, we analyze the DGSEM approximation of the initial boundary-value problem for the scalar, nonlinear viscous Burgers equation in one space dimension
\begin{equation}
BVP\;\left\{ 
 \begin{gathered}
  {u_t} + {f(u)_x} = {\left( {b(u){u_x}} \right)_x}\quad x \in [0,L] \hfill \\
  u(0,t)=u(L,t)=0 \hfill \\
  u(x,0) = {u_0}(x) \hfill 
 \end{gathered}  \right.
 \label{eq:visBurgersBVP}
\end{equation}
where $f(u) = u^2/2$ and $b(u)> 0$ is a positive viscosity function.

In contrast to linear problems, where we use energy estimates, we use entropy to define stability for nonlinear problems \cite{tadmor2003}. For the Burgers equation, we use the entropy function
\begin{equation}
s(u) \equiv \frac{u^2}{2},
\label{eq:Burgers_entropy)}
\end{equation}
and the entropy variable
\begin{equation}
w(u) \equiv \frac{ds}{du} = u\,.
\label{eq:BurgersEntropyVariable}
\end{equation}
The entropy variable contracts the advective term of the Burgers equation, i.e. 
\begin{equation}
w(u)\,f(u)_x = u\,f(u)_x = u\,\frac{df}{du}\,u_x = u\,u\,u_x = \left(\frac{u^3}{3}\right)_x\equiv f_x^\ent(u)\,,
\label{eq:BurgEntropyFluxDeriv}
\end{equation}
where $f^\ent(u)=\frac{u^3}{3}$ is the entropy flux. Similarly, $wu_{t}=s_{t}$.

For the following continuous analysis, we rewrite the second order problem as before into a first order system. We formally rewrite the viscous flux $b(u) u_x$ in terms of the derivative of the entropy variable
\begin{equation}
b(u) u_x = \hat{b}(u) w_x\,,
\end{equation}
where, due to the specific choice of entropy variables, $w=u$ and $\hat{b}(u)=b(u)> 0$. The system then reads as 
\begin{equation}
 \begin{split}
  {u_t} + {f(u)_x} &= {\left( {\hat{b}(u)\,{q}} \right)_x} \\
 q &= w_x \,.
 \end{split}
 \label{eq:mixed_visBurgers}
\end{equation}

\subsubsection{Continuous Entropy Analysis in 1D}

The entropy for \cref{eq:visBurgersBVP} is bounded in time by the initial entropy. To show this, we multiply the first equation in \cref{eq:mixed_visBurgers} by the entropy variable $w$ rather than the solution $u$,  and the second equation by the viscous flux $f_v \equiv\hat{b}(u)\,q$, and integrate over the domain to get the two weak forms
\begin{equation}
 \begin{split}
  \iprod{w(u),u_t} + \iprod{w(u),f(u)_x} &= \iprod{w(u), f_{v,x}} \\
 \iprod{q, f_v} &= \iprod{w_x,f_v}\,,
 \end{split}
 \label{eq:mixed_visBurgers_inetgral1}
\end{equation}
where $\iprod{\cdot,\cdot}$ is the $L^{2}$ inner product on the interval $[0,L]$. We then use the relations for the entropy function \cref{eq:Burgers_entropy)}, entropy variable \cref{eq:BurgersEntropyVariable}, the entropy flux \cref{eq:BurgEntropyFluxDeriv}, and integration by parts for the viscous volume integral in the first equation of \cref{eq:mixed_visBurgers_inetgral1} to get
\begin{equation}
 \begin{split}
  \iprod{s_t(u),1} + \iprod{f^\ent(u)_x,1} &= \left. w(u)\,f_v\right|_0^L - \iprod{w_x,f_v}\\
 \iprod{q, f_v} &= \iprod{w_x,f_v}\,.
 \end{split}
 \label{eq:mixed_visBurgers_inetgral2}
\end{equation}
We can insert the second equation of \cref{eq:mixed_visBurgers_inetgral2} into the first to eliminate the derivative of the entropy variable leaving
\begin{equation}
  \iprod{s_t(u),1} + \iprod{f^\ent(u)_x,1} = \left. w(u)\,f_v\right|_0^L - \iprod{q, f_v}\,.
\label{eq:mixed_visBurgers_entropy}
\end{equation}
We note that the viscosity coefficient $\hat{b}(u)$ is always positive and thus $- \iprod{q, f_v}=- \iprod{q, \hat{b}(u)\,q}$ is guaranteed nonpositive. Therefore we can bound \cref{eq:mixed_visBurgers_entropy} as
\begin{equation}
 \frac{d}{dt} \overline{s}\leqslant \left. \left\{-f^\ent(u)+w(u)\,f_v\right\}\right|_0^L\,,
\label{eq:mixed_visBurgers_entropy1}
\end{equation}
where the total entropy is
\begin{equation}
\overline{s}\equiv\int\limits_0^L s(u)\,dx =  \iprod{s(u),1}\,.
\end{equation}
When we insert the specific form of the entropy function \cref{eq:Burgers_entropy)}, entropy variable, \cref{eq:BurgersEntropyVariable}, and entropy flux for the Burgers equation into \cref{eq:mixed_visBurgers_entropy1} ,
\begin{equation}
 \frac{d}{dt}\int\limits_0^L \frac{u^2}{2}\,dx \leqslant  \left.\left\{-  \left(\frac{u^3}{3}\right) + u\,\hat{b}(u)\,q\right\}\right|_0^L\,,
\end{equation}
so that when the boundary conditions, $u=0$, are applied,
\begin{equation}
 \frac{d}{dt}\int\limits_0^L \frac{u^2}{2}\,dx \leqslant  0\,,
\end{equation}
which says that the total entropy over the domain at any time is bounded by its initial value.

\subsubsection{An Entropy-Stable DGSEM in 1D}

Using the notation introduced in  \cref{sec:EnergyStableDGSEM}, we construct the DGSEM using the usual steps: (i) we multiply the equations in \cref{eq:mixed_visBurgers} with test functions and integrate over an element, (ii) we use integration by parts to separate boundary and volume contributions, (iii) we approximate the quantities with Lagrange polynomials of degree N constructed with $N+1$ GL nodes, (iv) we approximate the inner products with quadrature rules with the same N+1 GL nodes as for the Lagrange ansatz, (v) we insert numerical surface fluxes at the element interfaces, and (vi) we use summation-by-parts for the first equation to get the strong form DGSEM of the viscous Burgers equation
\begin{equation}\begin{split}
\frac{\Delta x_k}{2}\iprod{U_t,\phi}_N +  \left.(F^*-F)\phi\right|_{-1}^1 + \iprod{F_\xi,\phi}_N &= \left.(F_{v}^{*}-F_{v})\phi\right|_{-1}^1+ \iprod{F_{v,\xi},\phi}_N\\
\frac{\Delta x_k}{2} \iprod{Q, \psi}_N &= \left. W^*\psi\right|_{-1}^1 - \iprod{W,\psi_\xi}_N\,.
\end{split}
\label{eq:mixed_visBurgers_DGSEM1}
\end{equation}
In \cref{eq:mixed_visBurgers_DGSEM1}, $U,Q$, $W$ and $F$ are the polynomial approximations of the solution, the solution gradient, the entropy variable, and the advective flux collocated at the Gauss-Lobatto nodes. Furthermore, we introduce the shorthand notation of the viscous flux $F_v$, which is a polynomial approximation of the term ${\hat{b}(u){q}} $, i.e. $F_v = {\mathbb{I}^N}({\hat{b}(U){Q}})$. The quantities marked with $*$ are the numerical surface flux approximations, which, as before, depend on state values (and gradients) from the left and right of the element interface.

Unfortunately, the standard form of the DGSEM, \cref{eq:mixed_visBurgers_DGSEM1}, is unstable, independent of the choice of the numerical surface fluxes because of the way the advective terms are discretized. The problem, again, is the aliasing introduced in the discretization of the volume terms of the advective fluxes $\iprod{F_\xi,\phi}_N$. The aliasing error causes a sink/source in the discrete entropy and cannot be bounded. The result is that entropy can pile up during coarsely resolved simulations and can even lead to a blow up of the solution. 

In \cite{gassner_skew_burgers}, the fix to the aliasing problem was a specific reformulation of the volume integrals by using the  skew-symmetry strategy, where the discrete volume integral of the advective terms was replaced with  
\begin{equation}
\label{eq:Burgers_entropyconserving}
\iprod{F_\xi,\phi}_N = \frac{1}{2}\,\iprod{{\mathbb{I}^N}(U^2)_\xi,\phi}_N \approx  \frac{1}{3}\,\iprod{{\mathbb{I}^N}(U^2)_\xi + {\mathbb{I}^N}(U\,U_\xi),\phi}_N\,.
\end{equation}
We note that for well resolved approximations, the difference between the two volume integral approximations is spectrally small. However, with severe under resolution, the difference is important. In fact, the second term in \cref{eq:Burgers_entropyconserving} is constructed so that it mimics the continuous entropy analysis of the previous section. That is, when we insert the collocation interpolant of the entropy variable $W$($=U$ for the Burgers equation) for the test function $\phi$, it was shown in \cite{gassner_skew_burgers} and also in \cref{app:Proof__volint_entropy_Burgers} here, that the volume contributions with the skew-symmetric derivative approximation can be rewritten into an equivalent surface integral contribution
\begin{equation}
\label{eq:discrete_entropystable_volumeintegral_Burgers}
\frac{1}{3}\iprod{{\mathbb{I}^N}(U^2)_\xi + {\mathbb{I}^N}(U\,U_\xi),W}_N = \left. F^\ent\right|_{-1}^1\,,
\end{equation}
where $F^\ent$ is the interpolant of the continuous entropy flux, $f^\ent$ with nodes at the Gauss-Lobatto points. 

Unfortunately, the strategy of using skew-symmetric forms of the advective terms does not directly generalize to the compressible nonlinear Euler equations.  Carpenter and Fisher et al., e.g. \cite{skew_sbp2,carpenter_esdg} have, however, established a general discrete framework to determine the volume integral for general hyperbolic conservation laws so that crucial conditions like \cref{eq:discrete_entropystable_volumeintegral_Burgers} hold. It is therefore instructive to see the relationship between the 
split form and the two point flux form introducted in \cite{skew_sbp2,carpenter_esdg} in a case where they are equivalent.

The key is the summation-by-parts property of the underlying spatial operator and the existence of a numerical two-point flux $F^{\ec} = F^{\ec}(U,V)$ with the property
\begin{equation}
\label{eq:ec_condition_numflux0}
F^{\ec}(U,V)\,\left(W(U)-W(V)\right) - \left(F(U)\,W(U) - F(V)\,W(V)\right)+\left(F^\ent(U) - F^\ent(V)\right)= 0\,,
\end{equation}
or in shorthand notation,
\begin{equation}
\label{eq:ec_condition_numflux}
F^{\ec}\,\jump{W} - \jump{F\,W}+\jump{F^\ent} = 0\,,
\end{equation}
which gives entropy conservation when used in a first order finite volume discretization of the nonlinear hyperbolic conservation law, e.g. \cite{tadmor2003}. 

Carpenter and Fisher et al. use the summation-by-parts property and \cref{eq:ec_condition_numflux} to construct a high order accurate volume integral. To bring their scheme into the form used in this work, we introduce a special derivative projection operator 
\begin{equation}
\mathbb{D}(F)^{\ec}(\xi)\equiv \sum\limits_{k=0}^N \ell_k'(\xi)\,2\,F^{\ec}(U(\xi),U(\xi_k)).
\end{equation}
We note that $\mathbb{D}(F)^{\ec}$ is a function that depends on $\xi$, but is in general not a polynomial of degree $N$. Like Carpenter and Fisher et al. we replace the standard volume term approximation with this derivative projection 
 \begin{equation}
\label{eq:Burgers_discrete_entropyconserving}
\iprod{F_\xi,\phi}_N \approx \iprod{\mathbb{D}(F)^{\ec},\phi}_N\,,
\end{equation}
which has the property (see \cref{app:Proof__volint_entropy_Burgers})
\begin{equation}
\label{eq:discrete_entropystable_volumeintegral}
\iprod{\mathbb{D}(F)^{\ec},W}_N = \left. F^\ent\right|_{-1}^1\,,
\end{equation}
which is analogous to the property of the skew-symmetric volume integral \cref{eq:discrete_entropystable_volumeintegral_Burgers}. 

We now show that the two strategies are indeed equivalent. For the Burgers equation, the entropy conserving two-point flux has the form \cite{Gassner:2013ol}
\begin{equation}
F^{\ec}_\burg(U,V) = \frac{1}{6}\left(U^2 + U\,V+V^2\right)\,,
\label{eq:BurgersTwoPointFlux}
\end{equation}
so that the derivative projection operator is
\begin{equation}
\mathbb{D}(F)^{\ec}_\burg(\xi) = \sum\limits_{k=0}^N \ell_k'(\xi)\,\frac{1}{3}\,\left(U^2(\xi) + U(\xi)\,U(\xi_k)+U^2(\xi_k)\right)\,.
\end{equation}
When projected onto the test function,
\begin{equation}
 \begin{split}
\iprod{\mathbb{D}(F)_\burg^{\ec},\phi}_N  &= \sum\limits_{n=0}^N \omega_n\,\phi(\xi_n)\,\mathbb{D}(F)^{\ec}_\burg(\xi_n)\\ 
&= \sum\limits_{n=0}^N \omega_n\,\phi(\xi_n)\,\sum\limits_{k=0}^N \ell_k'(\xi_n)\,\frac{1}{3}\,\left(U^2(\xi_n) + U(\xi_n)\,U(\xi_k)+U^2(\xi_k)\right)\\
&= \sum\limits_{n=0}^N \omega_n\,\phi(\xi_n)\,\frac{1}{3}\,\left[U^2(\xi_n)\,\sum\limits_{k=0}^N \ell_k'(\xi_n)+U(\xi_n)\,\sum\limits_{k=0}^N \ell_k'(\xi_n)U(\xi_k)+\sum\limits_{k=0}^N\ell_k'(\xi_n)U^2(\xi_k)\right]\\
&= \sum\limits_{n=0}^N \omega_n\,\phi(\xi_n)\,\frac{1}{3}\,\left[U(\xi_n)\,\sum\limits_{k=0}^N \ell_k'(\xi_n)U(\xi_k)+\sum\limits_{k=0}^N\ell_k'(\xi_n)U^2(\xi_k)\right]\\
&=  \frac{1}{3}\,\iprod{\left({\mathbb{I}^N}\left(U^2\right)\right)_\xi + {\mathbb{I}^N}(U\,U_\xi),\phi}_N\,,
 \end{split}
 \label{eq:splitFormBurgEquivalence}
\end{equation}
where we use the consistency of the polynomial derivative, i.e. that the derivative of the constant one $\sum\limits_{k=0}^N \ell_k'(\xi_n) =0$. Thus, for the Burgers equation, the use of $\mathbb{D}(F)^{\ec}_\burg$ is equivalent to the approximation of the split form,
\begin{equation}
 uu_{\xi}=\frac{1}{3}\left\{\left(u^{2}\right)_{\xi}+ uu_{\xi}\right\}\,.
\end{equation}
The goal now is to show that in combination with a stable discretization of the advection terms, using BR1 for the viscous terms leads to a stable approximation in the sense that it is entropy conserving. 

We get a stable approximation by replacing the standard volume integral of the advection term with the entropy conserving version \cref{eq:Burgers_discrete_entropyconserving} and use an entropy conserving two-point flux $F^{\ec,*}=F^{\ec}$ satisfying \cref{eq:ec_condition_numflux} as the numerical surface flux for the advection terms
\begin{equation}
 \begin{split}
\frac{\Delta x_k}{2}\iprod{U_t,\phi}_N +  \left.(F^{\ec,*}-F)\phi\right|_{-1}^1 + \iprod{\mathbb{D}(F)_\burg^{\ec},\phi}_N &= \left.(F_{v}^{*}-F_v)\phi\right|_{-1}^1+ \iprod{F_{v,\xi},\phi}_N \\
\frac{\Delta x_k}{2} \iprod{Q, \psi}_N &= \left. W^*\psi\right|_{-1}^1 - \iprod{W,\psi_\xi}_N\,.
 \end{split}
 \label{eq:mixed_visBurgers_ECDGSEM1}
\end{equation}

To show that entropy is conserved, we set $\phi=W$ and $\psi=F_v$ in \cref{eq:mixed_visBurgers_ECDGSEM1}
\begin{equation}
 \begin{split}
\frac{\Delta x_k}{2}\iprod{U_t,W}_N +  \left.(F^{\ec,*}-F)W\right|_{-1}^1 + \iprod{\mathbb{D}(F)_\burg^{\ec},W}_N &= \left.(F_{v}^{*}-F_v)W\right|_{-1}^1+ \iprod{F_{v,\xi},W}_N \\
\frac{\Delta x_k}{2} \iprod{Q, F_v}_N &= \left. W^*F_v\right|_{-1}^1 - \iprod{W,F_{v,\xi}}_N\,.
 \end{split}
 \label{eq:mixed_visBurgers_EC1}
\end{equation}

From the collocation of interpolation and numerical integration (i.e. a point-wise multiplication at each GL node) and the assumption of analytical time integration,
\begin{equation}
\frac{\Delta x_k}{2}\iprod{U_t,W}_N = \frac{\Delta x_k}{2}\iprod{S_t,1}_N\,,
\end{equation}
for the first term in the first equation of \cref{eq:mixed_visBurgers_EC1}, where $S$ is the discrete polynomial approximation of the entropy $s(u)$. Next, we use \cref{eq:discrete_entropystable_volumeintegral} to replace the volume contribution of the advective flux by a surface contribution. Finally, we use the second equation of \cref{eq:mixed_visBurgers_EC1} to replace the viscous volume term contribution in the first to get
\begin{equation}
\frac{\Delta x_k}{2}\iprod{S_t,1}_N +  \left.\left\{(F^{\ec,*}-F)W+ F^\ent\right\}\right|_{-1}^1 = \left.\left\{(F_v^*-F_v)W+W^*F_v\right\}\right|_{-1}^1 - \frac{\Delta x_k}{2} \iprod{Q, F_v}_N.
\label{eq:mixed_visBurgers_EC2}
\end{equation}
The last term on the right of \cref{eq:mixed_visBurgers_EC2} is always non-positive since
\begin{equation}
\begin{split}
- \frac{\Delta x_k}{2} \iprod{Q, F_v}_N &= - \frac{\Delta x_k}{2} \sum_{n=0}^N Q(\xi_n) F_v(\xi_n)\,\omega_n \\&= - \frac{\Delta x_k}{2} \sum_{n=0}^N Q(\xi_n) \hat{b}(U(\xi_n)) Q(\xi_n)\omega_{n},
\end{split}
\end{equation}
and the viscosity coefficient $\hat{b}$ is always positive. 

With the element ID explicit, we have the first intermediate result that on an element $e^{k}$
\begin{equation}
\frac{\Delta x_k}{2}\iprod{S^k_t,1}_N  \leqslant -\left\{  \left.(F^{\ec,*}-F^k)W^k+ F^{\ent,k}\right\}\right|_{-1}^1 + \left.\left\{(F_v^*-F_v^k)W^k+W^*F_v^k\right\}\right|_{-1}^1\,.
\label{eq:mixed_visBurgers_EC3}
\end{equation}
When we sum over all elements, we get the temporal change of the total discrete entropy
\begin{equation}
\frac{d}{dt} \overline{S} \equiv \sum\limits_{k=1}^K\frac{\Delta x_k}{2}\iprod{S^k_t,1}_N= \sum\limits_{k=1}^K\frac{\Delta x_k}{2}\sum\limits_{n=0}^N S_t^k(\xi_n)\omega_n\,,
\end{equation}
and the inequality
\begin{equation}
\frac{d}{dt} \overline{S} \leqslant \BL-\BR+\BI\,,
\label{eq:dsdtWithBIForBurg}
\end{equation}
where
\begin{equation}
\BR=\left.\left\{(F^{BC}-F^K)W^K+ F^{\ent,K}\right\}-\left\{(F_v^{BC}-F_v^K)W^K+W^{BC}F_v^K\right\}\right|_{\xi=1}\,,
\end{equation}

\begin{equation}
\BL=\left.\left\{(F^{BC}-F^1)W^1+ F^{\ent,1}\right\}-\left\{(F_v^{BC}-F_v^{1})W^1+W^{BC}F_v^{1}\right\}\right|_{\xi=-1}\,,
\end{equation}
where we replaced quantities at the boundaries by their boundary (BC) values and
\begin{equation}
\BI =  \sum\limits_\interiorfaces  \left\{F^{\ec,*}\jump{W}-\jump{F\,W}+ \jump{F^{\ent}}\right\}-\left\{F_v^*\jump{W}-\jump{F_v\,W}+W^*\jump{F_v}\right\}\,,
\end{equation} 
with the 1D jump defined in \cref{eq:jump1D}.

Regrouping the BL and BR terms gives
\begin{equation}
\BL-\BR  = \left.\left\{- F^\ent+W^{BC}\,F_v\right\}\right|_{0}^L - \left.(F^{BC}-F)W\right|_{0}^L + \left. (F_v^{BC}-F_v)W\right|_{0}^L,
\end{equation}
where we introduced the boundary notation $|_{0}^L$ to indicate that this is the evaluation at the left and right physical boundaries. The first term is analogous to the boundary term in the continuous problem, \cref{eq:mixed_visBurgers_entropy1}. The second and third terms are penalty terms that account for the weak boundary condition implementation at the physical boundaries. We note that the boundary condition implementation has to be dissipative for stability, i.e. boundary conditions must be specified so that 
\begin{equation}
\begin{split}
- \left.(F^{BC}-F)W\right|_{0}^L &\leqslant 0\,,\\
 \left. (F_v^{BC}-F_v)W\right|_{0}^L &\leqslant 0\,.
\end{split}
\label{eq:Burgers_physicalBC_dissipation}
\end{equation}
With dissipative boundary implementations, the contribution at the boundary can be estimated as
\begin{equation}
\BL-\BR  \leqslant \left.\left\{- F^\ent+W^{BC}\,F_v\right\}\right |_{0}^L\,,
\end{equation}
which matches the influence of the boundaries in the continuous problem, and which we also assume is implemented so that it is dissipative.

What is left to bound are the interface conditions, BI, in \cref{eq:dsdtWithBIForBurg}. The first part of BI includes contributions from the advective part of the PDE. Note that the specific numerical surface flux $F^{\ec,*}$ is defined so that the first term identically vanishes, c.f. \cref{eq:ec_condition_numflux}, thus all that remains is the viscous contribution,
\begin{equation}
\BI =  - \sum\limits_\interiorfaces F_v^*\jump{W}-\jump{F_v\,W}+W^*\jump{F_v}\,.
\end{equation} 
In the BR1 scheme, the numerical surface fluxes are the arithmetic means of the arguments. For the viscous Burgers equation, 
\begin{equation}
F_v^* = \average{F_v}\,,\quad
W^* = \average{W}\,.
\label{eq:BurgBr1}
\end{equation}
Using the identity for the jump of a product of two quantities
\begin{equation}
\jump{ab}=\average{a}\jump{b}+\jump{a}\average{b}
\end{equation}
 and inserting the numerical surface fluxes \cref{eq:BurgBr1} for the viscous terms, we see that the internal boundary contribution is
\begin{equation}
\begin{split}
\BI &=  - \sum\limits_\interiorfaces F_v^*\jump{W}-\jump{F_v\,W}+W^*\jump{F_v} \\
    &=  - \sum\limits_\interiorfaces F_v^*\jump{W}-\average{F_v}\,\jump{W}-\jump{F_v}\,\average{W}+W^*\jump{F_{v}} \\
    &=  - \sum\limits_\interiorfaces  \average{F_v}\jump{W}-\average{F_v}\,\jump{W}-\jump{F_v}\,\average{W}+\average{W}\jump{F_{v}} \\
    &= 0\,.
\end{split}
\end{equation}
Therefore, the temporal change of the total discrete entropy is 
\begin{equation}
\frac{d}{dt} \overline{S} \leqslant \left.\left\{- F^\ent+W^{BC}\,F_v\right\}\right |_{0}^L,
\end{equation}
which mimics the continuous estimate \cref{eq:mixed_visBurgers_entropy1}
\begin{equation}
\frac{d}{dt} \overline{s}\leqslant \left.\left\{ -f^\ent(u)+w(u)\,f_{v}\right\}\right|_0^L,
\end{equation}
except for the additional dissipation introduced at the physical boundaries, \cref{eq:Burgers_physicalBC_dissipation}. 

\begin{rem}
The entropy-stability of the approximation is not specific to the Burgers equation. It holds for any scalar advection-diffusion problem with entropy variables $w$ and entropy flux $f^{\ent}$. 
\end{rem}

\subsection{BR1 is Entropy-Stable: 3D Nonlinear Compressible Navier-Stokes Equations}

For the compressible Navier-Stokes equations \cref{eq:nse} 
\begin{equation}
{\statevec u_t} + \spacevec{\nabla}_x\cdot\bigstatevec{f} = \overRe\,\spacevec{\nabla}_x\cdot\bigstatevec{f}_v\left({\statevec u,{\spacevec\nabla _x}\statevec u}\right)\,,
\end{equation}
entropy-stability is unfortunately not equal to nonlinear stability, but it does give a stronger estimate than the linear stability of  \cref{sec:BR1_isstable_3Dlinearized_NSE} \cite{merriam1989,tadmor2003}.  Analogous to the nonlinear Burgers equation, \cref{sec:BR1_isstable_1D_Burger}, we introduce an entropy pair $(s,\spacevec{f}^\ent)$, this time with the scalar entropy function 
\begin{equation}
s = s(\statevec{u}) = -\frac{\rho \varsigma}{\gamma-1} =-\frac{\rho\left(\ln{p}-\gamma\ln{\rho}\right)}{\gamma-1}\,,
\end{equation}
where $\varsigma = \ln{p}-\gamma\ln{\rho}$ is the physical entropy, and the entropy flux
\begin{equation}
\spacevec{f}^\ent = s\,\spacevec{v}\,.
\end{equation}
The entropy variables are  
\begin{equation}
\statevec{w} = \frac{\partial s}{\partial\statevec{u}} = \left[
      \left(\frac{\gamma - \varsigma}{\gamma -1} - \frac{\rho|\spacevec{v}|^2}{2p}\right)\quad
      \left(\frac{\rho v_1}{p}\right)\quad 
      \left(\frac{\rho v_2}{p}\right)\quad
      \left(\frac{\rho v_3}{p}\right)\quad
      \left(\frac{-\rho}{p}\right)\right]^T,
\end{equation}
with the property
\begin{equation}
\statevec{k}^T\frac{\partial^2 s}{\partial\statevec{u}^2}\statevec{k} >0,\quad \forall \statevec{k}\neq 0,
\end{equation}
if $\rho>0$ and $p>0$ \cite{carpenter_esdg,tadmor2003,Dutt:1988}. The positivity requirement on the density and the temperature, $T\propto p/\rho$, guarantees a one-to-one mapping between conservative and entropy variables; it is this requirement that makes entropy-stability not a true nonlinear statement. Thus, entropy-stable discretizations can (and do) produce invalid solutions with negative density or temperature and need further strategies to guarantee positivity. 

The entropy pair contracts the Euler terms, meaning that it satisfies the relations
\begin{equation}
\statevec{w}^T\,\statevec{u}_t = \left(\frac{\partial s}{\partial\statevec{u}}\right)^T\statevec{u}_t  = s_t(\statevec{u}),
\label{eq:wsContraction}
\end{equation}
and
\begin{equation}
\statevec{w}^T\, \spacevec{\nabla}_x\cdot\bigstatevec{f} = \spacevec\nabla_x\cdot\spacevec{f}^\ent.
\end{equation}
Furthermore, we will see that the viscous flux should be rewritten in terms of the gradient of the entropy variables
\begin{equation}
\bigstatevec{f}_v\left({\statevec u,{\spacevec\nabla _x}\statevec u}\right) = \bigmatrix{B}^\ent\,\spacevec\nabla_x\,\statevec{w},
\end{equation}
where $\bigmatrix{B}^\ent$ satisfies the properties
\begin{equation}
\label{eq:spd_viscous_terms}
\mmatrix{B}^\ent_{ij} = (\mmatrix{B}^\ent_{ji})^T,\qquad \sum\limits_{i=1}^d\sum\limits_{j=1}^d\left(\frac{\partial \statevec w}{\partial x_i}\right)^T\,\mmatrix{B}^\ent_{ij}\,\left(\frac{\partial \statevec w}{\partial x_j}\right)\geqslant 0,\quad \forall \statevec w,
\end{equation}
if $p>0$ and $\mu>0$ \cite{carpenter_esdg,tadmor2006,Dutt:1988}.

\subsubsection{Continuous Entropy Analysis in 3D} 

To motivate the bounds needed for the approximation to be entropy-stable, we start with the continuous entropy analysis using the Navier-Stokes equations in the form
\begin{equation}
\begin{split}
{\statevec u_t} + \spacevec\nabla_x\cdot\bigstatevec{f} &= \overRe\,\spacevec\nabla_x\cdot\left(\bigmatrix{B}^\ent\,\spacevec\nabla_x\statevec{w}\right)=\overRe\,\spacevec\nabla_x\cdot\bigstatevec{f}_v ,\\
\bigstatevec{q} &=\spacevec\nabla_x\statevec{w}.
\end{split}
\end{equation}
We multiply the first equation with the entropy variables and the second equation with the viscous flux and integrate over the domain to get the weak forms
\begin{equation}
\begin{split}
\iprod{\statevec{w}(u), {\statevec u_t}} + \iprod{\statevec{w}(u),\spacevec\nabla_x\cdot\bigstatevec{f}} &= \overRe\,\iprod{\statevec{w}(u),\spacevec\nabla_x\cdot\bigstatevec{f}_v} ,\\
\iprod{\bigstatevec{q}, \bigstatevec{f}_v}&=\iprod{\spacevec\nabla_x\statevec{w},\bigstatevec{f}_v}.
\end{split}
\label{eq:weakformNS1}
\end{equation}
Next we use the properties of the entropy pair to contract the left side of the first equation of \cref{eq:weakformNS1} and use integration by parts on the right hand side to get
\begin{equation}
\begin{split}
\iprod{s_t(\statevec{u}),1} + \iprod{\spacevec\nabla_x\cdot\spacevec{f}^\ent,1} &= \overRe\,\int\limits_{\partial\Omega}\statevec{w}^T(u)\,\left(\bigstatevec{f}_v\cdot\spacevec{n}\right) \dS  - \overRe\,\iprod{\spacevec\nabla_x\statevec{w}(u),\bigstatevec{f}_v} ,\\
\iprod{\bigstatevec{q}, \bigstatevec{f}_v}&=\iprod{\spacevec\nabla_x\statevec{w},\bigstatevec{f}_v}.
\end{split}
\label{eq:weakformNS2}
\end{equation}
Inserting the second equation of \cref{eq:weakformNS2} into the first gives
\begin{equation}
\begin{split}
\iprod{s_t(\statevec{u}),1} + \iprod{\spacevec\nabla_x\cdot\spacevec{f}^\ent,1} &= \overRe\,\int\limits_{\partial\Omega}\statevec{w}^T(u)\,\left(\bigstatevec{f}_v\cdot\spacevec{n}\right) \dS  - \overRe\,\iprod{\bigstatevec{q}, \bigmatrix{B}^\ent\,\bigstatevec{q}}.
\end{split}
\end{equation}
We use the property \cref{eq:spd_viscous_terms} and integrate the flux divergence on the left side to get the estimate
\begin{equation}
\label{eq:continuous_NSE_entropy_estimate}
\begin{split}
\frac{d}{dt}\overline{s}  \leqslant \int\limits_{\partial\Omega}\left( -\spacevec{f}^\ent \cdot\spacevec{n} + \overRe\,\statevec{w}^T(u)\left(\bigstatevec{f}_v\cdot\spacevec{n}\right)\right)\dS,
\end{split}
\end{equation}
where 
\begin{equation}
\overline{s} = \iprod{s(\statevec{u}),1}=\int\limits_{\Omega} s(\statevec{u}) \,\mathrm{dV}
\end{equation}
is again the total entropy. Boundary conditions then need to be specified so that the bound on the entropy depends only on the boundary data. We will assume here that boundary data is given so that the right hand side is non-positive so that the entropy will not increase in time. For more on boundary conditions for the Navier-Stokes equations, see, e.g. \cite{Nordstrom:2005qy}.

\subsubsection{An Entropy-Stable DGSEM in 3D}

As in  \cref{sec:ESDGEM_NSE3D}, we transform the Navier-Stokes equations into reference space
\begin{equation}
\begin{split}
\mathcal J\,\statevec u_t + \spacevec\nabla_\xi\cdot\bigcontravec{f} &= \overRe\,\spacevec\nabla_\xi\cdot\bigcontravec{f}_v,\\
\mathcal J\,\bigstatevec{q} &=\bigmatrix{M}\,\spacevec\nabla_\xi\statevec{w},
\end{split}
\end{equation}
where
\begin{equation}
\bigcontravec{f} = \bigmatrix{M}^T\,\bigstatevec{f},\quad \bigcontravec{f}_v = \bigmatrix{M}^T\,\bigstatevec{f}_v,
\end{equation}
or alternatively 
\begin{equation}
\begin{split}
\mathcal J\,\statevec u_t + \spacevec\nabla_\xi\cdot\left(\bigmatrix{M}^T\,\bigstatevec{f}\right) &= \overRe\,\spacevec\nabla_\xi\cdot\left(\frac{1}{\mathcal J}\bigmatrix{M}^T\,\bigmatrix{B}^\ent\,\bigmatrix{M}\,\spacevec\nabla_\xi\statevec{w}\right).
\end{split}
\label{eq:AltViscousfluxform}
\end{equation}
The standard strong form of the DGSEM is derived as described in \cref{sec:standard_DGSEM} and is
\begin{equation}
\begin{split}
\iprod{J\,\statevec U_t,\boldsymbol{\phi}}_N + &\iprod{\spacevec\nabla_\xi\cdot\mathbb{I}^N\left(\bigcontravec{F}\right),\boldsymbol{\phi}}_N +\int\limits_{\partial E,N} \,\boldsymbol{\phi}^T\left(\contrastatevec{F}^*-\left(\bigcontravec{F}\cdot\hat{n}\right)\right)\,\dS\\
&= \overRe\int\limits_{\partial E,N} \boldsymbol{\phi}^T\left(\contrastatevec{F}_{v}^{*}-\left(\bigcontravec{F}_v\cdot\hat{n}\right)\right)\,\dS + \overRe\,\iprod{\spacevec\nabla_\xi\cdot\bigcontravec{F}_v,\boldsymbol{\phi}}_N \\
\iprod{J\,\bigstatevec{Q},\bigstatevec{\boldsymbol{\psi}}}_N &= \int\limits_{\partial E,N} \left(\statevec{W}^*\right)^T\bigmatrix{M}^T\left(\bigstatevec{\boldsymbol{\psi}}\cdot\hat{n}\right)\,\dS - \iprod{\statevec{W},\spacevec\nabla_\xi\cdot\mathbb{I}^N\left(\bigmatrix{M}^T\,\bigstatevec{\boldsymbol{\psi}}\right)}_N\,,
\end{split}
\label{eq:OrigNonlinNSDGSEM}
\end{equation}
where $\statevec U,\statevec W,\statevec Q$ are the polynomial interpolants of the conservative variables, the entropy variables and the gradient of the entropy variables with nodes at the Gauss-Lobatto points. Furthermore, we use the polynomial interpolations of the contravariant Euler and the viscous flux $\bigcontravec{F}, \bigcontravec{F}_v$. As always, the quantities in the discrete surface integrals marked with $*$ are the numerical fluxes and entropy variables that couple the elements and depend on values (and gradients) from the ``left'' and ``right'' of the interface. 

Similar to the DGSEM discretization of the Burgers equation in \cref{sec:BR1_isstable_1D_Burger}, the standard DGSEM discretization of the Euler terms in \cref{eq:OrigNonlinNSDGSEM} is not entropy-stable and is infected with aliasing errors arising from the nonlinearity of the Euler fluxes. In contrast to the Burgers equation where a two point flux \cref{eq:BurgersTwoPointFlux} is equivalent to the explicit split form seen in \cref{eq:splitFormBurgEquivalence}, an explicit split form of the nonlinear Euler flux divergence that gives discrete entropy-stability is not known. Again, we use the framework of Carpenter and Fisher et al., e.g. \cite{skew_sbp2,carpenter_esdg}, that allows us to get a discrete entropy-stable approximation of general nonlinear hyperbolic conservation laws even when the explicit split form is not known, as long as a two-point entropy conservative numerical flux function is known. 

For an entropy-stable approximation to the advective terms, we approximate the divergence of the flux with an entropy conservative flux
\begin{equation}
\label{eq:entropy_stable_Euler_vol}
\spacevec\nabla_\xi\cdot\mathbb{I}^{N}\left(\bigcontravec{F}(\xi,\eta,\zeta)\right)\approx \spacevec{\mathbb{D}} (\bigcontravec{F} )^{\ec}(\xi,\eta,\zeta),
\end{equation}
where
\begin{equation}
\begin{split}
\spacevec{\mathbb{D}} (\bigcontravec{F} )^{\ec}(\xi,\eta,\zeta)\equiv2\sum_{m=0}^N 
&\quad \ell'_m(\xi)   \contrastatevec{F}^{1,\ec}(\xi,\eta,\zeta;\xi_m,\eta,\zeta)\\[-1ex]       
&+     \ell'_m(\eta)  \contrastatevec{F}^{2,\ec}(\xi,\eta,\zeta;\xi,\eta_m,\zeta)\\[1ex]    
&+     \ell'_m(\zeta) \contrastatevec{F}^{3,\ec}(\xi,\eta,\zeta;\xi,\eta,\zeta_m),                                                                                                                                              
\end{split}
\label{eq:DFEC}
\end{equation}
and the contravariant numerical volume fluxes are
\begin{equation}
\begin{split}
&\contrastatevec{F}^{l,\ec}(\xi,\eta,\zeta; \alpha,\beta,\gamma)\equiv\bigstatevec{F}^{\ec}(U(\xi,\eta,\zeta),U(\alpha,\beta,\gamma))
\cdot \frac{1}{2}\left[J\spacevec{a}^{\,l}(\xi,\eta,\zeta) + J\spacevec{a}^{\,l}(\alpha,\beta,\gamma)\right],\quad l=1,2,3.
\end{split}
\label{eq:ECFluxWithMetrics}
\end{equation}
An example of an entropy conserving two point flux is described in \cref{app: Euler_EC_Flux}. Note that since $\bigstatevec F^{\ec}$ is symmetric in its arguments, so is $\bigcontravec F ^{\ec}$.

If we make the replacement \cref{eq:entropy_stable_Euler_vol} into the discrete volume integral for the advective flux in \cref{eq:OrigNonlinNSDGSEM} and replace the test function with the Lagrange basis,
\begin{equation}
\begin{split}
\iprod{\spacevec{\mathbb{D}} (\bigcontravec{F} )^{\ec}\!\!\!,\ell_i(\xi)\ell_j(\eta)\ell_k(\zeta)}_N&=
\sum\limits_{a,b,c=0}^N\omega_a\omega_b\omega_c\ell_i(\xi_a)\ell_j(\eta_b)\ell_k(\zeta_c)\,
  \spacevec{\mathbb{D}} (\bigcontravec{F} )^{\ec}(\xi_a,\eta_b,\zeta_c)\\     
&=\omega_i\omega_j\omega_k\,\spacevec{\mathbb{D}} (\bigcontravec{F} )^{\ec}(\xi_i,\eta_j,\zeta_k) 
 =\omega_i\omega_j\omega_k\,\spacevec{\mathbb{D}} (\bigcontravec{F} )^{\ec}_{ijk},
\end{split}
\end{equation}
where
\begin{equation}
\label{eq:entropy-cons_volint}
\begin{split}
\spacevec{\mathbb{D}} (\bigcontravec{F} )^{\ec}_{ijk}\equiv
&\quad 2\sum_{m=0}^N D_{im}\,\left(\bigstatevec{F}^{\ec}(U_{ijk}, U_{mjk})\cdot\average{J\spacevec{a}^{\,1}}_{(i,m)jk}\right)\\[0.1cm]
&+     2\sum_{m=0}^N D_{jm}\,\left(\bigstatevec{F}^{\ec}(U_{ijk}, U_{imk})\cdot\average{J\spacevec{a}^{\,2}}_{i(j,m)k}\right)\\[0.1cm]
&+     2\sum_{m=0}^N D_{km}\,\left(\bigstatevec{F}^{\ec}(U_{ijk}, U_{ijm})\cdot\average{J\spacevec{a}^{\,3}}_{ij(k,m)}\right),
\end{split}
\end{equation}
is exactly the discrete form used in Carpenter and Fisher et al., e.g. \cite{skew_sbp2,carpenter_esdg} and Gassner et al. \cite{Gassner:2016ye} with the averaging of the metric terms, e.g.
\begin{equation}
\average{J\spacevec{a}^{\,1}}_{(i,m)jk}\equiv \frac{1}{2}\left\{J\spacevec{a}^{\,1}_{ijk}+J\spacevec{a}^{\,1}_{mjk}\right\}.
\end{equation}

\begin{rem}
In \cite{Gassner:2016ye} it was shown that standard (strong form) DGSEM volume integral $\iprod{\spacevec\nabla_\xi\cdot\bigcontravec{F},\phi}_N$ is recovered if we replace $\contrastatevec{F}^{l,\ec},l=1,2,3$ in \cref{eq:DFEC} with 
\begin{equation}
\contrastatevec{F}^{l,\stand}(\xi,\eta,\zeta; \alpha,\beta,\gamma)\equiv
\frac{1}{2}\left[\bigstatevec{F}(U(\xi,\eta,\zeta))     \cdot J\spacevec{a}^{\,l}(\xi,\eta,\zeta) 
                +\bigstatevec{F}(U(\alpha,\beta,\gamma))\cdot J\spacevec{a}^{\,l}(\alpha,\beta,\gamma)\right]\,.
\end{equation}
Note that the standard discretization differs in that it has the average of the flux times the metric terms, whereas in the entropy-stable variant \cref{eq:ECFluxWithMetrics} has the specially averaged entropy-conservative flux \emph{times} the average of the metric terms. The separation of the average of the metric terms corresponds, c.f , e.g. \cite{Kopriva:2014yq}, to a de-aliasing of the metric terms (which are variable functions themselves when elements have curved sides) and affects stability, see \cref{app:Proof__volint_entropy_Euler}.
\end{rem}

The goal now is to show that the DGSEM with BR1 is entropy-stable when used in combination with an entropy-stable discretization of the advective (Euler) terms. To show stability, we make the replacement \cref{eq:entropy_stable_Euler_vol} into \cref{eq:OrigNonlinNSDGSEM} and assume the use of entropy conservative numerical fluxes that satisfy the equivalent of \cref{eq:ec_condition_numflux} at the element boundaries. The final discretization is therefore
\begin{equation}
\label{eq:disc_nse1}
\begin{split}
\iprod{J\,\statevec U_t,\boldsymbol{\phi}}_N + &\iprod{\spacevec{\mathbb{D}} (\bigcontravec{F} )^{\ec},\boldsymbol{\phi}}_N +\int\limits_{\partial E,N} \boldsymbol{\phi}^T\left(\contrastatevec{F}^{\ec,*}-\left(\bigcontravec{F}\cdot\hat{n}\right)\right)\,\dS \\
&= \overRe\int\limits_{\partial E,N} \boldsymbol{\phi}^T\left(\contrastatevec{F}_{v}^{*}-\left(\bigcontravec{F}_v\cdot\hat{n}\right)\right)\,\dS + \overRe\,\iprod{\spacevec\nabla_\xi\cdot\bigcontravec{F}_v,\boldsymbol{\phi}}_N \\
\iprod{J\,\bigstatevec{Q},\bigstatevec{\boldsymbol{\psi}}}_N &= \int\limits_{\partial E,N} (\statevec{W}^*)^T\bigmatrix{M}^T\left(\bigstatevec{\boldsymbol{\psi}}\cdot\hat{n}\right)\,\dS - \iprod{\statevec{W},\spacevec\nabla_\xi\cdot\mathbb{I}^N\left(\bigmatrix{M}^T\,\bigstatevec{\boldsymbol{\psi}}\right)}_N\,.
\end{split}
\end{equation}

We now mimic the continuous entropy analysis as closely as possible and replace $\boldsymbol{\phi} \leftarrow \statevec W$ and $\bigstatevec{\boldsymbol{\psi}} \leftarrow \bigstatevec{F}_v$ in \cref{eq:disc_nse1} to get
\begin{equation}
\begin{split}
\iprod{J\,\statevec U_t,\statevec{W}}_N + &\iprod{\spacevec{\mathbb{D}} (\bigcontravec{F} )^{\ec},\statevec{W}}_N +\int\limits_{\partial E,N} \statevec{W}^T\left(\contrastatevec{F}^{\ec,*}-\left(\bigcontravec{F}\cdot\hat{n}\right)\right)\,\dS\\
&= \overRe\int\limits_{\partial E,N} \statevec{W}^T\left(\contrastatevec{F}_{v}^{*}-\left(\bigcontravec{F}_v\cdot\hat{n}\right)\right)\,\dS + \overRe\,\iprod{\spacevec\nabla_\xi\cdot\bigcontravec{F}_v,\statevec{W}}_N \\
\iprod{J\,\bigstatevec{Q},\bigstatevec{F}_v}_N &= \int\limits_{\partial E,N} (\statevec{W}^*)^T\left(\bigcontravec{F}_v\cdot\hat{n}\right)\,\dS - \iprod{\statevec{W},\spacevec\nabla_\xi\cdot\bigcontravec{F}_v}_N\,,
\end{split}
\label{eq:startOfNSEntropyAnalysis}
\end{equation}
where 
$
\mathbb{I}^N\left(\bigmatrix{M}^T\,\bigstatevec{F}_v\right) = \bigcontravec{F}_v\,.
$

We look first at the first term of the second equation of \cref{eq:startOfNSEntropyAnalysis} and insert the alternate form of the viscous flux shown in \cref{eq:AltViscousfluxform}. We use the property \cref{eq:spd_viscous_terms} of the viscous flux matrices $\bigmatrix{B}^\ent$ to see that 
\begin{equation}
\label{eq:disc_visc_volint_estimate}
\iprod{J\,\bigstatevec{Q},\bigstatevec{F}_v}_N = \iprod{J\,\bigstatevec{Q},\bigmatrix{B}^\ent\,\bigstatevec{Q}}_N\geqslant \min\limits_{E,N}(J) \iprod{\bigstatevec{Q},\bigmatrix{B}^\ent\,\bigstatevec{Q}}_N\geqslant 0\,,
\end{equation}
as long as the interpolant of the element mapping Jacobian determinant, $J$, is non-negative at the Gauss-Lobatto nodes.

We next replace the volume integral for the advective flux by a surface integral 
\begin{equation}
\label{eq:disc_entropy_cons_Euler}
\iprod{\spacevec{\mathbb{D}}(\bigcontravec{F})^{\ec},\statevec{W}}_N =  \int\limits_{\partial E,N} (\contraspacevec{F}^\ent\cdot\hat{n})\,\dS,
\end{equation}
which we prove in \cref{app:Proof__volint_entropy_Euler}.

Finally, we insert the second equation of \cref{eq:startOfNSEntropyAnalysis} into the first and use the estimate \cref{eq:disc_visc_volint_estimate} to allow all but the time derivative term to be replaced with surface integrals and an inequality
\begin{equation}
\begin{split}
\iprod{J\,\statevec U_t,\statevec{W}}_N +&\int\limits_{\partial E,N} \left\{(\statevec{W})^T\left(\contrastatevec{F}^{\ec,*}-\left(\bigcontravec{F}\cdot\hat{n}\right)\right) + \left(\contraspacevec{F}^\ent\cdot\hat{n}\right)\right\}\dS\\
\leqslant \overRe&\int\limits_{\partial E,N} \left\{\statevec{W}^T\left(\contrastatevec{F}_{v}^{*}-\left(\bigcontravec{F}_v\cdot\hat{n}\right)\right)+\left(\statevec{W}^*\right)^T\,\left(\bigcontravec{F}_v\cdot\hat{n}\right)\right\}\dS.
\end{split}
\end{equation}

The time derivative term is the time rate of change of the entropy in the element. Assuming that the chain rule with respect to differentiation in time holds, we use the contraction property of the entropy variable \cref{eq:wsContraction} at each GL point within the element to see that
\begin{equation}
\iprod{J\,\statevec U_t,\statevec{W}}_N = \sum\limits_{i,j,k=0}^N \omega_{ijk}J_{ijk}\statevec{W}^{T}_{ijk}\frac{d \statevec{U}_{ijk}}{dt}\, = \sum\limits_{i,j,k=0}^N \omega_{ijk}J_{ijk}\,\frac{d S_{ijk}}{dt} = \iprod{J\,S_t,1}_N.
\end{equation}
 We get the total discrete entropy by summing over all elements
\begin{equation}
\frac{d}{dt}\overline{S} \equiv \sum\limits_{k=1}^K  \iprod{J^k\,S^k_t,1}_N,
\end{equation}

so the total discrete entropy satisfies the estimate
\begin{equation}
\begin{split}
\frac{d}{dt}\overline{S}\leqslant &- \sum\limits_\interiorfaces  \int\limits_{N} \left\{\left(\contrastatevec{F}_{n}^{\ec,*}\right)^T\jump{\statevec W} - \jump{\left(\contrastatevec{F}_{n}\right)^T\,\statevec{W}} + \jump{\contravec{F}^\ent_{n}}\right\}\dS\\
                                                &+\overRe\sum\limits_\interiorfaces \int\limits_{N} \left\{\left(\contrastatevec{F}_{v,n}^{*}\right)^T\jump{\statevec{W}} - \jump{(\contrastatevec{F}_{v,n})^T\,\statevec{W}}+\left(\statevec{W}^*\right)^T\,\jump{\contrastatevec{F}_{v,n}}\right\}\dS +\PBT\,,
\end{split}
\end{equation}
where $PBT$ are the physical boundary terms and where we use the \emph{slave--master} definition of the jump \cref{eq:jump_master_slave}, with all interior faces in master element side orientation. To ensure stability, the numerical surface flux of the Euler terms, $\contrastatevec{F}_n^{\ec,*}$ is chosen so that at each quadrature point on the interior element faces the argument of the advective surface quadratures vanishes, i.e.
\begin{equation}
\label{eq:ec_surfint_euler}
\left(\contrastatevec{F}_{n}^{\ec,*}\right)^T\jump{\statevec{W}} - \jump{\contrastatevec{F}_{n}^T\,\statevec{W}} + \jump{\contravec{F}^\ent_{n}}= 0\,.
\end{equation}
Compare this generalized condition to the Burgers case \cref{eq:ec_condition_numflux} and see \cref{app: Euler_EC_Flux} for a choice for the entropy conserving numerical surface flux for the compressible Euler equations.  
See also \cref{app:entropy-conservation_3D_curvilinear_surface} to see how to relate the Cartesian flux conditions to the contravariant flux condition. 

The choice of the numerical surface fluxes for the viscous terms for the BR1 discretization are again simply the averages
\begin{equation}
\contrastatevec{F}_{v,n}^{*}=\average{\contrastatevec{F}_{v,n}}\,,\quad \statevec{W}^* = \average{\statevec{W}}\,.
\end{equation}
With the identity (see \cref{eq:KnightsFormerlyKnownAsLemma1})
\begin{equation}
 \jump{(\contrastatevec{F}_{v,n})^T\,\statevec{W}} = \average{(\contrastatevec{F}_{v,n})^T}\jump{\statevec{W}} + \jump{(\contrastatevec{F}_{v,n})^T}\average{\statevec{W}},
\end{equation}
we see that the contribution of the viscous numerical fluxes at the interior faces vanishes exactly
\begin{equation}
\overRe\sum\limits_\interiorfaces \int\limits_{N} \left\{\left(\contrastatevec{F}_{v,n}^{*}\right)^T\jump{\statevec{W}} - \jump{\left(\contrastatevec{F}_{v,n}\right)^T\,\statevec{W}}+\left(\statevec{W}^*\right)^T\,\jump{\contrastatevec{F}_{v,n}}\right\}\dS = 0\,.
\end{equation}

The final discrete entropy statement therefore 
\begin{equation}
\begin{split}
\frac{d}{dt}\overline{S}\leqslant
&\quad \sum\limits_\boundaryfaces \int\limits_{N} - \contravec{F}_n^\ent + \overRe \statevec{W}^T\, \contrastatevec{F}_{v,n}\,\dS\\
&+     \overRe\sum\limits_\boundaryfaces \int\limits_{N}  -\left(\contrastatevec{F}_{n}^{\ec,*}\right)^T\statevec{W} + (\contrastatevec{F}_{v,n}^{*})^T\statevec{W} +  \left(\statevec{W}^{*}\right)^T\,\left(\contrastatevec{F}_{v,n}\right)\,\dS\, ,
\end{split}
\end{equation}
which precisely mimics the continuous estimate \cref{eq:continuous_NSE_entropy_estimate} except for the additional dissipation at the physical boundaries due to the numerical surface fluxes $\bigcontravec{F}_{n}^{\ec,*},\bigcontravec{F}_{v,n}^{*},\statevec{W}^{*}$ evaluated at the boundary. Note that to guarantee stability, the choice of these auxiliary physical boundary terms must be determined to ensure that their effect is dissipative. From another point of view, the additional term gives constraints on the boundary fluxes from which to derive stable boundary conditions.

\section{Conclusions}
\label{sec:conclusions}
In this work stability of the BR1 scheme for general advection diffusion systems is shown. The focus of the paper is on the compressible NSE, however due to the generic structure of the BR1 scheme (i.e. it does not depend on the particular form of the viscous PDE terms) the stability estimate extends to all advection-diffusion problems, if stable DGSEM approximations of the advection terms are known. 

The paper is split in two parts, the linear energy analysis and the nonlinear entropy analysis. Whereas it is sufficient in the linear case to consider skew-symmetric split form approximations of the advection terms, special entropy-stable DGSEM formulations are necessary for the nonlinear case. To the best knowledge of the authors, we also present for the first time the full proof of entropy-stability for a DGSEM of the compressible Euler terms on three dimensional curvilinear meshes.

\textblue{Furthermore, we show that a Bassi and Rebay type approximation can be provably stable if the metric identities are discretely satisfied, a two-point average for the metric terms is used in the volume integral, the auxiliary gradients for the viscous terms are computed from gradients of entropy variables, and the BR1 scheme is used for the interface fluxes.}  

\textblue{Our analysis shows that even for three dimensional curvilinear grids, the BR1 fluxes do not add artificial dissipation at the interior element faces. Thus, the BR1 interface fluxes preserve the stability of the discretization of the advection terms and we get either energy stability or entropy-stability for the linear or nonlinear compressible NSE, respectively.}

\textred{An open source code that implements the scheme is available at github.com/project-fluxo. Numerical investigations show that it is robust in comparison to the standard DGSEM (even in comparison to ad hoc stabilisation techniques such as polynomial de-aliasing) for under-resolved compressible turbulence simulations, e.g. \cite{Gassner:2016ye}.} 

\textred{To further highlight the remarkable robustness of the stable DGSEM we conclude with a discussion of a numerical simulation of the Taylor-Green vortex with a Reynolds number of $Re=1600$ (see e.g. \cite{Gassner:2013qf,Gassner:2016ye} for details of the setup). We intentionally stress the numerical experiment by choosing the scheme \textit{without} added artificial stabilisation terms. That is, we use the BR1 scheme for the viscous terms and we use only the EC flux at the interface for the Euler terms, without adding any Riemann solver type dissipation. To further stress our high order method and to highlight its enhanced robustness, we intentionally under-resolve the physics of the vortex dominated flow.} 

\textred{A common challenging benchmark resolution is $64^3$ degrees of freedom. We choose a high polynomial degree of $N=7$ (eight GLL nodes in each direction inside the element) and $8^3$ elements. If we use the standard DGSEM without the entropy-stable modifications described here the simulation immediately crashes,  even when using Riemann solvers with added dissipation at the interface, see e.g. \cite{Gassner:2013qf}.  However, when using the entropy-stable approximation without added numerical dissipation at the element interfaces the simulation no longer crashes.} 

\begin{figure}
\centerline{\includegraphics{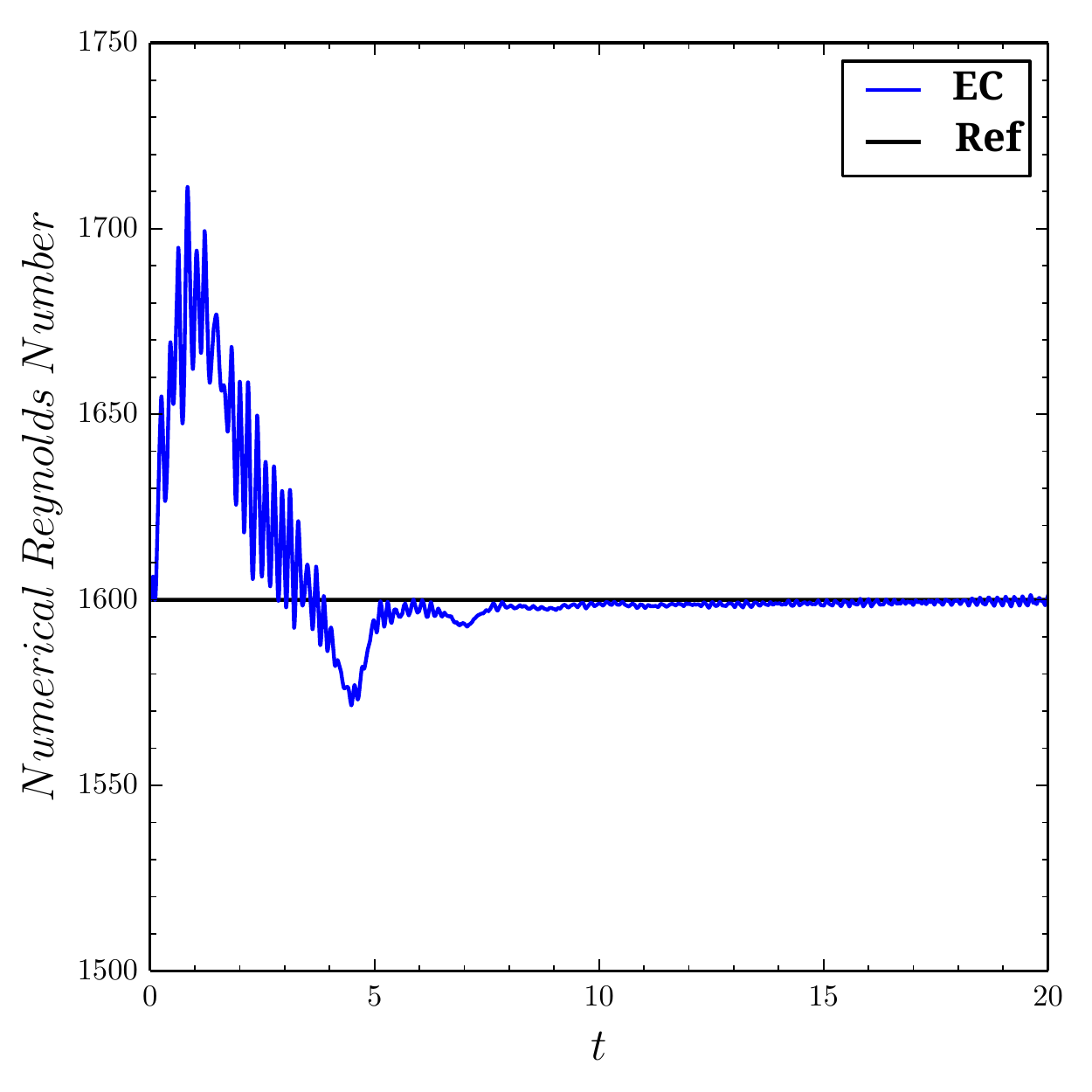}}
\caption{\label{fig:numericalReynoldsnumber}Numerical simulation of the Taylor-Green vortex with Mach number $M=0.1$ and Reynolds number $Re=1600$ showing temporal evolution of the numerical Reynolds number in comparison to the physical value. The DG discretization uses $8^3$ elements and a polynomial degree of $N=7$ resulting in $64^3$ degrees of freedom. The BR1 method is used for the viscous terms and the EC flux without additional dissipation is used for the Euler fluxes at the interfaces, making the scheme virtually dissipation free yet stable. The standard DGSEM immediately crashes for this test case.}
\end{figure}

\textred{In addition, the numerical scheme is virtually dissipation free. To see this,  we choose a Mach number of $M=0.1$ for the Taylor-Green vortex to allow for comparisons with the incompressible NSE. For incompressible flows, it is possible to relate the enstrophy ($ens$), the kinetic energy dissipation rate ($diss$) and the Reynolds number. Using this relationship, we estimate the numerical Reynolds number by computing the ratio of the discrete total enstrophy and the discrete total kinetic energy dissipation rate, i.e. $Re_{num}(t) \approx 2\,\frac{ens(t)}{diss(t)}$. The result is presented in Fig. \ref{fig:numericalReynoldsnumber}. As a reference, the physical Reynolds number $Re=1600$ is plotted as well. It is remarkable that for the fully developed flow at later times the numerical Reynolds number approaches the physical Reynolds number, which confirms the two properties proven here: (i) the method is virtually dissipation free; (ii) it is still stable when the standard scheme crashes.}

\acknowledgement{ This work was supported by a grant from the Simons Foundation (\#426393, David Kopriva). G.G has been supported by the European Research Council (ERC) under the European Union's Eights Framework Program Horizon 2020 with the research project \textit{Extreme}, ERC grant agreement no. 714487.

\appendix
\normalsize
\section{An Entropy Conserving Euler Flux}
\label{app: Euler_EC_Flux}

We explicitly present the Kinetic Energy Preserving and Entropy Conservative (KEPEC) numerical flux function for the compressible Euler equations. The KEPEC flux was first derived by Chandrashekar \cite{Chandrashekar2012}
\begin{equation}\label{Eq:EKEP}
\statevec{F}_1^{*,\rm ec} = {\small\begin{bmatrix}
\rho^{\ln}\avg{v_1} \\[0.15cm]
\rho^{\ln}\avg{v_1}^2 + \avg{p} \\[0.15cm]
\rho^{\ln}\avg{v_1}\avg{v_2}  \\[0.15cm]
\rho^{\ln}\avg{v_1}\avg{v_3} \\[0.15cm]
\frac{p^{\ln}\avg{v_1}}{(\gamma-1)} +\avg{p}\avg{v_1} + \frac{1}{2}\rho^{\ln}\avg{v_1}\uavg
\end{bmatrix}},
\end{equation}
where
\begin{equation}
\avg{p} = \frac{\avg{\rho}}{2\avg{\beta}},\qquad{p}^{\ln} = \frac{{\rho}^{\ln}}{2{\beta}^{\ln}},\qquad \uavg = 2\left(\avg{v_1}^2 + \avg{v_2}^2 + \avg{v_3}^2\right)-\left(\avg{v_1^2} + \avg{v_2^2} + \avg{v_3^2}\right),
\end{equation}
and $(\cdot)^{\ln}$ is the logarithmic mean
\begin{equation}
	(\cdot)\textsuperscript{ln} = \frac{\jump{\cdot}}{\jump{\ln(\cdot)}}\,. %
\end{equation}
A numerically stable method to compute the logarithmic mean when $(\cdot)_{\rm R} \approx (\cdot)_{\rm L}$ is given in \cite[App.~B]{IsmailRoe2009}.

It is also possible to add additional dissipation to the entropy conservative scheme and create an entropy-stable (ES) approximation. One example is with a matrix dissipation entropy-stable flux for the compressible Euler equations of the form
\begin{equation}
\statevec{F}^{*,\rm ec} = \statevec{F}^{*,\rm ec} -  \frac{1}{2}\R|\boldsymbol{\hat{\Lambda}}|\T\R^T\jump{\statevec{w}}.
\end{equation}
The average components of the dissipation term are given by
\begin{equation}\label{eq:EulerDissipation}
\begin{aligned}
\R &= {\small\begin{bmatrix}
1 & 1 & 0 & 0 & 1 \\[0.1cm]
\avg{v_1} - \bar{a} & \avg{v_1} & 0 & 0 & \avg{v_1} + \bar{a} \\[0.1cm]
\avg{v_2}  & \avg{v_2} & 1 & 0 & \avg{v_2} \\[0.1cm]
\avg{v_3} & \avg{v_3} & 0 & 1 & \avg{v_3} \\[0.1cm]
\bar{h} - \avg{v_1}\bar{a} & \frac{1}{2}\overline{\|\textbf{v}\|^2} & \avg{v_2} & \avg{v_3} & \bar{h} + \avg{v_1}\bar{a} \\[0.1cm]
\end{bmatrix}},\\[0.2cm]
\boldsymbol{\hat{\Lambda}} &= \textrm{diag}\left(\avg{v_1} - \bar{a},\avg{v_1},\avg{v_1},\avg{v_1},\avg{v_1} + \bar{a}\right),\\[0.2cm]
\T &= \textrm{diag}\left(\frac{\rho^{\ln}}{2\gamma},\frac{\rho^{\ln}(\gamma-1)}{\gamma},\avg{p},\avg{p},\frac{\rho^{\ln}}{2\gamma}\right),
\end{aligned}
\end{equation}
where
\begin{equation}
\bar{a} = \sqrt{\frac{\gamma \avg{p}}{\rho^{\ln}}},\qquad \bar{h} = \frac{\gamma}{2\beta^{\ln}(\gamma-1)} + \frac{1}{2}\uavg\,.
\end{equation}
We note that the selection of the discrete dissipation operator \cref{eq:EulerDissipation} creates a scheme that is able to exactly resolve stationary contact discontinuities. The proof of this property follows the same structure as that presented by Chandrashekar \cite{Chandrashekar2012}.

For completeness, we provide the entropy conserving (and stable) fluxes for the other Cartesian directions 
\begin{equation}\label{Eq:EKEP2}
\statevec{F}_2^{*,\rm ec} = {\small\begin{bmatrix}
\rho^{\ln}\avg{v_2} \\[0.15cm]
\rho^{\ln}\avg{v_1}\avg{v_2}  \\[0.15cm]
\rho^{\ln}\avg{v_2}^2 + \avg{p}  \\[0.15cm]
\rho^{\ln}\avg{v_2}\avg{v_3} \\[0.15cm]
\frac{p^{\ln}\avg{v_2}}{(\gamma-1)} +\avg{p}\avg{v_2} + \frac{1}{2}\rho^{\ln}\avg{v_2}\uavg
\end{bmatrix}},
\;
\statevec{F}_3^{*,\rm ec} = {\small\begin{bmatrix}
\rho^{\ln}\avg{v_3} \\[0.15cm]
\rho^{\ln}\avg{v_1}\avg{v_3}  \\[0.15cm]
\rho^{\ln}\avg{v_2}\avg{v_3}  \\[0.15cm]
\rho^{\ln}\avg{v_3}^2 + \avg{p} \\[0.15cm]
\frac{p^{\ln}\avg{v_3}}{(\gamma-1)} +\avg{p}\avg{v_3} + \frac{1}{2}\rho^{\ln}\avg{v_3}\uavg
\end{bmatrix}}.
\end{equation}
The matrix dissipation term remains similar where the matrix of right eigenvectors are given by
\begin{equation}\label{eq:EulerDissipation2}
\bm{\hat{\mathcal{R}}_2} = {\small\begin{bmatrix}
1 & 0 & 1 & 0 & 1 \\[0.1cm]
\avg{v_1}  & 1 & \avg{v_1} & 0 & \avg{v_1} \\[0.1cm]
\avg{v_2} - \bar{a} & 0 & \avg{v_2} & 0 & \avg{v_2} + \bar{a} \\[0.1cm]
\avg{v_3} & 0 & \avg{v_3} & 1 & \avg{v_3} \\[0.1cm]
\bar{h} - \avg{v_2}\bar{a} & \avg{v_1} & \frac{1}{2}\overline{\|\textbf{v}\|^2} & \avg{v_3} & \bar{h} + \avg{v_2}\bar{a} \\[0.1cm]
\end{bmatrix}} ,
\;
\bm{\hat{\mathcal{R}}_3} = {\small\begin{bmatrix}
1 & 0 & 0 & 1 & 1 \\[0.1cm]
\avg{v_1}  & 1 & 0 &\avg{v_1} & \avg{v_1} \\[0.1cm]
\avg{v_2}  & 0 & 1 & \avg{v_2} & \avg{v_2}  \\[0.1cm]
\avg{v_3} - \bar{a} & 0 & 0 & \avg{v_3} & \avg{v_3} + \bar{a} \\[0.1cm]
\bar{h} - \avg{v_3}\bar{a} & \avg{v_1} & \avg{v_2} & \frac{1}{2}\overline{\|\textbf{v}\|^2} & \bar{h} + \avg{v_3}\bar{a} \\[0.1cm]
\end{bmatrix}} 
\end{equation}
respectively. The diagonal matrix of eigenvalues uses the appropriate value of the velocity depending on the spatial direction. 

To create the contravariant entropy conservative fluxes we incorporate the average of the metric terms, e.g. at the $(\xi=1)$ element face
\begin{equation}
\tilde{\statevec{F}}^{*,\rm ec} = \avg{Ja_1^1}\statevec{F}_1^{*,\rm ec} + \avg{Ja_2^1}\statevec{F}_2^{*,\rm ec} + \avg{Ja_3^1}\statevec{F}_3^{*,\rm ec}.
\end{equation}
The dissipation terms remain unchanged for the contravariant entropy-stable approximations.

\section{Proofs of Entropy Conservation for Nonlinear Advection Terms}

\subsection{\bf Proof of Entropy Conservation for the One-Dimensional Volume Integral}
\label{app:Proof__volint_entropy_Burgers}
We first show the property \cref{eq:discrete_entropystable_volumeintegral}, rewritten here as
\begin{equation}
\iprod{\mathbb{D}(F)^{\ec},W}_N = \left. F^\ent\right|_{-1}^1.
\label{eq:Bdiscrete_entropystable_volumeintegral}
\end{equation}
For convenience we introduce the entropy potential 
\begin{equation}
\Psi^f_i \equiv W_i\,F_i - F^\ent_i
\label{eqScalarEntropyPotential}
\end{equation}
at each GL node $i=0,...,N$ to rewrite the entropy-conservation condition \cref{eq:ec_condition_numflux} on the two-point flux  
\begin{equation}
F^{\ec}\,\jump{W} - \jump{F\,W}+\jump{F^\ent} = 0
\end{equation}
as 
\begin{equation}
F^{\ec}\,\jump{W} = \jump{\Psi^f}.
\label{eq:FecJumpW}
\end{equation}
We then explicitly write the volume integral in \cref{eq:Bdiscrete_entropystable_volumeintegral} as a sum
\begin{equation}
\iprod{\mathbb{D}(F)^{\ec},W}_N = \sum\limits_{i=0}^N \omega_i W_i\,2\,\sum_{m=0}^N D_{im} F^{\ec}(U_i,U_m).
\end{equation}

Let us now introduce the summation-by-parts matrix $Q = MD$ with entries $Q_{im}\equiv \omega_{i}D_{im}$, which has the property
\begin{equation}
\label{eq:sbp_matrix}
Q+Q^T=B\,,
\end{equation}
where  $B=\mathrm{diag}([-1,0,...,0,1])$ is the boundary evaluation matrix. Because the derivative of a constant is exact, the $Q$ matrix satisfies
\begin{equation}
\sum\limits_{m=0}^N Q_{im} V_i = 0
\label{eq:ConsistencyofQ}
\end{equation}
for all $V_{i}$. 

In terms of $Q$,
\begin{equation}
\begin{split}
\iprod{\mathbb{D}(F)^{\ec},W}_N = \sum\limits_{i=0}^N W_i\,2\,\sum_{m=0}^N Q_{im} F^{\ec}(U_i,U_m)\,.
\end{split}
\end{equation}
But from \cref{eq:sbp_matrix}, $2\,Q_{im} = Q_{im} - Q_{mi} + B_{im}$ so
\begin{equation}
\label{eq:inter_proofa}
\begin{split}
\iprod{\mathbb{D}(F)^{\ec},W}_N &= \sum\limits_{i=0}^N \sum_{m=0}^N W_i\,(Q_{im} - Q_{mi} + B_{im}) F^{\ec}(U_i,U_m)\\&
 = \sum\limits_{i,m} {{W_i}{Q_{im}}{F^{\ec}}\left( {{U_i},{U_m}} \right)}  + \sum\limits_{i,m} {{W_i}{Q_{mi}}{F^{\ec}}\left( {{U_i},{U_m}} \right)}  + \sum\limits_{i,m} {{B_{im}}{W_i}{F^{\ec}}\left( {{U_i},{U_m}} \right)}\,.
\end{split}
\end{equation}
We now re-index the second sum, $i\leftrightarrow m$, use the fact that $F^{\ec}$ is symmetric in its arguments, and recombine the sums to get
\begin{equation}
\label{eq:inter_proof1}
\iprod{\mathbb{D}(F)^{\ec},W}_N = \sum\limits_{i=0}^N \sum_{m=0}^N \left\{(W_i - W_m)\,F^{\ec}(U_i,U_m)\,Q_{im} + B_{im}\,W_i\,F^{\ec}(U_i,U_m)\right\}\,.
\end{equation}
The definition of the entropy potential \cref{eq:FecJumpW} says that we can write the argument
\begin{equation}
(W_i - W_m)\,F^{\ec}(U_i,U_m) = \Psi^f_i - \Psi^f_m\,.
\label{eq:Psidiff1}
\end{equation}
We further note that $B_{im}$ only has entries for $i=m=0$ and $i=m=N$ so with \cref{eqScalarEntropyPotential} and the consistency condition on the entropy conserving flux $W_i\,F^{\ec}(U_i,U_i) = W_i\,F(U_i) = \Psi^f_i+F^\ent_i$,
\begin{equation}
B_{im}\,W_i\,F^{\ec}(U_i,U_m) = B_{im}\,\left(\Psi^f_i+F^\ent_i\right)\,.
\label{eq:BimPsi}
\end{equation}
Inserting \cref{eq:Psidiff1} and \cref{eq:BimPsi} into \cref{eq:inter_proof1}, 
\begin{equation}
\begin{split}
\iprod{\mathbb{D}(F)^{\ec},W}_N &=\sum\limits_{i=0}^N \sum_{m=0}^N  \left( \Psi^f_i - \Psi^f_m\right)\,Q_{im}  + B_{im}\,\left(\Psi^f_i+F^\ent_i\right)\\&
 = \sum\limits_i {\Psi _i^f\sum\limits_m {{Q_{im}}} }  - \sum\limits_{i,m}^{} {\Psi _m^f{Q_{im}}}  + \sum\limits_{i,m} {{B_{im}}\Psi _i^f}  + \sum\limits_{i,m} {{B_{im}}F_i^\ent} \,.
\end{split}
\label{eq:SplitUpSumsAgain}
\end{equation}
By \cref{eq:ConsistencyofQ}, the first term in the second line is zero. Next, with $B_{im}=Q_{im}+Q_{mi}$, re-indexing $i\leftrightarrow m$, and using \cref{eq:ConsistencyofQ},
\begin{equation}
\sum\limits_{i,m} {{B_{im}}\Psi _i^f}  = \sum\limits_{i,m} {{Q_{im}}\Psi _i^f}  + \sum\limits_{i,m} {{Q_{mi}}\Psi _i^f}  = \sum\limits_{i,m} {{Q_{im}}\Psi _m^f}\,.
\end{equation}
Therefore the second and third sums in the second line of \cref{eq:SplitUpSumsAgain} cancel.
What is left is what we set out to show, namely
\begin{equation}
\begin{split}
\iprod{\mathbb{D}(F)^{\ec},W}_N &= \sum\limits_{i=0}^N \sum_{m=0}^N   B_{im}\,F^\ent_i = F^\ent_N - F^\ent_0 =  \left. F^\ent\right|_{-1}^1\,.
\end{split}
\end{equation}
\begin{rem}
The result \cref{eq:Bdiscrete_entropystable_volumeintegral} is a general one that depends only on the summation-by-parts property. It is therefore not specific to the discontinuous Galerkin approximation, {\it per se}. 
 \end{rem}
\subsection{\bf Proof of Entropy Conservation at Interelement Interfaces for Curvilinear Elements}
\label{app:entropy-conservation_3D_curvilinear_surface}

We now show the property \cref{eq:ec_surfint_euler},
\begin{equation}
\sum\limits_\interiorfaces  \int\limits_{N} \left(\left(\contrastatevec{F}_n^{\ec,*}\right)^T\jump{\statevec{W}} - \jump{\left(\contrastatevec{F}_n\right)^T\,\statevec{W}} + \jump{\contravec{F}_n^\ent}\right)\,\dS = 0\,.
\end{equation}

For the approximate 2D surface integral, we evaluate the integrand at $(N+1)^2$ GL quadrature points, see \cref{eq:discrete_surfint}. Therefore we only need to prove that the integrand vanishes discretely at each interior face quadrature point. We assume that the following derivations are restricted to one interior face quadrature point, and skip quadrature point indices.

The general three-dimensional conditions on the \emph{Cartesian components} of the two-point numerical flux are
\begin{equation}
\label{eq:ec_conditions_3D_cart}
\begin{split}
&\left(\statevec F^{\ec,*} \right)^T\,\jump{\statevec W}  - \jump{\statevec F^T\,\statevec W} +\jump{ F^\ent} = 0,\\
&\left(\statevec G^{\ec,*} \right)^T\,\jump{\statevec W}  - \jump{\statevec G^T\,\statevec W} +\jump{ G^\ent} = 0,\\
&\left(\statevec H^{\ec,*} \right)^T\,\jump{\statevec W}  - \jump{\statevec H^T\,\statevec W} +\jump{ H^\ent} = 0,
\end{split}
\end{equation}
where we use the \emph{slave--master} jump definition from \cref{eq:jump_master_slave},
\begin{equation}
\jump{\statevec W}  =  \statevec W_{\mathrm{slave}} - \statevec W_ {\mathrm{master}} ,
\end{equation}
and each interior face has the master element side orientation, so that quadrature points of the slave and the master map on each other.

We make the assumption that the mesh is \emph{watertight}, i.e. that the normal vector and the surface element are continuous across element interfaces. For a conforming hexahedral mesh, the condition holds discretely at the surface quadrature points if we ensure the discrete metric identities \cref{eq:DiscreteMetricIdentities} and if the unit outward facing normal vector and surface element on the element side are constructed from the element metrics by
\begin{equation}
 \label{eq:normvec_surfelem}
 \shat  = \left|\sum\limits_{l=1}^3\left(J\spacevec{a}^{\,l}\right)  \hat{n}^l\right|\,,\quad  \spacevec{n}  =  \frac{1}{\shat } \sum\limits_{l=1}^3\left(J\spacevec{a}^{\,d}\right)  \hat{n}^l.
\end{equation} 
The continuity of the surface metric allows us to use only the metric of the master element side, so that we are able to move the metric into the jump.

We can write any normal contravariant surface flux using the Cartesian fluxes and the metric
\begin{equation}
\begin{split}
\left(\contrastatevec F_n\right)  &= \left(\shat \spacevec{n}\right) \cdot \bigstatevec F  =\shat \left (\statevec F n_1 +\statevec G n_2 + \statevec H n_3\right)  
                                 = \sum\limits_{l=1}^3\hat{n}^l \left(\left(Ja_1^{\,l}\right)   \statevec F +\left(Ja_2^{\,l}\right)   \statevec G + \left(Ja_3^{\,l}\right)   \statevec H \right ) \\ 
                               &= \left(\bigmatrix{M}^T \bigstatevec F \right) \cdot\hat{n}\,.
\end{split}
\end{equation} 

We combine the three Cartesian equations \cref{eq:ec_conditions_3D_cart} with the surface metric $\shat\spacevec{n}$, defined in \cref{eq:normvec_surfelem}, leading to
\begin{equation}
\begin{split}
 \shat \Bigg ( 
 &n_1 \left[(\statevec F^{\ec,*} )^T\,\jump{\statevec W}  - \jump{\statevec F^T\,\statevec W} +\jump{ F^\ent} \right]\\
+&n_2 \left[(\statevec G^{\ec,*} )^T\,\jump{\statevec W}  - \jump{\statevec G^T\,\statevec W} +\jump{ G^\ent} \right]\\
+&n_3 \left[(\statevec H^{\ec,*} )^T\,\jump{\statevec W}  - \jump{\statevec H^T\,\statevec W} +\jump{ H^\ent} \right] \Bigg)= 0\,.
\end{split}
\end{equation}
If we define the contravariant numerical flux as
\begin{equation}
 \contrastatevec{F}^{\ec,*}_{n} \equiv \shat \left(n_1\statevec F^{\ec,*}  + n_2\statevec G^{\ec,*}  + n_3\statevec H^{\ec,*}  \right)\,,
\end{equation} 
then
\begin{equation}
(\contrastatevec F^{\ec,*}_{n})^T\,\jump{\statevec W}   +\left(\shat \spacevec{n}\,\right)  \cdot
\left( - \jump{\bigstatevec F^T\,\statevec W} 
       + \jump{\spacevec F^\ent} \right)  = 0\,.
\end{equation}
Finally, using the continuity of the surface metric, we move it inside the jump terms, yielding
\begin{equation}
\begin{split}
 \left(\statevec F^{\ec,*}_{n}\right)^T\,\jump{\statevec W}  - \jump{\left(\left(\bigmatrix M^T\bigstatevec F\right)\cdot\hat{n}\right)^T\,\statevec W} +\jump{\left(\bigmatrix M^T\spacevec F^\ent\right)\cdot\hat{n}}  &    \\
  =\left(\contrastatevec F^{\ec,*}_{n}\right)^T\,\jump{\statevec W}  - \jump{\left(\contrastatevec F_n\right)^T\,\statevec W} +\jump{\contravec F_n^\ent}   &= 0 \,,\\
\end{split}
\end{equation}
proving that in  \cref{eq:ec_surfint_euler}, the integrand vanishes at each quadrature point of an interior face individually.

\subsection{\bf Proof of Entropy Conservation in 3D Curvilinear Coordinates}
\label{app:Proof__volint_entropy_Euler}

We show in this section that the property \cref{eq:disc_entropy_cons_Euler} holds, reproduced here for convenience as
\begin{equation}
\label{eq:3D_volint_whatwewanttoshow}
\iprod{\spacevec{\mathbb{D}}(\bigcontravec{F})^{\ec},\statevec{W}}_N =  \int\limits_{\partial E,N} \left(\contraspacevec{F}^\ent\cdot\hat{n}\right)\dS\,,
\end{equation}
provided that the discrete metric identities \cref{eq:DiscreteMetricIdentities} are satisfied.

Similar to what was done in \cref{app:Proof__volint_entropy_Burgers}, we introduce three entropy potentials, one for each Cartesian space direction
\begin{equation}
\begin{split}
\Psi^l_{ijk} &\equiv\left(\statevec F^{l}_{ijk}\right)^{T}\statevec W_{ijk} -  F^{l,\ent}_{ijk},\; l = 1,2,3,
\label{eq:PsiComponentDefinition}
\end{split}
\end{equation}
for each GL node $i,j,k=0,...,N$. We use them to rewrite the entropy-conservation condition on the two-point fluxes $F^{l,\ec,*}$, c.f. \cref{eq:ec_conditions_3D_cart}, into
\begin{equation}
\left(\statevec F^{l,\ec}_{(i,m)jk}\right)^T\,\jump{\statevec W}_{(i,m)jk}  = \jump{ \Psi^l}_{(i,m)jk} ,\; l = 1,2,3,
\label{eq:FEC*RelationCartesian}
\end{equation}
for $i,j,k,m=0,...,N$.

To match terms, we expand all the terms of the volume integral approximation
\begin{equation}
\label{eq:entropy-cons_volint2}
\begin{split}
\iprod{\spacevec{\mathbb{D}}(\bigcontravec{F})^{\ec},\statevec{W}}_N\equiv\sum\limits_{i,j,k=0}^N\omega_{ijk}\,\statevec W^T_{ijk}\Bigg [
&\quad 2\sum_{m=0}^N D_{im}\left(\bigstatevec F^{\ec}(U_{ijk}, U_{mjk})\cdot\average{J\spacevec{a}^{\,1}}_{(i,m)jk}\right)\\ 
&+     2\sum_{m=0}^N D_{jm}\left(\bigstatevec F^{\ec}(U_{ijk}, U_{imk})\cdot\average{J\spacevec{a}^{\,2}}_{i(j,m)k}\right)\\ 
&+     2\sum_{m=0}^N D_{km}\left(\bigstatevec F^{\ec}(U_{ijk}, U_{ijm})\cdot\average{J\spacevec{a}^{\,3}}_{ij(k,m)}\right) \Bigg]\,,                                                                                                                                                 
\end{split}
\end{equation}
(c.f. \cref{eq:entropy-cons_volint}) and of the surface integral approximation
\begin{equation}
\begin{split}
\int\limits_{\partial E,N} \left(\contraspacevec{F}^\ent\cdot\hat{n}\right)\,\dS =&\quad
  \sum\limits_{j,k=0}^N\omega_{jk}\left[ \left(\left(J\spacevec{a}^{\,1}\right)_{Njk}\cdot\spacevec F^\ent_{Njk}\right) 
                                        -\left(\left(J\spacevec{a}^{\,1}\right)_{0jk}\cdot\spacevec F^\ent_{0jk}\right) \right]\\
&+\sum\limits_{i,k=0}^N\omega_{ik}\left[ \left(\left(J\spacevec{a}^{\,2}\right)_{iNk}\cdot\spacevec F^\ent_{iNk}\right) 
                                        -\left(\left(J\spacevec{a}^{\,2}\right)_{i0k}\cdot\spacevec F^\ent_{i0k}\right) \right]\\ 
&+\sum\limits_{i,j=0}^N\omega_{ij}\left[ \left(\left(J\spacevec{a}^{\,3}\right)_{ijN}\cdot\spacevec F^\ent_{ijN}\right) 
                                       - \left(\left(J\spacevec{a}^{\,3}\right)_{ij0}\cdot\spacevec F^\ent_{ij0}\right) \right]\,.
\end{split}
\end{equation}
We then focus on the first ($\xi$ direction) term of the volume integral approximation, which allows us to follow the one dimensional proof as closely as possible. The sum can be written in terms of $Q_{im}=w_{i}D_{im}$,
\begin{equation}
\begin{split}
&\sum\limits_{i,j,k=0}^N\omega_{ijk}\,\statevec W^T_{ijk} 2\sum_{m=0}^N D_{im}\,
\bigstatevec F^{\ec}(U_{ijk}, U_{mjk})\cdot \average{J\spacevec{a}^{\,1}}_{(i,m)jk} \\
&\qquad\qquad=\sum\limits_{j,k=0}^N\omega_{jk}\sum\limits_{i=0}^N \statevec W^T_{ijk}2
  \sum_{m=0}^N \omega_iD_{im}\,\bigstatevec F^{\ec}(U_{ijk}, U_{mjk})\cdot\average{J\spacevec{a}^{\,1}}_{(i,m)jk} \\
&\qquad\qquad =\sum\limits_{j,k=0}^N\omega_{jk}\sum\limits_{i=0}^N \statevec W^T_{ijk}\,
 \sum_{m=0}^N 2Q_{im}\,\bigstatevec F^{\ec}(U_{ijk}, U_{mjk})\cdot\average{J\spacevec{a}^{\,1}}_{(i,m)jk}\,.
\end{split}
\label{eq:FirstXiDirectionVolTerm}
\end{equation}
Therefore, we can use the same steps as in one dimension: We use the summation-by-parts property $2\,Q_{im}=Q_{im} - Q_{mi} + B_{im}$, a re-indexing of $i$ and $m$ to subsume the $Q_{mi}$ term, and the facts that $F^{\ec}(U_{ijk}, U_{mjk})$ and the jump operator of the metric term $\average{Ja^1_1}_{(i,m)jk}$ are symmetric with respect to the index $i$ and $m$ to rewrite the $\xi$ direction contribution to the volume integral approximation as
\begin{equation}
\label{eq:inter_proof2}
\begin{split}
\sum\limits_{i=0}^N\statevec W^T_{ijk}\,2\sum_{m=0}^N& Q_{im}\,\bigstatevec F^{\ec}(U_{ijk}, U_{mjk})\cdot\average{J\spacevec{a}^{\,1}}_{(i,m)jk} \\
= \sum\limits_{i=0}^N \sum_{m=0}^N &\statevec W^T_{ijk}\,(Q_{im} - Q_{mi} + B_{im}) 
\bigstatevec F^{\ec}(U_{ijk},U_{mjk})\cdot\average{J\spacevec{a}^{\,1}}_{(i,m)jk}\\
= \sum\limits_{i=0}^N \sum_{m=0}^N &Q_{im}\left(\statevec W_{ijk} - \statevec W_{mjk}\right)^T\,
\bigstatevec F^{\ec}(U_{ijk},U_{mjk})\cdot\average{J\spacevec{a}^{\,1}}_{(i,m)jk}\, \\
&+ B_{im}\,\statevec W^T_{ijk}\,\bigstatevec F^{\ec}(U_{ijk},U_{mjk})\cdot\average{J\spacevec{a}^{\,1}}_{(i,m)jk}.
\end{split}
\end{equation}

Next, we use the relations \cref{eq:FEC*RelationCartesian} of the Cartesian two-point entropy-conserving flux $\statevec F^{l,\ec}$ to note that
\begin{equation}
\label{eq:FRelation1}
 \left(\statevec W_{ijk} - \statevec W_{mjk}\right)^T\,\statevec F^{l,\ec}(U_{ijk},U_{mjk}) =  \Psi^l_{ijk} -  \Psi^l_{mjk}\,, \; l = 1,2,3.
\end{equation}
We further note that $B_{im}$ only has entries for $i=m=0$ and $i=m=N$ (c.f. \cref{eq:BimPsi}) so
\begin{equation}
B_{im}\,\statevec W^T_{ijk}\,\statevec F^{l,\ec}(U_{ijk},U_{mjk}) = B_{im}\,\left( \Psi^l_{ijk}+  F^{l,\ent}_{ijk}\right)\,, \; l = 1,2,3.
\end{equation}
Finally, we can exploit the consistency of the two-point flux and \cref{eq:PsiComponentDefinition} so that for $i=m=0$ and $i=m=N$,
\begin{equation}
 \statevec W^T_{ijk}\,\statevec F^{l,\ec}(U_{ijk},U_{ijk}) = \statevec W^T_{ijk}\,\statevec F(U_{ijk}) =   \Psi^l_{ijk}+ F^{l,\ent}_{ijk}\,, \; l = 1,2,3.
 \label{eq:FRelation3}
\end{equation}

Inserting relations \cref{eq:FRelation1}-\cref{eq:FRelation3} into the last line of \cref{eq:inter_proof2} says that
\begin{equation}
\label{eq:inter_proof3}
\begin{split}
\sum\limits_{i=0}^N &\statevec W^T_{ijk}\,2\sum_{m=0}^N Q_{im}\,\bigstatevec F^{\ec}(U_{ijk}, U_{mjk})\cdot\average{J\spacevec{a}^{\,1}}_{(i,m)jk}\\
 &= \sum\limits_{i=0}^N \sum_{m=0}^N Q_{im}  \left(\spacevec \Psi_{ijk} - \spacevec \Psi_{mjk}\right)\cdot\average{J\spacevec{a}^{\,1}}_{(i,m)jk}  
 +B_{im}\,\left(\spacevec \Psi_{ijk}+\spacevec F^\ent_{ijk}\right)\cdot\average{J\spacevec{a}^{\,1}}_{(i,m)jk}.
\end{split}
\end{equation}
Now, since the derivative of a constant is zero,
\[
\begin{split}
\sum\limits_{i,m = 0}^N {{Q_{im}}{\spacevec \Psi _{ijk}}\cdot\average{J\spacevec{a}^{\,1}}_{(i,m)jk}}  &= \frac{1}{2}\sum\limits_{i = 0}^N {{\spacevec \Psi _{ijk}}\cdot \left(J\spacevec{a}^{\,1}\right)_{ijk}\sum\limits_{m = 0}^N {{Q_{im}}} }  + \frac{1}{2}\sum\limits_{i = 0}^N {\sum\limits_{m = 0}^N {{{\spacevec \Psi _{ijk}}\cdot \left(J\spacevec{a}^{\,1}\right)_{ijk}Q_{im}}} } \\&=  \frac{1}{2}\sum\limits_{i = 0}^N {\sum\limits_{m = 0}^N {{{\spacevec \Psi _{ijk}}\cdot \left(J\spacevec a^{1}\right)_{ijk}Q_{im}}} }.
\end{split}\]
Next, since $B_{im}=0$ unless $i=m=0$ or $i=m=N$,
\[
\sum\limits_{i=0}^N \sum_{m=0}^N B_{im}\spacevec F^\ent_{ijk}\cdot\average{J\spacevec{a}^{\,1}}_{(i,m)jk}=\sum\limits_{i=0}^N \sum_{m=0}^N B_{im}\spacevec F^\ent_{ijk}\cdot \left(J\spacevec{a}^{\,1}\right)_{ijk}.
 \]
Finally, $B_{im}$ being diagonal and symmetric, we can swap the $i\leftrightarrow m$ in $B_{im}\spacevec \Psi_{ijk}$
to rewrite the second line of \cref{eq:inter_proof3}
\begin{equation}
\label{eq:inter_proof4}
\begin{split}
&\sum\limits_{i=0}^N \sum\limits_{m=0}^N   Q_{im}\left(\spacevec \Psi_{ijk} - \spacevec \Psi_{mjk}\right)\cdot\average{J\spacevec{a}^{\,1}}_{(i,m)jk} 
 +  B_{im}\,\left(\spacevec \Psi_{ijk}+\spacevec F^\ent_{ijk}\right)\cdot\average{J\spacevec{a}^{\,1}}_{(i,m)jk} \\
&=\sum\limits_{i=0}^N \sum\limits_{m=0}^N  \frac{1}{2}Q_{im}\spacevec \Psi_{ijk} \cdot\left(J\spacevec{a}^{\,1}\right)_{mjk} +  B_{im}\spacevec F^\ent_{ijk}\cdot\left(J\spacevec{a}^{\,1}\right)_{ijk}  
+ (B_{im}-Q_{im})\spacevec \Psi_{mjk}\cdot\average{J\spacevec{a}^{\,1}}_{(i,m)jk}. \\
\end{split}
\end{equation}
The last step is to use the summation-by-parts property $B_{im}-Q_{im}=Q_{mi}$ and re-index to see that
\begin{equation}
\begin{split}
\sum\limits_{i=0}^N \sum\limits_{m=0}^N(B_{im}-Q_{im})\spacevec \Psi_{mjk}\cdot\average{J\spacevec{a}^{\,1}}_{(i,m)jk} 
&=\sum\limits_{i=0}^N \sum\limits_{m=0}^N  \frac{1}{2}(Q_{mi})\spacevec \Psi_{mjk}\cdot\left(J\spacevec{a}^{\,1}\right)_{ijk} \\
&= \sum\limits_{i=0}^N \sum\limits_{m=0}^N  \frac{1}{2}(Q_{im})\spacevec \Psi_{ijk}\cdot\left(J\spacevec{a}^{\,1}\right)_{mjk}\,.
\end{split}
\end{equation}
Therefore,
\begin{equation}
\label{eq:inter_proof4a}
\begin{split}
&\sum\limits_{i=0}^N \sum\limits_{m=0}^N   Q_{im}\left(\spacevec \Psi_{ijk} - \spacevec \Psi_{mjk}\right)\cdot\average{J\spacevec{a}^{\,1}}_{(i,m)jk} 
 +  B_{im}\,\left(\spacevec \Psi_{ijk}+\spacevec F^\ent_{ijk}\right)\cdot\average{J\spacevec{a}^{\,1}}_{(i,m)jk} 
\\ &=\sum\limits_{i=0}^N \sum\limits_{m=0}^N Q_{im}\spacevec \Psi_{ijk} \cdot\left(J\spacevec{a}^{\,1}\right)_{mjk} +  B_{im}\spacevec F^\ent_{ijk}\cdot\left(J\spacevec{a}^{\,1}\right)_{ijk} .
\end{split}
\end{equation}
To summarize, the first ($\xi$ direction) part of the volume integral approximation \cref{eq:FirstXiDirectionVolTerm} is
\begin{equation}
\begin{split}
\sum\limits_{i,j,k=0}^N&\omega_{ijk}\,\statevec W^T_{ijk} 2\sum_{m=0}^N D_{im}\,
\bigstatevec F^{\ec}(U_{ijk}, U_{mjk})\cdot \average{J\spacevec{a}^{\,1}}_{(i,m)jk}  \\
&= \sum\limits_{j,k=0}^N \omega_{jk}\left(\spacevec F^\ent_{Njk}\cdot\left(J\spacevec{a}^{\,1}\right)_{Njk}-\spacevec F^\ent_{0jk}\cdot\left(J\spacevec{a}^{\,1}\right)_{0jk}+ \sum\limits_{i,m=0}^N\omega_i D_{im}\spacevec \Psi_{ijk} \cdot\left(J\spacevec{a}^{\,1}\right)_{mjk} \right)\,,
\end{split}
\end{equation}
where we have returned $Q_{im}=\omega_i D_{im}$ and represented the boundary terms explicitly.

Similar results hold for the second and third parts of the volume integral approximation, leading to
\begin{equation}
\begin{split}
&\iprod{\spacevec{\mathbb{D}}(\bigcontravec{F})^{\ec},\statevec{W}}_N = \\
&\qquad\quad\sum\limits_{j,k=0}^N \omega_{jk}\left(       \spacevec F^\ent_{Njk} \cdot\left(J\spacevec{a}^{\,1}\right)_{Njk}
                                                         -\spacevec F^\ent_{0jk} \cdot\left(J\spacevec{a}^{\,1}\right)_{0jk}
+           \sum\limits_{i,m=0}^N \omega_i          D_{im}\spacevec\Psi_{ijk} \cdot\left(J\spacevec{a}^{\,1}\right)_{mjk} \right)\\
&\qquad+    \sum\limits_{i,k=0}^N \omega_{ik}\left(       \spacevec F^\ent_{iNk} \cdot\left(J\spacevec{a}^{\,2}\right)_{iNk}
                                                         -\spacevec F^\ent_{i0k} \cdot\left(J\spacevec{a}^{\,2}\right)_{i0k}
+           \sum\limits_{j,m=0}^N \omega_j          D_{jm}\spacevec\Psi_{ijk} \cdot\left(J\spacevec{a}^{\,2}\right)_{imk} \right)\\
&\qquad+    \sum\limits_{i,j=0}^N \omega_{ij}\left(       \spacevec F^\ent_{ijN} \cdot\left(J\spacevec{a}^{\,3}\right)_{ijN}
                                                         -\spacevec F^\ent_{ij0} \cdot\left(J\spacevec{a}^{\,3}\right)_{ij0}
+           \sum\limits_{k,m=0}^N \omega_k          D_{km}\spacevec\Psi_{ijk} \cdot\left(J\spacevec{a}^{\,3}\right)_{ijm} \right)\,.
\end{split}
\end{equation}

Regrouping the boundary terms and the volume terms, we have shown that 
\begin{equation}
\begin{split}
\iprod{\spacevec{\mathbb{D}}(\bigcontravec{F})^{\ec},\statevec{W}}_N=&\int\limits_{\partial E,N} \left(\contraspacevec{F}^\ent\cdot\hat{n}\right)\,\dS\\
&+ \sum\limits_{i,j,k=0}^N\omega_{ijk}\Psi^1_{ijk}\, \left\{\sum_{m=0}^N D_{im}\,\left(Ja^1_1\right)_{mjk}+D_{jm}\,\left(Ja^2_1\right)_{imk}+D_{km}\,\left(Ja^3_1\right)_{ijm}\right\}\\
&+ \sum\limits_{i,j,k=0}^N\omega_{ijk}\Psi^2_{ijk}\, \left\{\sum_{m=0}^N D_{im}\,\left(Ja^1_2\right)_{mjk}+D_{jm}\,\left(Ja^2_2\right)_{imk}+D_{km}\,\left(Ja^3_2\right)_{ijm}\right\}\\
&+ \sum\limits_{i,j,k=0}^N\omega_{ijk}\Psi^3_{ijk}\, \left\{\sum_{m=0}^N D_{im}\,\left(Ja^1_3\right)_{mjk}+D_{jm}\,\left(Ja^2_3\right)_{imk}+D_{km}\,\left(Ja^3_3\right)_{ijm}\right\},                                                                                                                     
\end{split}
\end{equation}
which gives the desired entropy condition \cref{eq:3D_volint_whatwewanttoshow}, provided that the metric identities \cref{eq:DiscreteMetricIdentities} hold discretely, i.e. that
\begin{equation}
\sum_{m=0}^N D_{im}\,\left(J  a^1_n\right)_{mjk}+D_{jm}\,\left(Ja^2_n\right)_{imk}+D_{km}\,\left(Ja^3_n\right)_{ijm}=\left(\sum_{l=1}^{3}\frac{\partial}{\partial\xi^{l}}\mathbb{I}^{N}\left(Ja^l_n\right)\right)_{ijk} = 0 \,,                                                                                                    
\end{equation}
for $n=1,2,3$ and all GL points $i,j,k,=0,...,N$ within an element. For the strategy on how to properly approximate the metric terms so that the metric identities are satisfied, see \cite{Kopriva:2006er}.

\bibliographystyle{plain}

\bibliography{dakBib}

\end{document}